\documentclass{amsbook}
\def\PrintBook{}
\def\UseRussian{}

\def\Defined{}
\scrollmode
\ifx\FilePrefix\undefined
\newcommand{\FilePrefix}{}
\fi
\ifx\UseRussian\Defined
\usepackage{cmap}
\usepackage[T2A,T2B]{fontenc}
\usepackage[cp1251]{inputenc}
\fi
\ifx\GJSFRA\Defined
\usepackage[top=1.905cm,bottom=1.905cm,inner=1.65cm,outer=1.65cm]{geometry}
\fi
\ifx\Presentation\Defined
\usepackage[margin=1cm]{geometry}
\fi
\ifx\CreateSpace\Defined
\usepackage[top=1.905cm,bottom=1.905cm,inner=1.905cm,outer=1.27cm]{geometry}
\def\Publisher{CreateSpace Independent Publishing Platform}
\fi
\ifx\PublishBook\Defined
\def\PrintPaper{}
\usepackage{setspace}
\ifx\UseRussian\undefined
\usepackage{pslatex}
\fi
\usepackage[margin=2cm]{geometry}
\fi
\usepackage{footmisc}
\usepackage[all]{xy}
\usepackage{color}
\ifx\PrintPaper\undefined
\definecolor{CoverColor}{rgb}{.82,.7,.55}
\definecolor{UrlColor}{rgb}{.9,0,.3}
\definecolor{SymbColor}{rgb}{.4,0,.9}
\definecolor{IndexColor}{rgb}{1,.3,.6}
\newcommand\BlueText[1]{\textcolor{blue}{#1}}

\else
\definecolor{UrlColor}{rgb}{.1,.1,.1}
\definecolor{SymbColor}{rgb}{.1,.1,.1}
\definecolor{IndexColor}{rgb}{.1,.1,.1}
\newcommand\BlueText[1]{#1}

\fi

\usepackage{chngcntr}
\usepackage{xr-hyper}
\usepackage[unicode]{hyperref}
\hypersetup{pdfdisplaydoctitle=true}
\hypersetup{colorlinks}
\hypersetup{citecolor=UrlColor}
\hypersetup{urlcolor=UrlColor}
\hypersetup{linkcolor=UrlColor}
\hypersetup{pdffitwindow=true}
\hypersetup{pdfnewwindow=true}
\hypersetup{pdfstartview={FitH}}
\ifx\UseRussian\Defined
\usepackage[english,russian]{babel}
\fi

\newcounter{Index}
\newcounter{Symbol}
\newcounter{Symbols}

\def\hyph{\penalty0\hskip0pt\relax-\penalty0\hskip0pt\relax}
\def\Hyph{-\penalty0\hskip0pt\relax}%

\def\ValueOff{off}%
\def\ValueOn{on}%
\def\Items#1{\ItemList#1,LastItem,}%
\def\LastItem{LastItem}%
\def\ItemList#1,{\def\ViewBook{#1}%
\ifx\ViewBook\LastItem%
\else%
\ifx\ViewBook\BookNumber%
\def\Semafor{on}%
\fi%
\expandafter\ItemList%
\fi%
}%

\newcommand{\ePrints}[1]%
{%
\def\Semafor{off}%
\Items{#1}%
}%

\ifx\PrintPaper\undefined
\newcommand{\gi}[1]{\boldsymbol{\textcolor{IndexColor}{#1}}}
\else
\newcommand{\gi}[1]{\boldsymbol{#1}}
\fi
\newcommand\gii{\gi i}

\newcommand\gij{\gi j}

\makeatletter
\newcommand{\NameDef}[1]{%
\expandafter\gdef\csname #1\endcsname%
}%
\newcommand{\xNameDef}[1]{%
\expandafter\xdef\csname #1\endcsname%
}%
\newcommand{\ShowSymbol}[2]{%
\@nameuse{ViewSymbol#1,,,#2}%
}%
\newcommand{\symb}[3]{%
\@ifundefined{ViewSymbol#2}{%
\addtocounter{Symbols}{1}%
\edef\SymbolId{\arabic{Symbols}}%
\xNameDef{ViewSymbol#2}{\SymbolId}%
\NameDef{ViewSymbol#2:::\SymbolId}{#1}%
\@namedef{RefSymbol}{:}%
}{%
\edef\Symbols{\@nameuse{ViewSymbol#2}}%
\def\aSymbolId{0}%
\@for\Symbol:=\Symbols\do{%
\protected@edef\TempA{#1}%
\protected@edef\TempB{\@nameuse{ViewSymbol#2:::\Symbol}}%
\ifx\TempA\TempB%
\edef\aSymbolId{\Symbol}%
\fi%
}%
\def\Zero{0}%
\ifx\aSymbolId\Zero%
\addtocounter{Symbols}{1}%
\edef\SymbolIds{\@nameuse{ViewSymbol#2},\arabic{Symbols}}%
\xNameDef{ViewSymbol#2}{\SymbolIds}%
\edef\SymbolId{\arabic{Symbols}}%
\NameDef{ViewSymbol#2:::\SymbolId}{#1}%
\else%
\def\SymbolId{\aSymbolId}%
\fi%
\addtocounter{Symbol}{1}%
\@namedef{RefSymbol}{\arabic{Symbol}}%
}%
\@namedef{LabelSymbol}{\label{symbol: #2:\@nameuse{RefSymbol}}}%
\edef\RefIds{RefSymbol#2===\SymbolId}%
\@ifundefined{\RefIds}{%
\xNameDef{\RefIds}{\@nameuse{RefSymbol}}%
}{%
\xNameDef{\RefIds}{\@nameuse{\RefIds},\@nameuse{RefSymbol}}%
}%
\NameDef{ViewSymbol#2,,,#3}{\textcolor{SymbColor}{#1}}%
\def\Temp{#3}%
\def\One{1}%
\def\Two{2}%
\ifx\Temp\One%
$\@nameuse{ViewSymbol#2,,,#3}$%
\fi%
\ifx\Temp\Two%
\[\@nameuse{ViewSymbol#2,,,#3}\]%
\fi%
\@nameuse{LabelSymbol}%
}%
\newcommand{\DefEq}[2]{%
\@ifundefined{ViewEq#2}{%
\NameDef{ViewEq#2}{#1}%
}{%
\errmessage {second entry of DefEq: #2}%
}%
}%
\newcommand{\DefEquation}[2]{%
\DefEq%
{%
\begin{equation}%
#1%
\EqLabel{#2}%
\end{equation}%
}%
{#2}%
}%
\def\ViewParm#1{\protect\getParm#1,endParm,}%
\def\endParm{endParm}%
\def\getParm#1,{\def\temp{#1}%
\ifx\temp\endParm%
\else%
\ShowEq{#1}%
\expandafter\getParm%
\fi%
}%
\newcommand{\EquationParm}[2]{%
\@ifundefined{ViewEq#1[#2]}%
{%
\ViewParm{#2}%
\DefEquation{\ShowEq{#1}}{#1[#2]}%
}{}%
\ShowEq{#1[#2]}%
}%
\newcommand{\DrawEqParm}[3]{%
\ViewParm{#2}%
\@ifundefined{ViewEq#1(#2)}{%
\DefEq%
{%
\ShowEq{#1}%
}{#1(#2)}%
}{%
}%
\DrawEq{#1(#2)}{#3}%
}%
\newcommand\EqLabel[1]{\label{eq: #1}}%
\DefEq%
{%
}%
{=.}%
\DefEq%
{%
}%
{=c}%
\DefEq%
{%
}%
{=z}%
\newcommand\ShowEq[1]{%
\@ifundefined{ViewEq#1}{%
\message {error: missed ShowEq #1}%
  }{%
\@nameuse{ViewEq#1}%
}%
}%
\newcommand\DrawEq[2]{%
\@ifundefined{ViewEq#1}{%
\message {error: missed ShowEq #1}%
  }{%
\def\Temp{}%
\def\Tempa{#2}%
\ifx\Tempa\Temp%
\[%
\@nameuse{ViewEq#1}%
\]%
\else%
\def\Temp{-}%
\ifx\Tempa\Temp%
$\@nameuse{ViewEq#1}$%
\else%
\begin{equation}%
\@nameuse{ViewEq#1}%
\EqLabel{#1(#2)}%
\end{equation}%
\fi%
\fi%
}%
}%
\makeatother

\newcommand{\CRstar}{{}^*{}_*}
\newcommand{\RCstar}{{}_*{}^*}

\newcommand\DcrPartial[1]%
{%
\def\tempa{}%
\def\tempb{#1}%
\ifx\tempa\tempc%
(\partial\CRstar)%
\else%
(\partial_{\gi{#1}}\CRstar)%
\fi%
}%
\newcommand\rcDPartial[1]%
{%
\def\tempa{}%
\def\tempb{#1}%
\ifx\tempa\tempc%
(\RCstar\partial)%
\else%
(\RCstar\partial_{\gi{#1}})%
\fi%
}%
\newcommand\StandPartial[3]%
{%
\frac{\partial^{\gi{#3}} #1}{\partial #2}%
}%

\renewcommand{\uppercasenonmath}[1]{}

\makeatletter
\newcommand\@dotsep{4.5}
\def\@tocline#1#2#3#4#5#6#7
{\relax
\par \addpenalty\@secpenalty\addvspace{#2}%
\begingroup 
\@ifempty{#4}{%
\@tempdima\csname r@tocindent\number#1\endcsname\relax
}{%
\@tempdima#4\relax
}%
\parindent\z@ \leftskip#3\relax \advance\leftskip\@tempdima\relax
\rightskip\@pnumwidth plus1em \parfillskip-\@pnumwidth
#5\leavevmode\hskip-\@tempdima #6\relax
\leaders\hbox{$\m@th
\mkern \@dotsep mu\hbox{.}\mkern \@dotsep mu$}\hfill
\hbox to\@pnumwidth{\@tocpagenum{#7}}\par
\nobreak
\endgroup
}
\makeatother 

\ifx\PrintBook\undefined

\makeatletter
\renewcommand{\@indextitlestyle}{%
\twocolumn[\section{\indexname}]%
\def\IndexSpace{off}%
}
\makeatother 
\ifx\PrintPaper\undefined
\thanks{\href{mailto:Aleks\_Kleyn@MailAPS.org}{Aleks\_Kleyn@MailAPS.org}}
\ePrints{1102.1776,1201.4158}
\ifx\Semafor\ValueOff
\thanks{\ \ \ \url{http://AleksKleyn.dyndns-home.com:4080/}\ \ \ \ \ \url{http://arxiv.org/a/kleyn\_a\_1}}
\thanks{\ \ \ \url{http://sites.google.com/site/AleksKleyn/}\ \ \ \ \url{http://AleksKleyn.blogspot.com/}}
\fi
\fi
\else

\makeatletter
\def\@maketitle{%
  \cleardoublepage \thispagestyle{empty}%
  \begingroup \topskip\z@skip
  \null\vfil
  \begingroup
  \LARGE\bfseries \centering
  \openup\medskipamount
  \@title
  \par
  \ifx\subtitle\undefined
  \else
  \centerline{\ }
  \centerline{\emph\subtitle}
  \fi
  \ifx\subtitleA\undefined
  \else
  \centerline{\emph\subtitleA}
  \fi
  \ifx\edition\undefined
  \else
  \centerline{\emph\edition}
  \fi
  \par\vspace{24pt}%
  \def\and{\par\medskip}\centering
  \mdseries\authors\par\bigskip
  \endgroup
  \vfill
\noindent
\href{mailto:Aleks\_Kleyn@MailAPS.org}{Aleks\_Kleyn@MailAPS.org}
\newline
\url{http://AleksKleyn.dyndns-home.com:4080/}
\newline
\url{http://sites.google.com/site/AleksKleyn/}
\newline
\url{http://arxiv.org/a/kleyn\_a\_1}
\newline
\url{http://AleksKleyn.blogspot.com/}
  \newpage\thispagestyle{empty}
  \begin{center}
    \ifx\@empty\@subjclass\else\@setsubjclass\fi
    \ifx\@empty\@keywords\else\@setkeywords\fi
    \ifx\@empty\@translators\else\vfil\@settranslators\fi
    \ifx\@empty\thankses\else\vfil\@setthanks\fi
  \end{center}
  \vfil
  \@setabstract
\vfil
  \def\Temp{0000}
  \ifx\copyrightyear\Temp
  \else
  \begin{center}
\begin{tabular}{@{}c}
Copyright\ \copyright\ \copyrightyear\ \copyrightholder
\\
All rights reserved.
\end{tabular}
  \end{center}
  \fi
  \ifx\Publisher\undefined%
  \else
  \begin{center}
\begin{tabular}{@{}c}
\Publisher
\end{tabular}
  \end{center}
  \fi
  \ifx\ISBN\undefined%
  \else%
 \begin{center}
\begin{tabular}{@{}r@{\ }l}
ISBN:&\ISBN
\\
ISBN-13:&\ISBNa
\end{tabular}
  \end{center}
  \fi%
  \ifx\titleNote\undefined
  \else
  \par\vspace{24pt}%
  \centerline{\mdseries\titleNote}
	  \centerline{\Title}
	  \ifx\Subtitle\undefined
	  \else
	  \centerline{\emph\Subtitle}
	  \fi
	  \ifx\Edition\undefined
	  \else
	  \centerline{\Edition}
	  \fi
	  \centerline{\Authors}
  \fi
  \endgroup}
\renewcommand{\@indextitlestyle}{%
\twocolumn[\chapter{\indexname}]%
\def\IndexSpace{off}%
\let\@secnumber\@empty
\chaptermark{\indexname}%
}
\makeatother 
\email{\href{mailto:Aleks\_Kleyn@MailAPS.org}{Aleks\_Kleyn@MailAPS.org}}
\ePrints{1102.1776,1201.4158}
\ifx\Semafor\ValueOff
\urladdr{\url{http://AleksKleyn.dyndns-home.com:4080/}}
\urladdr{\url{http://sites.google.com/site/alekskleyn/}}
\urladdr{\url{http://arxiv.org/a/kleyn\_a\_1}}
\urladdr{\url{http://AleksKleyn.blogspot.com/}}
\fi
\fi

\ifx\SelectlEnglish\undefined
\ifx\UseRussian\undefined
\def\SelectlEnglish{}
\fi
\fi

\newcommand\arXivOldRef{http://arxiv.org/PS_cache/}
\newcommand\arXivRef{http://arxiv.org/pdf/}
\newcommand\AmazonRef{http://www.amazon.com/s/ref=nb_sb_noss?url=search-alias=aps&field-keywords=aleks+kleyn}
\newcommand\wRefDef[2]
{
\def\Tempa{#1}
\def\Tempb{0405.027}
\ifx\Tempa\Tempb
\def\wRef{\arXivOldRef gr-qc/pdf/0405/0405027v3.pdf}
\fi
\def\Tempb{0405.028}
\ifx\Tempa\Tempb
\def\wRef{\arXivOldRef gr-qc/pdf/0405/0405028v5.pdf}
\fi
\def\Tempb{0412.391}
\ifx\Tempa\Tempb
\def\wRef{\arXivOldRef math/pdf/0412/0412391v4.pdf}
\fi
\def\Tempb{0612.111}
\ifx\Tempa\Tempb
\def\wRef{\arXivOldRef math/pdf/0612/0612111v2.pdf}
\fi
\def\Tempb{0701.238}
\ifx\Tempa\Tempb
\def\wRef{\arXivOldRef math/pdf/0701/0701238v6.pdf}
\fi
\def\Tempb{0702.561}
\ifx\Tempa\Tempb
\def\wRef{\arXivOldRef math/pdf/0702/0702561v3.pdf}
\fi
\def\Tempb{0707.2246}
\ifx\Tempa\Tempb
\def\wRef{\arXivRef 0707.2246v2.pdf}
\fi
\def\Tempb{0803.3276}
\ifx\Tempa\Tempb
\def\wRef{\arXivRef 0803.3276v3.pdf}
\fi
\def\Tempb{0812.4763}
\ifx\Tempa\Tempb
\def\wRef{\arXivRef 0812.4763v7.pdf}
\fi
\def\Tempb{0906.0135}
\ifx\Tempa\Tempb
\def\wRef{\arXivRef 0906.0135v3.pdf}
\fi
 \def\Tempb{0909.0855}
\ifx\Tempa\Tempb
\def\wRef{\arXivRef 0909.0855v5.pdf}
\fi
 \def\Tempb{0912.3315}
\ifx\Tempa\Tempb
\def\wRef{\arXivRef 0912.3315v3.pdf}
\fi
 \def\Tempb{0912.4061}
\ifx\Tempa\Tempb
\def\wRef{\arXivRef 0912.4061v2.pdf}
\fi
 \def\Tempb{1001.4852}
\ifx\Tempa\Tempb
\def\wRef{\arXivRef 1001.4852.pdf}
\fi
 \def\Tempb{1003.3714}
\ifx\Tempa\Tempb
\def\wRef{\arXivRef 1003.3714v2.pdf}
\fi
 \def\Tempb{1003.1544}
\ifx\Tempa\Tempb
\def\wRef{\arXivRef 1003.1544v2.pdf}
\fi
 \def\Tempb{1006.2597}
\ifx\Tempa\Tempb
\def\wRef{\arXivRef 1006.2597v2.pdf}
\fi
 \def\Tempb{1011.3102}
\ifx\Tempa\Tempb
\def\wRef{\arXivRef 1011.3102.pdf}
\fi
 \def\Tempb{1104.5197}
\ifx\Tempa\Tempb
\def\wRef{\arXivRef 1104.5197.pdf}
\fi
 \def\Tempb{1105.4307}
\ifx\Tempa\Tempb
\def\wRef{\arXivRef 1105.4307.pdf}
\fi
 \def\Tempb{1107.1139}
\ifx\Tempa\Tempb
\def\wRef{\arXivRef 1107.1139.pdf}
\fi
 \def\Tempb{1107.5037}
\ifx\Tempa\Tempb
\def\wRef{\arXivRef 1107.5037.pdf}
\fi
 \def\Tempb{1111.6035}
\ifx\Tempa\Tempb
\def\wRef{\arXivRef 1111.6035.pdf}
\fi
 \def\Tempb{1202.6021}
\ifx\Tempa\Tempb
\def\wRef{\arXivRef 1202.6021.pdf}
\fi
 \def\Tempb{1211.6965}
\ifx\Tempa\Tempb
\def\wRef{\arXivRef 1211.6965.pdf}
\fi
 \def\Tempb{1302.7204}
\ifx\Tempa\Tempb
\def\wRef{\arXivRef 1302.7204v1.pdf}
\fi
 \def\Tempb{1305.4547}
\ifx\Tempa\Tempb
\def\wRef{\arXivRef 1305.4547.pdf}
\fi
 \def\Tempb{1310.5591}
\ifx\Tempa\Tempb
\def\wRef{\arXivRef 1310.5591.pdf}
\fi
 \def\Tempb{1502.04063}
\ifx\Tempa\Tempb
\def\wRef{\arXivRef 1502.04063v2.pdf}
\fi
 \def\Tempb{1505.03625}
\ifx\Tempa\Tempb
\def\wRef{\arXivRef 1505.03625v1.pdf}
\fi
 \def\Tempb{MRepro}
\ifx\Tempa\Tempb
\def\wRef{\arXivRef MReprov1.pdf}
\fi
 \def\Tempb{8433-5163}
\ifx\Tempa\Tempb
\def\wRef{\AmazonRef}
\fi
 \def\Tempb{8443-0072}
\ifx\Tempa\Tempb
\def\wRef{\AmazonRef}
\fi
 \def\Tempb{4776-3181}
\ifx\Tempa\Tempb
\def\wRef{\AmazonRef}
\fi
 \def\Tempb{5059-9176}
\ifx\Tempa\Tempb
\def\wRef{\AmazonRef}
\fi
 \def\Tempb{5114-6019}
\ifx\Tempa\Tempb
\def\wRef{\AmazonRef}
\fi
 \def\Tempb{5148-4632}
\ifx\Tempa\Tempb
\def\wRef{\AmazonRef}
\fi
 \def\Tempb{BRepro}
\ifx\Tempa\Tempb
\def\wRef{\AmazonRef}
\fi
 \def\Tempb{CACAA.01.291}
\ifx\Tempa\Tempb
\def\wRef{http://www.cliffordanalysis.com/}
\fi
 \def\Tempb{CACAA.02.97}
\ifx\Tempa\Tempb
\def\wRef{http://www.cliffordanalysis.com/}
\fi
 \def\Tempb{GJSFRA.13.1.39}
\ifx\Tempa\Tempb
\def\wRef{http://www.cliffordanalysis.com/}
\fi
\externaldocument[#1-#2-]{\FilePrefix #1.#2}[\wRef]
}
\newcommand\LanguagePrefix{}%
\makeatletter

\newcommand\StopLabelItem
{
\def\theenumi{\@arabic\c@enumi}
\def\labelenumi{(\theenumi)}
}
\makeatother
\ifx\SelectlEnglish\undefined
\newcommand\CurrentLanguage{Russian.}%
\author{Александр Клейн}
\newtheorem{theorem}{Теорема}[section]
\newtheorem{corollary}[theorem]{Следствие}
\newtheorem{convention}[theorem]{Соглашение}
\newtheorem{definition}[theorem]{Определение}
\newtheorem{example}[theorem]{Пример}
\newtheorem{xca}[theorem]{Exercise}
\newtheorem{remark}[theorem]{Замечание}
\newtheorem{question}[theorem]{Вопрос}
\newtheorem{lemma}[theorem]{Лемма}
\newtheorem{Example}[theorem]{Example}

\makeatletter
\newcommand\xRef[2]%
{%
\@ifundefined{xRefDef#1}{%
\ref{#2}%
}{%
\citeBib{#1}-\ref{#1-Russian-#2}%
}%
}%
\newcommand\xEqRef[2]%
{%
\@ifundefined{xRefDef#1}{%
\eqref{eq: #2}%
}{%
\citeBib{#1}-\eqref{#1-Russian-eq: #2}%
}%
}%
\makeatother
\ifx\PrintBook\undefined
\newcommand{\BibTitle}{%
\section{Список литературы}%
}
\else
\newcommand{\BibTitle}{%
\chapter{Список литературы}%
}
\fi
\else
\newcommand\CurrentLanguage{English.}%
\author{Aleks Kleyn}
\newtheorem{theorem}{Theorem}[section]

\newtheorem{example}[theorem]{Example}

\makeatletter
\newcommand\xRef[2]%
{%
\@ifundefined{xRefDef#1}{%
\ref{#2}%
}{%
\citeBib{#1}-\ref{#1-English-#2}%
}%
}%
\newcommand\xEqRef[2]%
{%
\@ifundefined{xRefDef#1}{%
\eqref{eq: #2}%
}{%
\citeBib{#1}-\eqref{#1-English-eq: #2}%
}%
}%
\makeatother
\ifx\PrintBook\undefined
\newcommand{\BibTitle}{%
\section{References}%
}
\else
\newcommand{\BibTitle}{%
\chapter{References}%
}
\fi
\fi
\theoremstyle{definition}
\theoremstyle{remark}

\ifx\PrintBook\undefined
\numberwithin{Hfootnote}{section}
\else
\numberwithin{section}{chapter}
\numberwithin{footnote}{chapter}
\numberwithin{Hfootnote}{chapter}
\fi

\ifx\Presentation\undefined
\numberwithin{equation}{section}
\numberwithin{figure}{section}
\numberwithin{table}{section}
\numberwithin{Item}{section}
\fi

\makeatletter
\newcommand\org@maketitle{}
\let\org@maketitle\maketitle
\def\maketitle{%
\hypersetup{pdftitle={\@title}}%
\hypersetup{pdfauthor={\authors}}%
\hypersetup{pdfsubject=\@keywords}%
\ifx\UseRussian\Defined
\pdfbookmark[1]{\@title}{TitleRussian}
\else
\pdfbookmark[1]{\@title}{TitleEnglish}
\fi
\org@maketitle
}
\def\make@stripped@name#1{%
\begingroup
\escapechar\m@ne
\global\let\newname\@empty
\protected@edef\Hy@tempa{\CurrentLanguage #1}%
\edef\@tempb{%
\noexpand\@tfor\noexpand\Hy@tempa:=%
\expandafter\strip@prefix\meaning\Hy@tempa
}%
\@tempb\do{%
\if\Hy@tempa\else
\if\Hy@tempa\else
\xdef\newname{\newname\Hy@tempa}%
\fi
\fi
}%
\endgroup
}%
\newenvironment{enumBib}{%
\BibTitle
\advance\@enumdepth \@ne
\edef\@enumctr{enum\romannumeral\the\@enumdepth}\list
{\csname biblabel\@enumctr\endcsname}{\usecounter
{\@enumctr}\def\makelabel##1{\hss\llap{\upshape##1}}}
}{%
\endlist
}

\makeatletter

\newcommand{\BiblioItem}[2]
{
\def\Semafor{off}
\@ifundefined{\LanguagePrefix ViewCite#1}{}{%
\def\Semafor{on}%
}%
\ifx\Semafor\ValueOff
\@ifundefined{xRefDef#1}{}{%
\def\Semafor{on}%
}%
\fi
\ifx\Semafor\ValueOn
\ifx\IndexState\ValueOff
\begin{enumBib}
\def\IndexState{on}
\fi
\item \label{\LanguagePrefix bibitem: #1}#2%
\fi
}
\makeatother

\newcommand{\CloseBiblio}
{
\ifx\IndexState\ValueOn
\end{enumBib}
\def\IndexState{off}
\fi
}

\makeatletter
\def\StartCite{[}%
\def\citeBib#1{\protect\showCiteBib#1,endCite,}%
\def\endCite{endCite}%
\def\showCiteBib#1,{\def\temp{#1}%
\ifx\temp\endCite
]%
\def\StartCite{[}%
\else
\StartCite\LanguagePrefix \ref{\LanguagePrefix bibitem: #1}%
\@ifundefined{\LanguagePrefix ViewCite#1}{%
\NameDef{\LanguagePrefix ViewCite#1}{}%
}{%
}%
\def\StartCite{, }%
\expandafter\showCiteBib%
\fi}%
\makeatother

\newcommand{\arp}{\ar @{-->}}
\newcommand{\bundle}[4]%
{%
\def\tempa{}%
\def\tempb{#3}%
\def\tempc{#1}%
\ifx\tempa\tempb%
\ifx\tempa\tempc%
#2%
\else%
\xymatrix{#2:#1\arp[r]&#4}%
\fi%
\else%
\ifx\tempa\tempc%
#2[#3]%
\else%
\xymatrix{#2[#3]:#1\arp[r]&#4}%
\fi%
\fi%
}%
\makeatletter
\newcommand{\AddIndex}[2]%
{%
\@ifundefined{RefIndex#2}{%
\xNameDef{RefIndex#2}{:}%
\@namedef{LabelIndex}{\label{index: #2::}}%
}{%
\addtocounter{Index}{1}%
\xNameDef{RefIndex#2}{\@nameuse{RefIndex#2},\arabic{Index}}%
\@namedef{LabelIndex}{\label{index: #2:\arabic{Index}}}%
}%
\@nameuse{LabelIndex}%
{\bf #1}%
}%
\newcommand{\Index}[2]%
{%
\def\Semafor{off}%
\@ifundefined{RefIndex#2}{%
}{%
\def\Semafor{on}
}%
\ifx\Semafor\ValueOn%
\def\tempa{}%
\def\tempb{#2}%
\ifx\IndexState\ValueOff%
\begin{theindex}%
\def\IndexState{on}%
\fi%
\ifx\IndexSpace\ValueOn%
\indexspace%
\def\IndexSpace{off}%
\fi%
\item #1%
\ifx\tempa\tempb%
\else%
\edef\PageRefs{\@nameuse{RefIndex#2}}
\def\Sep{\ }%
\@for\PageRef:=\PageRefs\do{%
\Sep
\pageref{index: #2:\PageRef}%
\def\Sep{,\ }%
}%
\fi%
\fi%
}%

\newcommand{\Symb}[2]
{
\def\Semafor{off}
\@ifundefined{ViewSymbol#2}{%
}{%
\def\Semafor{on}
}%
\ifx\Semafor\ValueOn
\ifx\IndexState\ValueOff
\begin{theindex}
\def\IndexState{on}
\fi
\ifx\IndexSpace\ValueOn
\indexspace
\def\IndexSpace{off}
\fi
\edef\Symbols{\@nameuse{ViewSymbol#2}}
\@for\Symbol:=\Symbols\do{%
\edef\Temp{ViewSymbol#2:::\Symbol}
\item $\displaystyle\textcolor{SymbColor}{\@nameuse{\Temp}}$
\ \ #1
\edef\PageRefs{\@nameuse{RefSymbol#2===\Symbol}}
\def\Sep{}%
\@for\PageRef:=\PageRefs\do{%
\Sep
\pageref{symbol: #2:\PageRef}%
\def\Sep{,\ }%
}%
}%
\fi
}

\newcommand{\Symba}[2]
{
\def\Semafor{off}
\@ifundefined{ViewSymbol#2}{%
}{%
\def\Semafor{on}
}%
\ifx\Semafor\ValueOn
\ifx\IndexState\ValueOff
\begin{theindex}
\def\IndexState{on}
\fi
\ifx\IndexSpace\ValueOn
\indexspace
\def\IndexSpace{off}
\fi
\item $\displaystyle\@nameuse{ViewSymbol#2}$\ \ #1
\edef\PageRefs{\@nameuse{RefSymbol#2}}
\def\Sep{}%
\@for\PageRef:=\PageRefs\do{%
\Sep
\pageref{symbol: #2:\PageRef}%
\def\Sep{,\ }%
}%
\fi
}

\makeatother

\newcommand{\SetIndexSpace}%
{%
\def\IndexSpace{on}%
}%

\newcommand{\CloseIndex}
{
\ifx\IndexState\ValueOn
\end{theindex}
\def\IndexState{off}
\fi
}

\def\LastMemo{LastMemo}%
\def\MemoList#1//{\def\temp{#1}%
\ifx\temp\LastMemo
\else%
\setlength{\parindent}{5mm}
\par
\BlueText{#1}%
\expandafter\MemoList%
\fi%
}     

\hypersetup{pdfsubject={English Russian Science Dictionary}}
\newcommand{\ERitem}[2]
	{\item[#1]#2}
\newcommand{\REitem}[2]
	{\item[#2]#1}

\newcommand{\ERexample}[6]
	{
		\def\tempa{}%
		\begin{Example}
		{\ }
		\par {\bf #1}
		\def\tempb{#3}%
		\ifx\tempa\tempb%
		\else%
			\begin{flushright}
			\renewcommand\LanguagePrefix{English.}%
			\begin{tabular}{r}
			\hline
			\makebox[200pt][r]{see \citeBib{#3}, p. #4}%
			\end{tabular}
			\end{flushright}
		\fi%
		\par #2
		\def\tempb{#5}%
		\ifx\tempa\tempb%
		\else%
			\begin{flushright}
			\renewcommand\LanguagePrefix{Russian.}%
			\begin{tabular}{r}
			\hline
			\makebox[200pt][r]{см. \citeBib{#5}, стр. #6}%
			\end{tabular}
			\end{flushright}
		\fi%
		\qed
		\end{Example}
	}
\newcommand{\REexample}[6]
	{
		\def\tempa{}%
		\begin{example}
		{\ }
		\par {\bf #2}
		\def\tempb{#5}%
		\ifx\tempa\tempb%
		\else%
			\begin{flushright}
			\renewcommand\LanguagePrefix{Russian.}%
			\begin{tabular}{r}
			\hline
			\makebox[200pt][r]{см. \citeBib{#5}, стр. #6}%
			\end{tabular}
			\end{flushright}
		\fi%
		\par #1
		\def\tempb{#3}%
		\ifx\tempa\tempb%
		\else%
			\begin{flushright}
			\renewcommand\LanguagePrefix{English.}%
			\begin{tabular}{r}
			\hline
			\makebox[200pt][r]{see \citeBib{#3}, p. #4}%
			\end{tabular}
			\end{flushright}
		\fi%
		\qed
		\end{example}
	}
\newcommand{\RName}[4]
	{\item \BlueText{\ifx#3\@empty\else#3 \fi{\bf #4}}
	\ifx#1\@empty\else#1 \fi#2}
\newcommand{\EName}[4]
	{\item \BlueText{\ifx#1\@empty\else#1 \fi{\bf #2}}
	\ifx#3\@empty\else#3 \fi#4}

\makeatletter
\def\@maketitle{%
  \cleardoublepage \thispagestyle{empty}%
  \begingroup \topskip\z@skip
  \null\vfil
  \begingroup
  \LARGE\bfseries \centering
  \openup\medskipamount
  English Russian Scientific Dictionary
  \par
  \ifx\edition\undefined
  \else
  \centerline{\emph\edition}
  \fi
  \par\vspace{24pt}%
 \centering
  \mdseries Aleks Kleyn
  \par\vspace{72pt}%
  {\bf Русско-английский научный словарь}
  \par
  \ifx\edition\undefined
  \else
  \centerline{\emph\edition}
  \fi
  \par\vspace{24pt}%
 \centering
  \mdseries Александр Клейн
  \par
  \bigskip
  \endgroup
  \vfill
\noindent
\href{mailto:Aleks\_Kleyn@MailAPS.org}{Aleks\_Kleyn@MailAPS.org}
\newline
\url{http://AleksKleyn.dyndns-home.com:4080/}
\newline
\url{http://sites.google.com/site/AleksKleyn/}
\newline
\url{http://arxiv.org/a/kleyn\_a\_1}
\newline
\url{http://AleksKleyn.blogspot.com/}
  \newpage\thispagestyle{empty}
  \begin{center}
    \ifx\@empty\@subjclass\else\@setsubjclass\fi
    \ifx\@empty\@keywords\else\@setkeywords\fi
    \ifx\@empty\@translators\else\vfil\@settranslators\fi
    \ifx\@empty\thankses\else\vfil\@setthanks\fi
  \end{center}
  \vfil
  \@setabstract
\vfil
  \def\Temp{0000}
  \ifx\copyrightyear\Temp
  \else
  \begin{center}
\begin{tabular}{@{}c}
Copyright\ \copyright\ \copyrightyear\ \copyrightholder
\\
All rights reserved.
\end{tabular}
  \end{center}
  \fi
  \ifx\Publisher\undefined%
  \else
  \begin{center}
\begin{tabular}{@{}c}
\Publisher
\end{tabular}
  \end{center}
  \fi
  \ifx\ISBN\undefined%
  \else%
 \begin{center}
\begin{tabular}{@{}r@{\ }l}
ISBN:&\ISBN
\\
ISBN-13:&\ISBNa
\end{tabular}
  \end{center}
  \fi%
  \ifx\titleNote\undefined
  \else
  \par\vspace{24pt}%
  \centerline{\mdseries\titleNote}
	  \centerline{\Title}
	  \ifx\Subtitle\undefined
	  \else
	  \centerline{\emph\Subtitle}
	  \fi
	  \ifx\Edition\undefined
	  \else
	  \centerline{\Edition}
	  \fi
	  \centerline{\Authors}
  \fi
  \endgroup}
\makeatother

\begin{document}
\title{English Russian Scientific Dictionary}
\begin{abstract}
English Russian and Russian English dictionaries presented
in this book are dedicated to help translate a scientific text
from one language to another. I also included the bilingual
name index into this book.

Задача англо-русского и русско-английского словарей, представленных
в этой книге, - это помощь в переводе научного текста. Я включил в эту книгу также двуязычный
именной указатель.
\end{abstract}

\maketitle
\tableofcontents

\shorthandoff{"}

				\chapter{Preface}

According to legend, all people spoke one language.
Then they decided to build up the tower of Babel
in order to  climb up to the sky. To stop project,
L-rd confused their languages...
Much water has flowed since then. A lot of different occupations appeared.
Occupation of interpreter is one of the most ancient and important.
At the same time translation of text
is art form.\footnote{
Although we can understand one another,
in spite of difference of languages,
people are not ready
to build up new tower of Babel. I believe people should learn to understand
each other, no matter what views they may have. Mankind must become mature
to get this goal.}

Lively text contains phraseology, slang.
It is hard to translate it one to one.
Translation of sonnets by Shakespeare from English into Russian
is very striking example. You can compare translations of sonnets by Shakespeare
made by S. J. Marshak and A. M. Finkel.
They are different works of art although original is common.
Although translations are different, each of them leaves
indelible track. Each of translations reflects individual perception
of interpreter, his emotional experience. The strength of translations is that
Marshak and Finkel are co-authors of Shakespeare.

There exists opinion that it is much easier to translate technical text.
Unambiguity of translation of  term is one of the reasons.
Another reason is existence of set of phrases; using of
such phrases make easier realization of the text.

However I do not agree with this point of view.
It is not enough to write good theorems for writing of the paper.
We expect that somebody reads the paper. And good style is important here.
This is more important when we translate a paper from one language to other.
Sometimes it is not easy task to put text into frame prepared
in advance.
Moreover, any text has emotional color which expresses
author's relationship.

I am looking for outside assistance during translation
when I meet either unusual
grammatical form, or expression which I cannot translate
unambiguously.
Today we can find such help because of development
of software and internet.
For instance, you can find a lot of forums
on web page
\begin{center}
\href{http://forum.lingvo.ru/actualforum.aspx}
{http://forum.lingvo.ru/actualforum.aspx}
\end{center}
and somebody in these forums can help you to translate text.
I frequently use help of the forum
"English-Russian-English translation"
\begin{center}
\href{http://forum.lingvo.ru/actualtopics.aspx?bid=18}
{http://forum.lingvo.ru/actualtopics.aspx?bid=18}
\end{center}
and I am very grateful to everybody who helps
me to translate my papers.

In Russian literature there is good custom to mention name of
interpreter. Because I turn for help to forum, it is hard for me
to follow this custom. Therefore I decided to
tell about people who help me in the book dedicated to
process of translation.

When I feel that I need additional assistance,
I'm looking for patern in books and articles
where I can find appropriate terminology or
sample how to write proper sentence.
In particular, I recommend book search in Google
\begin{center}
\href{http://books.google.com/bkshp?tab=wp}
{http://books.google.com/bkshp?tab=wp}
\end{center}

Very often when I translate I feel that I need
to change original sentence in order to make it more sufficient to my ideas
in both languages.
This is why I consider translation as
an inseparable part of work on the text
and creative process.

Initially I started this vocabulary as immediate helper to write papers.
Although there are a lot of scientific texts in internet, search for a particular term or definition
takes a lot of time.
And I decided to save translation of terms and definitions which
was important for me and I met in books published on both languages.
So the vocabulary has appeared.
I wrote vocabulary slowly.
When vocabulary became large enough I decided to make it available for others.
I also took into account the wishes expressed in response to
previous versions of vocabulary.
This vocabulary involves only physical and mathematical
terminology. Word may have an alternative translation or different meaning
in common speech.
However I did not include this information into vocabulary

I also included the bilingual
name index into this book.
I added into this list only names which I met in both Russian and English texts.
However there is one more problem
which I meet when I translate papers.
Sometimes the relationship between the writing and
and the pronunciation is not clear. It would be good idea to add transcription,
however I did not make it because I did not have reliable source.

				\chapter{Введение}

Согласно преданию, все люди разговаривали на одном языке.
Тогда они решили построить Вавилонскую башню,
чтобы взобраться по ней на небо. Чтобы остановить строительство,
Б-г смешал их языки...
Немало воды утекло с тех пор. Появилось много разных профессий.
Профессия переводчика одна из древних и важных.
И в тоже время перевод текста - это
искусство.\footnote{
Хотя мы можем понимать друг друга,
несмотря на различие языков,
люди не готовы
строить новую Вавилонскую башню. Я думаю, люди
должны научиться понимать друг друга независимо от того,
каких взглядов они придерживаются. А для этого человечество
ещё должно повзрослеть.}

Живой текст содержит фразеологизмы, сленг.
Всё то, что дословно не переведёшь. Очень яркий пример для меня
- это перевод сонетов Шекспира
на русский. Сравните переводы сонетов Шекспира,
сделанные С. Я. Маршаком и А. М. Финкелем.
Это разные произведения, хотя оригинал общий.
И хотя переводы различны, каждый из них оставляет
неизгладимый след. Каждый из переводов отражает личное восприятие
переводчика, его переживание. Сила переводов в том, что и
Маршак, и Финкель оказались соавторами Шекспира.

Существует представление, что техническую литературу переводить легче.
Одна из причин - однозначность перевода технических терминов.
Другая причина - существование некоторого набора штампов, использование
которых облегчает понимание текста.

Однако я не согласен с этой точкой зрения.
Недостаточно вывести хорошие теоремы, для того, чтобы написать статью.
Статья должна быть прочитана. И хороший стиль играет здесь не последнюю роль.
Особенно это важно, когда мы переводим статью на другой язык.
Поместить фразу в заранее приготовленный штамп порой нелегко.
Кроме того, любой текст имеет эмоциональную окраску, выражающую
отношение автора.

Я ищу постороннюю помощь при переводе
если я столкнулся либо с необычной
грамматической формой, либо с выражением, которое не переводится
однозначно.
Сегодня, благодаря развитию программных средств и интернета,
найти такую помощь легче, чем это было недавно.
Например, на интернет странице
\begin{center}
\href{http://forum.lingvo.ru/actualforum.aspx}
{http://forum.lingvo.ru/actualforum.aspx}
\end{center}
вы можете найти разнообразные форумы, где вам всегда
помогут перевести текст. Я очень часто пользуюсь
помощью форума "Англо-русско-английский перевод"
\begin{center}
\href{http://forum.lingvo.ru/actualtopics.aspx?bid=18}
{http://forum.lingvo.ru/actualtopics.aspx?bid=18}
\end{center}
и я приношу свою глубокую благодарность всем, кто помогает
мне в переводе моих статей.

В русской литературе существует хорошая традиция указывать имя
переводчика. Так как я обращаюсь за помощью к форуму, мне трудно
следовать этой традиции. Поэтому я принял решение
рассказать о тех, кто мне помогает, в книге, специально посвящённой
процессу перевода.

В тех случаях, когда я чувствую, что этой помощи мне не достаточно,
я ищу образцы текста в книгах и статьях,
где я могу найти правильную терминологию или образец употребления
той или иной фразы.
В частности, я рекомендую поиск в книгах на Google
\begin{center}
\href{http://books.google.com/bkshp?tab=wp}
{http://books.google.com/bkshp?tab=wp}
\end{center}

Нередко в процессе перевода я понимаю, что я должен изменить
исходную фразу для того, чтобы она более адекватно выражала мою мысль
на обоих языках.
Поэтому я рассматриваю перевод как
неотъемлемую часть работы над текстом
и творческий процесс.

Изначально этот словарь был задуман как непосредственный помощник при написании статей.
Хотя в интернете есть немало научных текстов, поиск того или иного термина или определения
занимает немало времени.
Поэтому я решил сохранять перевод терминов и определений, которые
представляли для меня интерес и которые я встречал
в книгах, изданных на обоих языках.
Так появился словарь.
Словарь составлялся постепенно.
Постепенно объём словаря рос.
Я решил сделать словарь доступным для других.
Поэтому я постарался так же учесть пожелания, высказанные в ответ на
предыдущие издания словаря.
Этот словарь включает в себя только физические и математические
термины. Если слово имеет альтернативный перевод или другой смысл
в обычной речи, то эта информация в словаре отсутствует.

Я также включил в книгу двуязычный
именной указатель. В этот список я также включал только те имена, которые
встречал и в русском, и в английском текстах. Однако здесь есть ещё одна проблема,
с которой я столкнулся при переводе статей. Далеко не всегда очевидна
связь между записью и произношением имени. Было бы неплохо добавить транскрипцию,
но я этого не сделал, так как у меня нет надёжного источника.

\chapter{English Russian Dictionary}%

\section{A}%
\begin{description}%

\ERitem%
   {$A$\Hyph valued function}%
   {$A$\Hyph значная функция}%

\ERitem%
   {Abelian group}%
   {абелевая группа}%

\ERitem%
   {absolute convergence}%
   {абсолютная сходимость}%

\ERitem%
   {absolute value}%
   {абсолютная величина; норма}%
\ERexample%
   {Absolute value on skew field $D$ is a map\[d\in D\rightarrow |d|\in R\]which satisfies the following axioms\begin{itemize}\item $|a|\ge 0$\item $|a|=0$ if, and only if, $a=0$\item $|ab|=|a|\ |b|$\item $|a+b|\le|a|+|b|$\end{itemize}}%
   {Норма на теле $D$ - это отображение\[d\in D\rightarrow |d|\in R\]такое, что\begin{itemize}\item $|a|\ge 0$\item $|a|=0$ равносильно $a=0$\item $|ab|=|a|\ |b|$\item $|a+b|\le|a|+|b|$\end{itemize}}%
   {}{}{}{}

\ERitem%
   {absorption of photon}%
   {поглощение фотона}%

\ERitem%
   {acceleration}%
   {ускорение}%

\ERitem%
   {accelerator}%
   {ускоритель}%

\ERitem%
   {according to theorem 2.1}%
   {согласно теореме 2.1}%
\ERexample%
   {According to theorem 2.1, triangles $ABC$ and $DBC$ are equal.}%
   {Согласно теореме 2.1 треугольники $ABC$ и $DBC$ равны.}%
   {}{}{}{}

\ERitem%
   {acute angle}%
   {острый угол}%

\ERitem%
   {additive group}%
   {аддитивная группа}%

\ERitem%
   {adjacent angle}%
   {смежный угол}%

\ERitem%
   {adjoin}%
   {присоединить}%
\ERexample%
   {To derive equations of motion of a charged particle we adjoin Lorentz equations to Maxwell equations.}%
   {Для того, чтобы вывести уравнения движения заряда, мы присоединим уравнения Лоренца к уравнения Максвелла.}%
   {}{}{}{}

\ERitem%
   {adjoint group}%
   {присоединённая группа}%

\ERitem%
   {algebra bundle}%
   {расслоенная алгебра}%

\ERitem%
   {algebraic}%
   {алгебраический}%

\ERitem%
   {algebraic complement of matrix}%
   {алгебраическое дополнение матрицы}%

\ERitem%
   {algebraic extension}%
   {алгебраическое расширение}%

\ERitem%
   {amplitude}%
   {амплитуда}%

\ERitem%
   {analytic function}%
   {аналитическая функция}%

\ERitem%
   {angle}%
   {угол}%

\ERitem%
   {angle of incidence}%
   {угол падения}%

\ERitem%
   {angle of reflection}%
   {угол отражения}%

\ERitem%
   {angle of refraction}%
   {угол преломления}%

\ERitem%
   {angular momentum}%
   {момент количества движения}%

\ERitem%
   {anholonomic coordinates}%
   {неголономные координаты}%

\ERitem%
   {anholonomity}%
   {неголономность}%

\ERitem%
   {anholonomity object}%
   {объект неголономности}%

\ERitem%
   {annihilation operator}%
   {оператор уничтожения}%

\ERitem%
   {annihilator}%
   {аннулятор}%

\ERitem%
   {apocentre}%
   {апоцентр}%

\ERitem%
   {approximation}%
   {приближение}%

\ERitem%
   {arccosecant}%
   {арккосеканс}%

\ERitem%
   {arccosine}%
   {арккосинус}%

\ERitem%
   {arccotangent}%
   {арккотангенс}%

\ERitem%
   {arcsecant}%
   {арксеканс}%

\ERitem%
   {arcsine}%
   {арксинус}%

\ERitem%
   {arctangent}%
   {арктангенс}%

\ERitem%
   {arity}%
   {арность}%

\ERitem%
   {as small as we please}%
   {сколь угодно малый}%

\ERitem%
   {associative}%
   {ассоциативный}%

\ERitem%
   {associative law}%
   {закон ассоциативности}%

\ERitem%
   {associativity}%
   {ассоциативность}%

\ERitem%
   {at first glance}%
   {на первый взгляд}%

\ERitem%
   {at first sight}%
   {на первый взгляд}%

\ERitem%
   {at least}%
   {по крайней мере}%
\ERexample%
   {At least in the neighborhood of the identity.}%
   {По крайней мере, в окрестности единичного элемента.}%
   {}{}{}{}

\ERitem%
   {attractor}%
   {аттрактор}%

\ERitem%
   {auto parallel line}%
   {автопараллельная кривая}%

\ERitem%
   {axiom of choice}%
   {аксиома выбора}%

\end{description}%

\section{B}%
\begin{description}%

\ERitem%
   {Banach algebra}%
   {банахова алгебра}%

\ERitem%
   {Banach space}%
   {банахово пространство}%

\ERitem%
   {base of fibered correspondence}%
   {база расслоенного соответствия}%

\ERitem%
   {base of topology}%
   {базис топологии}%

\ERitem%
   {basis for vector space}%
   {базис в векторном пространстве; базис векторного пространства}%

\ERitem%
   {basis of vector space}%
   {базис в векторном пространстве; базис векторного пространства}%

\ERitem%
   {behavior}%
   {поведение}%

\ERitem%
   {Bell's theorem}%
   {теорема Белла}%

\ERitem%
   {bijection}%
   {биекция}%

\ERitem%
   {bimodule}%
   {бимодуль}%

\ERitem%
   {binary}%
   {бинарный}%

\ERitem%
   {Bott periodicity}%
   {периодичность Ботта}%

\ERitem%
   {boundary}%
   {граница}%

\ERitem%
   {boundary conditions}%
   {граничные условия}%

\ERitem%
   {By Theorem 2,1}%
   {согласно теореме 2.1}%
\ERexample%
   {By Theorem 2.1, $a=b$.}%
   {Согласно теореме 2.1, $a=b$.}%
   {}{}{}{}

\end{description}%

\section{C}%
\begin{description}%

\ERitem%
   {canonical map}%
   {каноническое отображение}%

\ERitem%
   {Cartesian coordinate system}%
   {декартова система координат}%

\ERitem%
   {Cartesian power}%
   {декартова степень}%

\ERitem%
   {Cartesian product}%
   {декартово произведение}%

\ERitem%
   {category}%
   {категория}%

\ERitem%
   {Cauchy sequence}%
   {последовательность Коши}%

\ERitem%
   {causal relationship}%
   {причинно-следственная связь}%

\ERitem%
   {causal scalar field}%
   {причинное скалярное поле}%

\ERitem%
   {causal vector field}%
   {причинное векторное поле}%

\ERitem%
   {chain rule}%
   {правило дифференцирования сложной функции}%

\ERitem%
   {change of coordinates}%
   {замена координат}%

\ERitem%
   {change of variable}%
   {замена переменной}%

\ERitem%
   {chaos}%
   {хаос}%

\ERitem%
   {chart over $U$}%
   {тривиализация над $U$}%

\ERitem%
   {closure of the set}%
   {замыкание множества}%

\ERitem%
   {closure operator}%
   {оператор замыкания}%

\ERitem%
   {closure system}%
   {система замыканий}%

\ERitem%
   {cluster point}%
   {точка прикосновения}%

\ERitem%
   {coarsest topology}%
   {самая слабая топология}%

\ERitem%
   {cofactor of matrix}%
   {алгебраическое дополнение матрицы}%

\ERitem%
   {colloquia}%
   {коллоквиумы}%

\ERitem%
   {colloquium}%
   {коллоквиум}%

\ERitem%
   {column vector}%
   {вектор столбец}%

\ERitem%
   {combinatorics}%
   {комбинаторика}%

\ERitem%
   {commutative diagram}%
   {коммутативная диаграмма}%

\ERitem%
   {commutativity}%
   {коммутативность}%

\ERitem%
   {commutator}%
   {коммутатор}%

\ERitem%
   {commute}%
   {коммутирует}%
\ERexample%
   {Position and momentum operators do not commute.}%
   {Операторы положения и импульса не коммутируют.}%
   {}{}{}{}

\ERitem%
   {compact\hyph open topology}%
   {компактно\hyph открытая топология}%

\ERitem%
   {comparable topology}%
   {сравнимые топологии}%

\ERitem%
   {complete division ring}%
   {полное тело}%

\ERitem%
   {complete field}%
   {полное поле}%

\ERitem%
   {complete lattice}%
   {полная структура}%

\ERitem%
   {complete space}%
   {полное пространство}%

\ERitem%
   {complete system}%
   {полная система}%

\ERitem%
   {completely integrable}%
   {вполне интегрируемый}%

\ERitem%
   {complex field}%
   {поле комплексных чисел}%

\ERitem%
   {componentwise}%
   {покомпонентно}%

\ERitem%
   {conditions of integrability}%
   {условия интегрируемости}%

\ERitem%
   {conformal transformation}%
   {конформное преобразование}%

\ERitem%
   {congruence}%
   {конгруэнтность; согласованность}%

\ERitem%
   {conjugate quaternion}%
   {сопряжённый кватернион}%

\ERitem%
   {connected group}%
   {связная группа}%

\ERitem%
   {connection coefficients}%
   {коэффициенты связности}%

\ERitem%
   {conservation law}%
   {закон сохранения}%

\ERitem%
   {consider}%
   {рассматривать}%
\ERexample%
   {Consider correspondence from set $A$ to set $B$.}%
   {Рассмотрим соответствие $\Phi$ из множества $A$ в множество $B$.}%
   {}{}{}{}

\ERitem%
   {continuous in neighborhood}%
   {непрерывен в окрестности}%

\ERitem%
   {continuous in $x$}%
   {непрерывный по $x$}%

\ERitem%
   {contradiction}%
   {противоречие}%
\ERexample%
   {The contradiction completes the proof of the theorem.}%
   {Полученное противоречие доказывает теорему.}%
   {}{}{}{}

\ERitem%
   {contravariant}%
   {контравариантный}%

\ERitem%
   {convection}%
   {конвекция}%

\ERitem%
   {convention}%
   {соглашение}%
\ERexample%
   {We use the convention that we present any set of vectors of the vector space as a row.}%
   {Мы пользуемся соглашением, что в заданном векторном пространстве мы представляем любое семейство векторов в виде строки.}%
   {}{}{}{}

\ERitem%
   {converge}%
   {сходиться}%
\ERexample%
   {Filter $\mathfrak{F}$ converges to $x$.}%
   {Фильтр $\mathfrak{F}$ сходится к $x$.}%
   {}{}{}{}

\ERitem%
   {conversely}%
   {обратно}%

\ERitem%
   {convex function}%
   {выпуклая функция}%

\ERitem%
   {coordinate chart}%
   {координатная карта}%

\ERitem%
   {correlation}%
   {корреляция}%

\ERitem%
   {correspondence from $A$ to $B$}%
   {соответствие из $A$ в $B$}%

\ERitem%
   {cosecant}%
   {косеканс}%

\ERitem%
   {cosine}%
   {косинус}%

\ERitem%
   {cotangent}%
   {котангенс}%

\ERitem%
   {countable set}%
   {счётное множество}%

\ERitem%
   {counter}%
   {счётчик}%

\ERitem%
   {covariant}%
   {ковариантный}%

\ERitem%
   {cover}%
   {покрытие}%

\ERitem%
   {covering space}%
   {накрытие}%
\ERexample%
   {Consider the covering space $R\rightarrow S^1$ of the circle $S^1$ defined by $p(t)=(\sin t,\cos t)$ for any $t\in R$.}%
   {Рассмотрим накрытие $R\rightarrow S^1$ окружности $S^1$, определённое формулой $p(t)=(\sin t,\cos t)$ для любого $t\in R$.}%
   {}{}{}{}

\ERitem%
   {Cramer's Rule}%
   {правило Крамера}%

\ERitem%
   {creation operator}%
   {оператор рождения}%

\ERitem%
   {crystal lattice}%
   {кристаллическая решётка}%

\ERitem%
   {curve}%
   {кривая}%

\ERitem%
   {curvilinear coordinates}%
   {криволинейные координаты}%

\ERitem%
   {cycle}%
   {цикл}%

\ERitem%
   {cyclic group}%
   {циклическая группа}%

\end{description}%

\section{D}%
\begin{description}%

\ERitem%
   {$D$\Hyph vector space}%
   {$D$\Hyph векторное пространство}%

\ERitem%
   {decomposition of map}%
   {разложение отображения}%

\ERitem%
   {define}%
   {определяет}%
\ERexample%
   {This equation defines the inverse transformation.}%
   {Это уравнение определяет обратное преобразование.}%
   {}{}{}{}

\ERitem%
   {degree of map}%
   {степень отображения}%

\ERitem%
   {denominator}%
   {знаменатель}%
\ERexample%
   {Let us reduce items to a common denominator.}%
   {Приведём слагаемые к общему знаменателю.}%
   {}{}{}{}

\ERitem%
   {dependence}%
   {зависимость}%

\ERitem%
   {derivative of second or greater order with respect}%
   {производная второго или более высокого порядка по }%

\ERitem%
   {develop equation}%
   {вывести уравнение}%

\ERitem%
   {diagram of correspondences}%
   {диаграмма соответствий}%

\ERitem%
   {difference}%
   {разность}%

\ERitem%
   {difference module}%
   {фактор модуль}%

\ERitem%
   {differentiability}%
   {дифференцируемость}%

\ERitem%
   {differentiable function}%
   {дифференцируемая функция}%

\ERitem%
   {differentiable in the Fr\'echet sense}%
   {дифференцируемый по Фреше}%

\ERitem%
   {differentiable in the G\^ateaux sense}%
   {дифференцируемый по Гато}%

\ERitem%
   {differentiate the function with respect to $x$}%
   {дифференцировать фун\-к\-цию по $x$}%

\ERitem%
   {diffraction}%
   {дифракция}%

\ERitem%
   {diffusion}%
   {диффузия}%

\ERitem%
   {discrete space}%
   {дискретное пространство}%

\ERitem%
   {discrete topology}%
   {дискретная топология}%

\ERitem%
   {distributive law}%
   {закон дистрибутивности}%

\ERitem%
   {distributive property of multiplication over addition}%
   {ди\-с\-т\-ри\-бу\-тив\-ность умножения относительно сложения}%

\ERitem%
   {division ring}%
   {тело (кольцо с делением)}%

\ERitem%
   {domain}%
   {область определения}%

\ERitem%
   {Doppler shift}%
   {эффект Доплера}%

\ERitem%
   {downstairs}%
   {внизу}%
\ERexample%
   {We sum over any index which appears twice in the same term, once upstairs and once downstairs.}%
   {Подразумевается сумма по любому индексу, появляющемуся дважды в одном и том же слагаемом, один раз вверху, другой - внизу.}%
   {}{}{}{}

\ERitem%
   {dual module}%
   {дуальный модуль}%

\ERitem%
   {dual space}%
   {дуальное пространство}%

\ERitem%
   {dynamics}%
   {динамика}%

\end{description}%

\section{E}%
\begin{description}%

\ERitem%
   {eccentricity}%
   {эксцентриситет}%

\ERitem%
   {eclipse}%
   {затмение}%

\ERitem%
   {elementary particle}%
   {элементарная частица}%

\ERitem%
   {emission of photon}%
   {излучение фотона}%

\ERitem%
   {endomorphism}%
   {эндоморфизм}%

\ERitem%
   {energy}%
   {энергия}%

\ERitem%
   {engine}%
   {двигатель}%

\ERitem%
   {enhanced}%
   {расширенный}%

\ERitem%
   {entropy}%
   {энтропия}%

\ERitem%
   {envelope of a family of plane curves}%
   {огибающая семейства плоских кривых}%

\ERitem%
   {enveloping algebra}%
   {обвёртывающая алгебра}%

\ERitem%
   {equation is satisfied identically}%
   {уравнение удовлетворяется тождественно}%

\ERitem%
   {equivalence class}%
   {класс эквивалентности}%

\ERitem%
   {equivalence relation}%
   {эквивалентность}%

\ERitem%
   {Erlanger Program}%
   {Эрлангенская программа}%

\ERitem%
   {essential parameters in a set of functions}%
   {существенные параметры семейства функций}%

\ERitem%
   {Euclidean metric}%
   {эвклидова метрика}%

\ERitem%
   {Euclidean space}%
   {эвклидово пространство}%

\ERitem%
   {evaluating by equating $x$ to the $a$}%
   {подстановка $a$ вместо $x$}%

\ERitem%
   {event horizon}%
   {горизонт событий}%

\ERitem%
   {event space}%
   {пространство событий}%

\ERitem%
   {evidence}%
   {очевидность}%

\ERitem%
   {evidently}%
   {очевидно}%
\ERexample%
   {Evidently $x=1$ is the root of the equation.}%
   {Очевидно, $x=1$ является корнем уравнения.}%
   {}{}{}{}

\ERitem%
   {exact sequence of modules}%
   {точная последовательность модулей}%

\ERitem%
   {extension field}%
   {расширение поля}%

\ERitem%
   {extension of correspondence}%
   {продолжение соответствия}%

\ERitem%
   {exterior differential}%
   {внешний дифференциал}%

\ERitem%
   {exterior product}%
   {внешнее произведение}%

\ERitem%
   {external algebra}%
   {внешняя алгебра}%

\ERitem%
   {external power}%
   {внешняя степень}%

\ERitem%
   {extremal}%
   {экстремальный}%

\ERitem%
   {extreme}%
   {экстремальный}%

\ERitem%
   {extreme line}%
   {экстремальная кривая}%

\end{description}%

\section{F}%
\begin{description}%

\ERitem%
   {factor}%
   {множитель; разложить на множители; сомножитель}%
\ERexample%
   {To factor a polynomial means to find two or more po\-ly\-no\-mials whose product is the given polynomial.}%
   {Чтобы разложить многочлен на множители, необходимо найти два или более многочленов, произведение которых есть данный многочлен.}%
   {}{}{}{}

\ERitem%
   {factor group}%
   {факторгруппа}%

\ERitem%
   {factorization}%
   {разложение на множители}%

\ERitem%
   {fiber}%
   {слой}%

\ERitem%
   {fibered correspondence}%
   {расслоенное соответствие}%

\ERitem%
   {fibered map}%
   {морфизм расслоений}%

\ERitem%
   {fibered product}%
   {расслоенное произведение}%

\ERitem%
   {field-strength tensor}%
   {тензор напряжённости поля}%

\ERitem%
   {filter}%
   {фильтр}%

\ERitem%
   {filter base}%
   {базис фильтра}%

\ERitem%
   {finest topology}%
   {самая сильная топология}%

\ERitem%
   {finite dimensional}%
   {конечномерный}%

\ERitem%
   {finite set}%
   {конечное множество}%

\ERitem%
   {Finsler metric}%
   {финслеровая метрика}%

\ERitem%
   {Finslerian metric}%
   {финслеровая метрика}%

\ERitem%
   {force}%
   {сила}%

\ERitem%
   {the Fr\'echet derivative}%
   {производная Фреше}%

\ERitem%
   {the Fr\'echet differential}%
   {дифференциал Фреше}%

\ERitem%
   {free representation}%
   {свободное представление}%

\ERitem%
   {frequency}%
   {частота}%

\ERitem%
   {friction}%
   {сила трения; трение}%

\ERitem%
   {function $f$ of $x$}%
   {функция $f$ от $x$}%

\ERitem%
   {functional}%
   {функционал}%

\ERitem%
   {functor}%
   {функтор}%

\ERitem%
   {fundamental sequence}%
   {фундаментальная последовательность}%

\end{description}%

\section{G}%
\begin{description}%

\ERitem%
   {$G$-principal bundle}%
   {главное $G$-расслоение}%

\ERitem%
   {galaxy}%
   {галактика}%

\ERitem%
   {the G\^ateaux derivative}%
   {производная Гато}%

\ERitem%
   {the G\^ateaux differential}%
   {дифференциал Гато}%

\ERitem%
   {gauge invariance}%
   {калибровочная инвариантность}%

\ERitem%
   {general relativity}%
   {общая теория относительности}%

\ERitem%
   {generally speaking}%
   {вообще говоря}%

\ERitem%
   {generated}%
   {порождённый}%
\ERexample%
   {Algebra $A$ generated by the set $S$ is a $K$\Hyph algebra}%
   {Алгебра $A$, порождённая множеством $S$, является $K$\Hyph алгеброй.}%
   {}{}{}{}

\ERitem%
   {generator}%
   {образующая}%

\ERitem%
   {geometry}%
   {геометрия}%

\ERitem%
   {Global Positioning System}%
   {глобальная система позиционирования}%

\ERitem%
   {gluing functions}%
   {функции склеивания}%

\ERitem%
   {Gram-Schmidt orthogonalization procedure}%
   {процесс ортогонализации Грама--Шмидта}%

\ERitem%
   {graph}%
   {граф}%

\ERitem%
   {graph theory}%
   {теория графов}%

\ERitem%
   {gravity probe}%
   {гравитационный зонд}%

\ERitem%
   {group bundle}%
   {расслоенная группа}%

\end{description}%

\section{H}%
\begin{description}%

\ERitem%
   {Hadamard inverse}%
   {обращение Адамара}%

\ERitem%
   {has relevance to}%
   {имеет отношение к}%

\ERitem%
   {head of vector}%
   {конец вектора}%

\ERitem%
   {helical structure}%
   {спиральная структура}%

\ERitem%
   {helicity}%
   {спиральность}%

\ERitem%
   {hermitian form}%
   {эрмитова форма}%

\ERitem%
   {highest common factor of $p$ and $q$}%
   {наибольший общий делитель $p$ и $q$}%

\ERitem%
   {holonomic coordinates}%
   {голономные координаты}%

\ERitem%
   {homeomorphic}%
   {гомеоморфный}%

\ERitem%
   {homeomorphism}%
   {гомеоморфизм}%

\ERitem%
   {homogeneous}%
   {однородный}%

\ERitem%
   {homogeneous Lorentz group}%
   {однородная группа Лоренца}%

\ERitem%
   {homology}%
   {гомология}%

\ERitem%
   {homomorphism}%
   {гомоморфизм}%

\ERitem%
   {homotopic}%
   {гомотопный}%

\ERitem%
   {homotopy}%
   {гомотопия}%

\ERitem%
   {hyperfine splitting}%
   {сверхтонкое расщепление}%

\ERitem%
   {hyperplane}%
   {гиперплоскость}%

\end{description}%

\section{I}%
\begin{description}%

\ERitem%
   {identical particles}%
   {тождественные частицы}%

\ERitem%
   {identification}%
   {отождествление}%

\ERitem%
   {identification morphism}%
   {морфизм отождествления}%

\ERitem%
   {identity}%
   {единичный элемент}%

\ERitem%
   {iff}%
   {тогда и только тогда, когда}%
\ERexample%
   {$a=0$ iff $a_{\gi i}^{\gi j}=0$ for any ${\gi i}$, ${\gi j}$.}%
   {$a=0$ тогда и только тогда, когда $a_{\gi i}^{\gi j}=0$ для любых ${\gi i}$, ${\gi j}$.}%
   {}{}{}{}

\ERitem%
   {image under map}%
   {образ при отображении}%
\ERexample%
   {We define the image of the set $A$ under correspondence $\Phi$ according to law\[A\Phi=\{b:(a,b)\in\Phi,a\in A\}\]}%
   {Мы определим образ множества $A$ при соответствии $\Phi$ согласно равенству\[A\Phi=\{b:(a,b)\in\Phi,a\in A\}\]}%
   {}{}{}{}

\ERitem%
   {in a similar way}%
   {подобным образом}%
\ERexample%
   {In a similar way, we can introduce a coordinate reference frame.}%
   {Подобным образом мы можем определить координатную систему отсчёта.}%
   {}{}{}{}

\ERitem%
   {in general}%
   {вообще говоря}%
\ERexample%
   {However in general this product is not linear map.}%
   {Однако, вообще говоря, это отображением не является линейным.}%
   {}{}{}{}

\ERitem%
   {indicatrix}%
   {индикатриса}%

\ERitem%
   {inequation}%
   {неравенство}%

\ERitem%
   {infinitesimal}%
   {бесконечно малая величина; бесконечно малый}%

\ERitem%
   {inhomogeneous}%
   {неоднородный}%

\ERitem%
   {inhomogeneous Lorentz group}%
   {неоднородная группа Лоренца}%

\ERitem%
   {injection}%
   {инъекция}%

\ERitem%
   {insulator}%
   {диэлектрик; изолятор}%

\ERitem%
   {interaction}%
   {взаимодействие}%
\ERexample%
   {Based on differential geometry, general relativity describes gravitational interaction.}%
   {Общая теория относительности описывает гравитационное взаимодействие, опираясь на дифференциальную геометрию.}%
   {}{}{}{}

\ERitem%
   {interaction picture}%
   {представление взаимодействия}%

\ERitem%
   {interference}%
   {интерференция}%

\ERitem%
   {inverse transformation}%
   {обратное преобразование}%

\ERitem%
   {irreducible representation}%
   {неприводимое представление}%

\ERitem%
   {is related to}%
   {имеет отношение к}%

\ERitem%
   {isotropic vector}%
   {изотропный вектор}%

\ERitem%
   {it is evident that}%
   {очевидно, что}%
\ERexample%
   {From (2.2), it is evident that any solution of (2.7) satisfies (2.9).}%
   {На основании (2.2) очевидно, что любое решение уравнения (2.7) удовлетворяет (2.9).}%
   {}{}{}{}

\end{description}%

\section{J}%
\begin{description}%

\ERitem%
   {Jacobian}%
   {якобиан}%

\ERitem%
   {Jacobian matrix}%
   {матрица Якоби}%

\end{description}%

\section{K}%
\begin{description}%

\ERitem%
   {kernel}%
   {ядро (отображения)}%

\ERitem%
   {Kerr metric}%
   {метрика Керра}%

\ERitem%
   {kinematics}%
   {кинематика}%

\ERitem%
   {Klein bottle}%
   {бутылка Клейна}%

\ERitem%
   {knot}%
   {узел}%

\end{description}%

\section{L}%
\begin{description}%

\ERitem%
   {Lagrangian}%
   {лагранжиан}%

\ERitem%
   {lattice}%
   {структура (алгебраическая система)}%

\ERitem%
   {left side of equation}%
   {левая часть равенства}%

\ERitem%
   {left-distributive}%
   {дистрибутивен слева}%

\ERitem%
   {L'H\^ospital's rule}%
   {правило Лопиталя}%

\ERitem%
   {lift of correspondence}%
   {лифт соответствия; подъём соответствия}%

\ERitem%
   {lift of morphism}%
   {лифт морфизма; подъём морфизма}%

\ERitem%
   {lift of vector field}%
   {лифт векторного поля; подъём векторного поля}%

\ERitem%
   {limit}%
   {предел}%

\ERitem%
   {limit of correspondence with respect to the filter}%
   {предел соответствия по фильтру}%

\ERitem%
   {limit of sequence}%
   {предел последовательности}%

\ERitem%
   {limit point}%
   {предельная точка}%

\ERitem%
   {limit set}%
   {предельное множество}%

\ERitem%
   {linearly dependent}%
   {линейно зависимые}%

\ERitem%
   {linearly independent}%
   {линейно независимые}%

\ERitem%
   {little group}%
   {малая группа}%

\ERitem%
   {locally compact space}%
   {локально компактное пространство}%

\ERitem%
   {loop (quasigroup with unit element)}%
   {лупа (квазигруппа с единицей)}%

\ERitem%
   {Lorentz transformation}%
   {преобразование Лоренца}%

\ERitem%
   {lower index}%
   {нижний индекс}%

\end{description}%

\section{M}%
\begin{description}%

\ERitem%
   {the Mach principle}%
   {принцип Маха}%

\ERitem%
   {manifolds with affine connections}%
   {пространство аффинной связности}%

\ERitem%
   {mapping}%
   {отображение}%

\ERitem%
   {mass}%
   {масса}%

\ERitem%
   {massive particle}%
   {массивная частица}%

\ERitem%
   {massless particle}%
   {безмассовая частица}%

\ERitem%
   {mathematical}%
   {математический}%

\ERitem%
   {mathematician}%
   {математик}%

\ERitem%
   {mathematics}%
   {математика}%

\ERitem%
   {mean value theorem}%
   {теорема о конечных приращениях}%

\ERitem%
   {measure}%
   {измерять}%

\ERitem%
   {measurement}%
   {измерение}%

\ERitem%
   {metric-affine manifold}%
   {аффинно-метрическое многообразие}%

\ERitem%
   {Milky Way}%
   {Млечный Путь}%

\ERitem%
   {mixed system}%
   {смешанная система}%

\ERitem%
   {Moebius band}%
   {лист Мёбиуса}%

\ERitem%
   {momentum}%
   {импульс}%

\ERitem%
   {momentum operator}%
   {оператор импульса}%

\ERitem%
   {monic polynomial}%
   {приведенный многочлен; унитарный многочлен}%

\ERitem%
   {monomial}%
   {одночлен}%

\ERitem%
   {monotone function}%
   {монотонная функция}%

\ERitem%
   {monotonic function}%
   {монотонная функция}%

\ERitem%
   {multiple root}%
   {кратный корень}%

\ERitem%
   {multiplication}%
   {умножение}%

\ERitem%
   {multiplication table}%
   {таблица умножения}%

\ERitem%
   {multiplicative group}%
   {мультипликативная группа}%

\ERitem%
   {multiplicity of $x$ in $f$}%
   {кратность $x$ в $f$}%
\ERexample%
   {If the multiplicity of $a$ is greater than $1$, then $a$ is called a multiple root.}%
   {Если кратность $a$ больше, чем $1$, то $a$ называется кратным корнем.}%
   {}{}{}{}

\ERitem%
   {multiply by $2$}%
   {умножить на $2$}%

\ERitem%
   {multiply by $b$}%
   {умножить на $b$}%

\ERitem%
   {muon}%
   {мюон}%

\ERitem%
   {mutually orthogonal}%
   {взаимно ортогональные; ортогональные друг другу; попарно ортогональные}%

\ERitem%
   {mutually perpendicular}%
   {взаимно перпендикулярные; перпендикулярные друг другу}%

\end{description}%

\section{N}%
\begin{description}%

\ERitem%
   {name index}%
   {именной указатель}%

\ERitem%
   {natural mapping}%
   {естественное отображение}%

\ERitem%
   {natural morphism}%
   {естественный морфизм}%

\ERitem%
   {necessary and sufficient}%
   {необходимо и достаточно}%
\ERexample%
   {In order that any $x\in A$ is root of the system of linear equations\[a_i^jx^i=0\]necessary and sufficient $a^{\gij}_{\gii}=0$.}%
   {Для того, чтобы любое $x\in A$ было корнем системы линейных уравнений\[a_i^jx^i=0\]необходимо и достаточно, чтобы $a^{\gij}_{\gii}=0$.}%
   {}{}{}{}

\ERitem%
   {neighborhood}%
   {окрестность}%

\ERitem%
   {neutrino}%
   {нейтрино}%

\ERitem%
   {neutron}%
   {нейтрон}%

\ERitem%
   {neutron star}%
   {нейтронная звезда}%

\ERitem%
   {non-Abelian group}%
   {неабелевая группа}%

\ERitem%
   {nondegenerate form}%
   {невырожденная форма}%

\ERitem%
   {nontrivial}%
   {нетривиальный}%

\ERitem%
   {norm}%
   {норма}%
\ERexample%
   {Norm on $\Omega$\Hyph group $A$ is a map\[d\in A\rightarrow \|d\|\in R\]such that\begin{itemize}\item $\|a\|\ge 0$\item $\|a\|=0$ iff, $a=0$\item $\|a+b\|\le\|a\|+\|b\|$\item $\|-a\|=\|a\|$\end{itemize}}%
   {Норма на $\Omega$\Hyph группе $A$ - это отображение\[d\in A\rightarrow \|d\|\in R\]такое, что\begin{itemize}\item $\|a\|\ge 0$\item $\|a\|=0$ равносильно $a=0$\item $\|a+b\|\le\|a\|+\|b\|$\item $\|-a\|=\|a\|$\end{itemize}}%
   {}{}{}{}

\ERitem%
   {normed space}%
   {нормированное пространство}%

\ERitem%
   {nucleus}%
   {ядро (атома)}%

\ERitem%
   {numerator}%
   {числитель}%

\end{description}%

\section{O}%
\begin{description}%

\ERitem%
   {obtain by differentiating}%
   {получить дифференцированием}%

\ERitem%
   {obtuse angle}%
   {тупой угол}%

\ERitem%
   {Occam's razor}%
   {бритва Оккама}%

\ERitem%
   {operate}%
   {действовать}%
\ERexample%
   {Operating on equation (1) with operator $V$ yields an integral equation.}%
   {Подействовав на уравнение (1) оператором $V$, получим интегральное уравнение.}%
   {}{}{}{}

\ERitem%
   {opposite preordering}%
   {противоположная предупорядоченность}%

\ERitem%
   {ordered set}%
   {упорядоченное множество}%

\ERitem%
   {ordering}%
   {упорядоченность}%

\ERitem%
   {orthonormal basis}%
   {ортонормированный базис}%

\end{description}%

\section{P}%
\begin{description}%

\ERitem%
   {parallel transport}%
   {параллельный перенос}%

\ERitem%
   {parallelepiped}%
   {параллелепипед}%

\ERitem%
   {parity}%
   {чётность}%

\ERitem%
   {partial differential equation}%
   {уравнение в частных производных}%

\ERitem%
   {partial ordering}%
   {частичная упорядоченность}%

\ERitem%
   {partition of unity}%
   {разбиение единицы; разложение единицы}%

\ERitem%
   {Pasch's axiom}%
   {аксиома Паша}%

\ERitem%
   {passage to the limit}%
   {предельный переход}%

\ERitem%
   {pericentre}%
   {перицентр}%

\ERitem%
   {perturbation}%
   {возмущение}%

\ERitem%
   {pfaffian derivative}%
   {пфаффова производная}%

\ERitem%
   {phenomena}%
   {явления}%

\ERitem%
   {phenomenon}%
   {явление}%

\ERitem%
   {photon}%
   {фотон}%

\ERitem%
   {physical}%
   {физический}%

\ERitem%
   {physicist}%
   {физик}%

\ERitem%
   {physics}%
   {физика}%

\ERitem%
   {point}%
   {точечный; точка}%

\ERitem%
   {polology}%
   {полология}%

\ERitem%
   {polyadditive map}%
   {полиаддитивное отображение}%

\ERitem%
   {polylinear form}%
   {полилинейная форма}%

\ERitem%
   {polynomial}%
   {многочлен; полином}%

\ERitem%
   {polyvector}%
   {поливектор}%

\ERitem%
   {poolback bundle}%
   {обратный образ расслоения; прообраз расслоения}%

\ERitem%
   {positive definite form}%
   {положительно определённая форма}%

\ERitem%
   {positive integer}%
   {натуральное число}%

\ERitem%
   {power of set}%
   {мощность множества}%

\ERitem%
   {precession}%
   {прецессия}%

\ERitem%
   {preordering}%
   {предупорядоченность}%

\ERitem%
   {prime ideal}%
   {простой идеал}%

\ERitem%
   {principal bundle}%
   {главное расслоение}%

\ERitem%
   {probability}%
   {вероятность}%

\ERitem%
   {problem}%
   {задача}%

\ERitem%
   {proceeding in this way}%
   {продолжая таким образом; продолжив этот процесс}%

\ERitem%
   {projection}%
   {проекция}%

\ERitem%
   {projective plane}%
   {проективная плоскость}%

\ERitem%
   {proof by induction}%
   {доказательство по индукции}%

\ERitem%
   {propagation}%
   {распространение}%

\ERitem%
   {proper state}%
   {собственное состояние}%

\ERitem%
   {proper value}%
   {собственное значение}%

\ERitem%
   {pseudo\Hyph Euclidean space}%
   {псевдоевклидовое пространство}%

\ERitem%
   {pulsar}%
   {пульсар}%

\end{description}%

\section{Q}%
\begin{description}%

\ERitem%
   {quadratic equation}%
   {квадратное уравнение}%

\ERitem%
   {quantum entanglement}%
   {квантовая запутанность}%

\ERitem%
   {quark}%
   {кварк}%

\ERitem%
   {quasar}%
   {квазар}%

\ERitem%
   {quasigroup}%
   {квазигруппа}%

\ERitem%
   {quaternion}%
   {кватернион}%

\ERitem%
   {quotient bundle}%
   {фактор расслоение}%

\ERitem%
   {quotient group}%
   {факторгруппа}%

\ERitem%
   {quotient ring}%
   {факторкольцо}%

\ERitem%
   {quotient set}%
   {фактор множество}%

\ERitem%
   {quotient topology}%
   {фактортопология}%

\end{description}%

\section{R}%
\begin{description}%

\ERitem%
   {radiation belt}%
   {радиационный пояс}%

\ERitem%
   {range}%
   {множество значений; область значений}%

\ERitem%
   {rational field}%
   {поле рациональных чисел}%

\ERitem%
   {real function}%
   {числовая функция}%

\ERitem%
   {real valued function}%
   {числовая функция}%

\ERitem%
   {reciprocal image}%
   {обратный образ}%

\ERitem%
   {reduced Cartesian product of bundles}%
   {приведенное декартово произведение расслоений}%

\ERitem%
   {reduced fibered correspondence}%
   {приведенное расслоенное соответствие}%

\ERitem%
   {reduction of similar terms}%
   {приведение подобных}%

\ERitem%
   {reference frame}%
   {система отсчёта}%

\ERitem%
   {reflexive}%
   {рефлексивный}%

\ERitem%
   {regression}%
   {регрессия}%

\ERitem%
   {relation}%
   {отношение}%

\ERitem%
   {relevance}%
   {важность; значимость; существенность}%
\ERexample%
   {Distinction between Lorentz and Poincar\'e groups is of no relevance here.}%
   {Различие между группами Лоренца и Пуанкаре для нас сейчас не важно.}%
   {}{}{}{}

\ERitem%
   {remainder}%
   {остаток}%

\ERitem%
   {repeated root}%
   {кратный корень}%

\ERitem%
   {represent}%
   {представлять}%

\ERitem%
   {residue}%
   {вычет}%

\ERitem%
   {resistance}%
   {сопротивление}%

\ERitem%
   {resistivity}%
   {удельное сопротивление}%

\ERitem%
   {resonance}%
   {резонанс}%

\ERitem%
   {restriction of correspondence $\Phi$ to set $C$}%
   {сужение соответствия $\Phi$ на множество $C$}%

\ERitem%
   {retract}%
   {ретракт}%

\ERitem%
   {retraction}%
   {ретракция}%

\ERitem%
   {right angle}%
   {прямой угол}%

\ERitem%
   {right side of equation}%
   {правая часть равенства}%

\ERitem%
   {right-distributive}%
   {дистрибутивен справа}%

\ERitem%
   {ring of characteristic $p$}%
   {кольцо характеристики $p$}%

\ERitem%
   {ring of integers}%
   {кольцо целых чисел}%

\ERitem%
   {rotation group}%
   {группа вращения}%

\ERitem%
   {row vector}%
   {вектор строка}%

\end{description}%

\section{S}%
\begin{description}%

\ERitem%
   {scale}%
   {масштаб}%

\ERitem%
   {scattering}%
   {рассеяние}%

\ERitem%
   {scattering amplitude}%
   {амплитуда рассеяния}%

\ERitem%
   {scattering from crystal}%
   {рассеяние на кристалле}%

\ERitem%
   {scattering theory}%
   {теория рассеяния}%

\ERitem%
   {secant}%
   {секанс}%

\ERitem%
   {section}%
   {сечение}%

\ERitem%
   {segment}%
   {отрезок}%

\ERitem%
   {seismology}%
   {сейсмология}%

\ERitem%
   {semiconductor}%
   {полупроводник}%

\ERitem%
   {semigroup}%
   {полугруппа}%

\ERitem%
   {series}%
   {ряд (бесконечная сумма)}%

\ERitem%
   {set of functions}%
   {семейство функций}%

\ERitem%
   {set of power of continuum}%
   {множество мощности континуум}%

\ERitem%
   {sieve of Eratosthenes}%
   {решето Эратосфена}%

\ERitem%
   {similar triangles}%
   {подобные треугольники}%

\ERitem%
   {simple polyvector}%
   {простой поливектор}%

\ERitem%
   {simple ring}%
   {простое кольцо}%

\ERitem%
   {simple root}%
   {простой корень}%

\ERitem%
   {simplex}%
   {симплекс}%

\ERitem%
   {simply connected}%
   {односвязный}%

\ERitem%
   {sine}%
   {синус}%

\ERitem%
   {single transitive representation}%
   {однотранзитивное представление}%

\ERitem%
   {skew field}%
   {тело (кольцо с делением)}%

\ERitem%
   {skew product of vectors}%
   {косое произведение векторов}%

\ERitem%
   {skew-symmetric form}%
   {кососимметричная форма}%

\ERitem%
   {skew-symmetric tensor}%
   {кососимметричный тензор}%

\ERitem%
   {solar eclipse}%
   {затмение Солнца}%

\ERitem%
   {solution of differential equation}%
   {решение дифференциального уравнения}%

\ERitem%
   {solve for the $c$}%
   {разрешить относительно $c$}%
\ERexample%
   {Equation may be solved for the $c$.}%
   {Уравнение можно разрешить  относительно $c$.}%
   {}{}{}{}

\ERitem%
   {spacelike vector}%
   {пространственноподобный вектор}%

\ERitem%
   {special relativity}%
   {специальная теория относительности}%

\ERitem%
   {sphere}%
   {сфера}%

\ERitem%
   {spherical triangle}%
   {сферический треугольник}%

\ERitem%
   {splittable algebra}%
   {расщепляемая алгебра}%

\ERitem%
   {square root}%
   {квадратный корень}%

\ERitem%
   {star\Hyph shaped domain}%
   {звёздная область}%

\ERitem%
   {stationary state}%
   {стационарное состояние}%

\ERitem%
   {strictly monotone function}%
   {строго монотонная функция}%

\ERitem%
   {strictly monotonic function}%
   {строго монотонная функция}%

\ERitem%
   {strongly monotone function}%
   {строго монотонная функция}%

\ERitem%
   {strongly monotonic function}%
   {строго монотонная функция}%

\ERitem%
   {structural constants}%
   {структурные константы}%

\ERitem%
   {structure constants}%
   {структурные константы}%

\ERitem%
   {subordinate}%
   {подчинённый}%
\ERexample%
   {For each open cover $U_a$ of $X$ there is a partition of unity $\{\varphi_b\}$ subordinate to the cover.}%
   {Для любого открытого покрытия $U_a$ многообразия $X$ существует разложение единицы $\{\varphi_b\}$, подчинённое покрытию.}%
   {}{}{}{}

\ERitem%
   {subtrahend}%
   {вычитаемое}%

\ERitem%
   {summation convention}%
   {правило суммирования}%

\ERitem%
   {summer solstice}%
   {летнее солнцестояние}%

\ERitem%
   {surjection}%
   {сюръекция}%

\ERitem%
   {symmetric}%
   {симметричный}%

\ERitem%
   {symmetry}%
   {симметрия}%

\ERitem%
   {synchronization procedure}%
   {процедура синхронизации}%

\ERitem%
   {synchrotron}%
   {синхротрон}%

\ERitem%
   {system of total differential equations}%
   {система уравнений в полных дифференциалах}%

\end{description}%

\section{T}%
\begin{description}%

\ERitem%
   {tail of vector}%
   {начало вектора}%

\ERitem%
   {tangent}%
   {тангенс}%

\ERitem%
   {Taylor series}%
   {Ряд Тейлора}%

\ERitem%
   {tensor of order $2$}%
   {тензор валентности $2$}%

\ERitem%
   {tensor product}%
   {тензорное произведение}%

\ERitem%
   {there exist}%
   {существует}%
\ERexample%
   {There exists a positive integer $N$ such that equations $F_1$, ..., $F_N$ are compatible.}%
   {Существует такое положительное целое число $N$, что уравнения $F_1$, ..., $F_N$ совместны.}%
   {}{}{}{}

\ERitem%
   {tidal acceleration}%
   {приливное ускорение}%

\ERitem%
   {timelike vector}%
   {времениподобный вектор}%

\ERitem%
   {to be coarser than}%
   {минорировать}%
\ERexample%
   {Filter $\mathfrak{F}$ is coarser than filter $\mathfrak{B}$.}%
   {Фильтр $\mathfrak{F}$ минорирует фильтр $\mathfrak{B}$.}%
   {}{}{}{}
\ERexample%
   {Topology $\mathfrak{T}_1$ is coarser than topology $\mathfrak{T}_2$.}%
   {Топология $\mathfrak{T}_1$ минорирует топологию $\mathfrak{T}_2$.}%
   {}{}{}{}

\ERitem%
   {to be finer than}%
   {мажорировать}%
\ERexample%
   {Filter $\mathfrak{F}$ is finer than filter $\mathfrak{B}$.}%
   {Фильтр $\mathfrak{F}$ мажорирует фильтр $\mathfrak{B}$.}%
   {}{}{}{}
\ERexample%
   {Topology $\mathfrak{T}_1$ is finer than topology $\mathfrak{T}_2$.}%
   {Топология $\mathfrak{T}_1$ мажорирует топологию $\mathfrak{T}_2$.}%
   {}{}{}{}

\ERitem%
   {topological space}%
   {топологическое пространство}%

\ERitem%
   {topology}%
   {топология}%

\ERitem%
   {torsion}%
   {кручение}%

\ERitem%
   {torus}%
   {тор}%

\ERitem%
   {total differential}%
   {полный дифференциал}%

\ERitem%
   {total ordering}%
   {линейная упорядоченность}%

\ERitem%
   {total space}%
   {тотальное пространство расслоения}%

\ERitem%
   {trajectory}%
   {траектория}%

\ERitem%
   {transition function}%
   {функция перехода}%

\ERitem%
   {transitive}%
   {транзитивный}%

\ERitem%
   {triangle}%
   {треугольник}%

\ERitem%
   {trigonometrical}%
   {тригонометрический}%

\ERitem%
   {trigonometry}%
   {тригонометрия}%

\ERitem%
   {trivial}%
   {тривиальный}%

\ERitem%
   {Tunguska Cosmic Body}%
   {тунгусский метеорит}%

\ERitem%
   {tuple}%
   {кортеж}%

\ERitem%
   {turbulence}%
   {турбулентность}%

\ERitem%
   {twin representation}%
   {парное представление}%

\ERitem%
   {two-sided ideal}%
   {двусторонний идеал}%

\end{description}%

\section{U}%
\begin{description}%

\ERitem%
   {uniform continuity}%
   {равномерная непрерывность}%

\ERitem%
   {uniform space}%
   {равномерное пространство}%

\ERitem%
   {uniformly continuous function}%
   {равномерно непрерывная функция}%

\ERitem%
   {unit sphere}%
   {единичная сфера}%

\ERitem%
   {unless otherwise stated}%
   {если не оговорено противное}%

\ERitem%
   {up to notation}%
   {с точностью до обозначений}%
\ERexample%
   {I can repeat, up to notation, proof of theorem 2.1.}%
   {Я могу с точностью до обозначений повторить доказательство теоремы 2.1.}%
   {}{}{}{}

\ERitem%
   {upper index}%
   {верхний индекс}%

\ERitem%
   {upstairs}%
   {вверху}%
\ERexample%
   {We sum over any index which appears twice in the same term, once upstairs and once downstairs.}%
   {Подразумевается сумма по любому индексу, появляющемуся дважды в одном и том же слагаемом, один раз вверху, другой - внизу.}%
   {}{}{}{}

\end{description}%

\section{V}%
\begin{description}%

\ERitem%
   {valued division ring}%
   {нормированное тело}%

\ERitem%
   {valued field}%
   {нормированное поле}%

\ERitem%
   {valued skew field}%
   {нормированное тело}%

\ERitem%
   {variable}%
   {переменная}%

\ERitem%
   {vector bundle}%
   {векторное расслоение}%

\ERitem%
   {vector function}%
   {вектор\hyph функция}%

\ERitem%
   {vector space}%
   {векторное пространство}%

\ERitem%
   {verify directly}%
   {непосредственная проверка доказывает}%
\ERexample%
   {We verify directly that $A$  is linear map.}%
   {Непосредственная проверка показывает, что $A$ - линейный оператор.}%
   {}{}{}{}
\ERexample%
   {We verify the statement of the theorem directly.}%
   {Мы можем доказать утверждение теоремы непосредственной проверкой.}%
   {}{}{}{}

\ERitem%
   {vertex}%
   {вершина}%

\ERitem%
   {vertical}%
   {вертикальный}%

\ERitem%
   {vertices}%
   {вершины}%

\ERitem%
   {viscosity}%
   {вязкость}%

\end{description}%

\section{W}%
\begin{description}%

\ERitem%
   {wave}%
   {волна}%

\ERitem%
   {winter solstice}%
   {зимнее солнцестояние}%

\ERitem%
   {without loss of generality}%
   {не нарушая общности; не уменьшая общности}%

\end{description}%

\section{Z}%
\begin{description}%

\ERitem%
   {Zermelo proposition}%
   {теорема Цермело}%

\ERitem%
   {zero divisor}%
   {делитель нуля}%

\end{description}%

\chapter{Русско английский словарь}%

\section{A}%
\begin{description}%

\REitem%
   {$A$\Hyph valued function}%
   {$A$\Hyph значная функция}%

\end{description}%

\section{D}%
\begin{description}%

\REitem%
   {$D$\Hyph vector space}%
   {$D$\Hyph векторное пространство}%

\end{description}%

\section{А}%
\begin{description}%

\REitem%
   {Abelian group}%
   {абелевая группа}%

\REitem%
   {absolute value}%
   {абсолютная величина}%

\REitem%
   {absolute convergence}%
   {абсолютная сходимость}%

\REitem%
   {auto parallel line}%
   {автопараллельная кривая}%

\REitem%
   {additive group}%
   {аддитивная группа}%

\REitem%
   {axiom of choice}%
   {аксиома выбора}%

\REitem%
   {Pasch's axiom}%
   {аксиома Паша}%

\REitem%
   {algebraic}%
   {алгебраический}%

\REitem%
   {algebraic complement of matrix; cofactor of matrix}%
   {алгебраическое дополнение матрицы}%

\REitem%
   {algebraic extension}%
   {алгебраическое расширение}%

\REitem%
   {amplitude}%
   {амплитуда}%

\REitem%
   {scattering amplitude}%
   {амплитуда рассеяния}%

\REitem%
   {analytic function}%
   {аналитическая функция}%

\REitem%
   {annihilator}%
   {аннулятор}%

\REitem%
   {apocentre}%
   {апоцентр}%

\REitem%
   {arccosecant}%
   {арккосеканс}%

\REitem%
   {arccosine}%
   {арккосинус}%

\REitem%
   {arccotangent}%
   {арккотангенс}%

\REitem%
   {arcsecant}%
   {арксеканс}%

\REitem%
   {arcsine}%
   {арксинус}%

\REitem%
   {arctangent}%
   {арктангенс}%

\REitem%
   {arity}%
   {арность}%

\REitem%
   {associativity}%
   {ассоциативность}%

\REitem%
   {associative}%
   {ассоциативный}%

\REitem%
   {attractor}%
   {аттрактор}%

\REitem%
   {metric-affine manifold}%
   {аффинно-метрическое многообразие}%

\end{description}%

\section{Б}%
\begin{description}%

\REitem%
   {base of fibered correspondence}%
   {база расслоенного соответствия}%

\REitem%
   {basis for vector space; basis of vector space}%
   {базис в векторном пространстве}%

\REitem%
   {basis for vector space; basis of vector space}%
   {базис векторного пространства}%

\REitem%
   {base of topology}%
   {базис топологии}%

\REitem%
   {filter base}%
   {базис фильтра}%

\REitem%
   {Banach algebra}%
   {банахова алгебра}%

\REitem%
   {Banach space}%
   {банахово пространство}%

\REitem%
   {massless particle}%
   {безмассовая частица}%

\REitem%
   {infinitesimal}%
   {бесконечно малая величина}%

\REitem%
   {infinitesimal}%
   {бесконечно малый}%

\REitem%
   {bijection}%
   {биекция}%

\REitem%
   {bimodule}%
   {бимодуль}%

\REitem%
   {binary}%
   {бинарный}%

\REitem%
   {Occam's razor}%
   {бритва Оккама}%

\REitem%
   {Klein bottle}%
   {бутылка Клейна}%

\end{description}%

\section{В}%
\begin{description}%

\REitem%
   {relevance}%
   {важность}%
\REexample%
   {Distinction between Lorentz and Poincar\'e groups is of no relevance here.}%
   {Различие между группами Лоренца и Пуанкаре для нас сейчас не важно.}%
   {}{}{}{}

\REitem%
   {upstairs}%
   {вверху}%
\REexample%
   {We sum over any index which appears twice in the same term, once upstairs and once downstairs.}%
   {Подразумевается сумма по любому индексу, появляющемуся дважды в одном и том же слагаемом, один раз вверху, другой - внизу.}%
   {}{}{}{}

\REitem%
   {column vector}%
   {вектор столбец}%

\REitem%
   {row vector}%
   {вектор строка}%

\REitem%
   {vector space}%
   {векторное пространство}%

\REitem%
   {vector bundle}%
   {векторное расслоение}%

\REitem%
   {vector function}%
   {вектор\hyph функция}%

\REitem%
   {probability}%
   {вероятность}%

\REitem%
   {vertical}%
   {вертикальный}%

\REitem%
   {upper index}%
   {верхний индекс}%

\REitem%
   {vertex}%
   {вершина}%

\REitem%
   {vertices}%
   {вершины}%

\REitem%
   {mutually orthogonal}%
   {взаимно ортогональные}%

\REitem%
   {mutually perpendicular}%
   {взаимно перпендикулярные}%

\REitem%
   {interaction}%
   {взаимодействие}%
\REexample%
   {Based on differential geometry, general relativity describes gravitational interaction.}%
   {Общая теория относительности описывает гравитационное взаимодействие, опираясь на дифференциальную геометрию.}%
   {}{}{}{}

\REitem%
   {exterior product}%
   {внешнее произведение}%

\REitem%
   {exterior differential}%
   {внешний дифференциал}%

\REitem%
   {external algebra}%
   {внешняя алгебра}%

\REitem%
   {external power}%
   {внешняя степень}%

\REitem%
   {downstairs}%
   {внизу}%
\REexample%
   {We sum over any index which appears twice in the same term, once upstairs and once downstairs.}%
   {Подразумевается сумма по любому индексу, появляющемуся дважды в одном и том же слагаемом, один раз вверху, другой - внизу.}%
   {}{}{}{}

\REitem%
   {perturbation}%
   {возмущение}%

\REitem%
   {wave}%
   {волна}%

\REitem%
   {generally speaking; in general}%
   {вообще говоря}%
\REexample%
   {However in general this product is not linear map.}%
   {Однако, вообще говоря, это отображением не является линейным.}%
   {}{}{}{}

\REitem%
   {completely integrable}%
   {вполне интегрируемый}%

\REitem%
   {timelike vector}%
   {времениподобный вектор}%

\REitem%
   {develop equation}%
   {вывести уравнение}%

\REitem%
   {convex function}%
   {выпуклая функция}%

\REitem%
   {residue}%
   {вычет}%

\REitem%
   {subtrahend}%
   {вычитаемое}%

\REitem%
   {viscosity}%
   {вязкость}%

\end{description}%

\section{Г}%
\begin{description}%

\REitem%
   {galaxy}%
   {галактика}%

\REitem%
   {geometry}%
   {геометрия}%

\REitem%
   {hyperplane}%
   {гиперплоскость}%

\REitem%
   {$G$-principal bundle}%
   {главное $G$-расслоение}%

\REitem%
   {principal bundle}%
   {главное расслоение}%

\REitem%
   {Global Positioning System}%
   {глобальная система позиционирования}%

\REitem%
   {holonomic coordinates}%
   {голономные координаты}%

\REitem%
   {homeomorphism}%
   {гомеоморфизм}%

\REitem%
   {homeomorphic}%
   {гомеоморфный}%

\REitem%
   {homology}%
   {гомология}%

\REitem%
   {homomorphism}%
   {гомоморфизм}%

\REitem%
   {homotopy}%
   {гомотопия}%

\REitem%
   {homotopic}%
   {гомотопный}%

\REitem%
   {event horizon}%
   {горизонт событий}%

\REitem%
   {gravity probe}%
   {гравитационный зонд}%

\REitem%
   {boundary}%
   {граница}%

\REitem%
   {boundary conditions}%
   {граничные условия}%

\REitem%
   {graph}%
   {граф}%

\REitem%
   {rotation group}%
   {группа вращения}%

\end{description}%

\section{Д}%
\begin{description}%

\REitem%
   {engine}%
   {двигатель}%

\REitem%
   {two-sided ideal}%
   {двусторонний идеал}%

\REitem%
   {operate}%
   {действовать}%
\REexample%
   {Operating on equation (1) with operator $V$ yields an integral equation.}%
   {Подействовав на уравнение (1) оператором $V$, получим интегральное уравнение.}%
   {}{}{}{}

\REitem%
   {Cartesian coordinate system}%
   {декартова система координат}%

\REitem%
   {Cartesian power}%
   {декартова степень}%

\REitem%
   {Cartesian product}%
   {декартово произведение}%

\REitem%
   {zero divisor}%
   {делитель нуля}%

\REitem%
   {diagram of correspondences}%
   {диаграмма соответствий}%

\REitem%
   {dynamics}%
   {динамика}%

\REitem%
   {discrete topology}%
   {дискретная топология}%

\REitem%
   {discrete space}%
   {дискретное пространство}%

\REitem%
   {left-distributive}%
   {дистрибутивен слева}%

\REitem%
   {right-distributive}%
   {дистрибутивен справа}%

\REitem%
   {distributive property of multiplication over addition}%
   {дистрибутивность умножения относительно сложения}%

\REitem%
   {diffraction}%
   {дифракция}%

\REitem%
   {the G\^ateaux differential}%
   {дифференциал Гато}%

\REitem%
   {the Fr\'echet differential}%
   {дифференциал Фреше}%

\REitem%
   {differentiate the function with respect to $x$}%
   {дифференцировать функцию по $x$}%

\REitem%
   {differentiable function}%
   {дифференцируемая функция}%

\REitem%
   {differentiability}%
   {дифференцируемость}%

\REitem%
   {differentiable in the G\^ateaux sense}%
   {дифференцируемый по Гато}%

\REitem%
   {differentiable in the Fr\'echet sense}%
   {дифференцируемый по Фреше}%

\REitem%
   {diffusion}%
   {диффузия}%

\REitem%
   {insulator}%
   {диэлектрик}%

\REitem%
   {proof by induction}%
   {доказательство по индукции}%

\REitem%
   {dual space}%
   {дуальное пространство}%

\REitem%
   {dual module}%
   {дуальный модуль}%

\end{description}%

\section{Е}%
\begin{description}%

\REitem%
   {unit sphere}%
   {единичная сфера}%

\REitem%
   {identity}%
   {единичный элемент}%

\REitem%
   {unless otherwise stated}%
   {если не оговорено противное}%

\REitem%
   {natural mapping}%
   {естественное отображение}%

\REitem%
   {natural morphism}%
   {естественный морфизм}%

\end{description}%

\section{З}%
\begin{description}%

\REitem%
   {dependence}%
   {зависимость}%

\REitem%
   {problem}%
   {задача}%

\REitem%
   {associative law}%
   {закон ассоциативности}%

\REitem%
   {distributive law}%
   {закон дистрибутивности}%

\REitem%
   {conservation law}%
   {закон сохранения}%

\REitem%
   {change of coordinates}%
   {замена координат}%

\REitem%
   {change of variable}%
   {замена переменной}%

\REitem%
   {closure of the set}%
   {замыкание множества}%

\REitem%
   {eclipse}%
   {затмение}%

\REitem%
   {solar eclipse}%
   {затмение Солнца}%

\REitem%
   {star\Hyph shaped domain}%
   {звёздная область}%

\REitem%
   {winter solstice}%
   {зимнее солнцестояние}%

\REitem%
   {denominator}%
   {знаменатель}%
\REexample%
   {Let us reduce items to a common denominator.}%
   {Приведём слагаемые к общему знаменателю.}%
   {}{}{}{}

\REitem%
   {relevance}%
   {значимость}%

\end{description}%

\section{И}%
\begin{description}%

\REitem%
   {emission of photon}%
   {излучение фотона}%

\REitem%
   {measurement}%
   {измерение}%

\REitem%
   {measure}%
   {измерять}%

\REitem%
   {insulator}%
   {изолятор}%

\REitem%
   {isotropic vector}%
   {изотропный вектор}%

\REitem%
   {has relevance to; is related to}%
   {имеет отношение к}%

\REitem%
   {name index}%
   {именной указатель}%

\REitem%
   {momentum}%
   {импульс}%

\REitem%
   {indicatrix}%
   {индикатриса}%

\REitem%
   {interference}%
   {интерференция}%

\REitem%
   {injection}%
   {инъекция}%

\end{description}%

\section{К}%
\begin{description}%

\REitem%
   {gauge invariance}%
   {калибровочная инвариантность}%

\REitem%
   {canonical map}%
   {каноническое отображение}%

\REitem%
   {category}%
   {категория}%

\REitem%
   {quadratic equation}%
   {квадратное уравнение}%

\REitem%
   {square root}%
   {квадратный корень}%

\REitem%
   {quasar}%
   {квазар}%

\REitem%
   {quasigroup}%
   {квазигруппа}%

\REitem%
   {quantum entanglement}%
   {квантовая запутанность}%

\REitem%
   {quark}%
   {кварк}%

\REitem%
   {quaternion}%
   {кватернион}%

\REitem%
   {kinematics}%
   {кинематика}%

\REitem%
   {equivalence class}%
   {класс эквивалентности}%

\REitem%
   {covariant}%
   {ковариантный}%

\REitem%
   {colloquium}%
   {коллоквиум}%

\REitem%
   {colloquia}%
   {коллоквиумы}%

\REitem%
   {ring of characteristic $p$}%
   {кольцо характеристики $p$}%

\REitem%
   {ring of integers}%
   {кольцо целых чисел}%

\REitem%
   {combinatorics}%
   {комбинаторика}%

\REitem%
   {commutative diagram}%
   {коммутативная диаграмма}%

\REitem%
   {commutativity}%
   {коммутативность}%

\REitem%
   {commutator}%
   {коммутатор}%

\REitem%
   {commute}%
   {коммутирует}%
\REexample%
   {Position and momentum operators do not commute.}%
   {Операторы положения и импульса не коммутируют.}%
   {}{}{}{}

\REitem%
   {compact\hyph open topology}%
   {компактно\hyph открытая топология}%

\REitem%
   {convection}%
   {конвекция}%

\REitem%
   {congruence}%
   {конгруэнтность}%

\REitem%
   {head of vector}%
   {конец вектора}%

\REitem%
   {finite set}%
   {конечное множество}%

\REitem%
   {finite dimensional}%
   {конечномерный}%

\REitem%
   {contravariant}%
   {контравариантный}%

\REitem%
   {conformal transformation}%
   {конформное преобразование}%

\REitem%
   {coordinate chart}%
   {координатная карта}%

\REitem%
   {correlation}%
   {корреляция}%

\REitem%
   {tuple}%
   {кортеж}%

\REitem%
   {cosecant}%
   {косеканс}%

\REitem%
   {cosine}%
   {косинус}%

\REitem%
   {skew product of vectors}%
   {косое произведение векторов}%

\REitem%
   {skew-symmetric form}%
   {кососимметричная форма}%

\REitem%
   {skew-symmetric tensor}%
   {кососимметричный тензор}%

\REitem%
   {cotangent}%
   {котангенс}%

\REitem%
   {connection coefficients}%
   {коэффициенты связности}%

\REitem%
   {multiplicity of $x$ in $f$}%
   {кратность $x$ в $f$}%
\REexample%
   {If the multiplicity of $a$ is greater then $1$, $a$ is called a multiple root.}%
   {Если кратность $a$ больше, чем $1$, то $a$ называется кратным корнем.}%
   {}{}{}{}

\REitem%
   {multiple root; repeated root}%
   {кратный корень}%

\REitem%
   {curve}%
   {кривая}%

\REitem%
   {curvilinear coordinates}%
   {криволинейные координаты}%

\REitem%
   {crystal lattice}%
   {кристаллическая решётка}%

\REitem%
   {torsion}%
   {кручение}%

\end{description}%

\section{Л}%
\begin{description}%

\REitem%
   {Lagrangian}%
   {лагранжиан}%

\REitem%
   {left side of equation}%
   {левая часть равенства}%

\REitem%
   {summer solstice}%
   {летнее солнцестояние}%

\REitem%
   {total ordering}%
   {линейная упорядоченность}%

\REitem%
   {linearly dependent}%
   {линейно зависимые}%

\REitem%
   {linearly independent}%
   {линейно независимые}%

\REitem%
   {Moebius band}%
   {лист Мёбиуса}%

\REitem%
   {lift of vector field}%
   {лифт векторного поля}%

\REitem%
   {lift of morphism}%
   {лифт морфизма}%

\REitem%
   {lift of correspondence}%
   {лифт соответствия}%

\REitem%
   {locally compact space}%
   {локально компактное пространство}%

\REitem%
   {loop (quasigroup with unit element)}%
   {лупа (квазигруппа с единицей)}%

\end{description}%

\section{М}%
\begin{description}%

\REitem%
   {to be finer than}%
   {мажорировать}%
\REexample%
   {Filter $\mathfrak{F}$ is finer than filter $\mathfrak{B}$.}%
   {Фильтр $\mathfrak{F}$ мажорирует фильтр $\mathfrak{B}$.}%
   {}{}{}{}
\REexample%
   {Topology $\mathfrak{T}_1$ is finer than topology $\mathfrak{T}_2$.}%
   {Топология $\mathfrak{T}_1$ мажорирует топологию $\mathfrak{T}_2$.}%
   {}{}{}{}

\REitem%
   {little group}%
   {малая группа}%

\REitem%
   {mass}%
   {масса}%

\REitem%
   {massive particle}%
   {массивная частица}%

\REitem%
   {scale}%
   {масштаб}%

\REitem%
   {mathematician}%
   {математик}%

\REitem%
   {mathematics}%
   {математика}%

\REitem%
   {mathematical}%
   {математический}%

\REitem%
   {Jacobian matrix}%
   {матрица Якоби}%

\REitem%
   {Kerr metric}%
   {метрика Керра}%

\REitem%
   {to be coarser than}%
   {минорировать}%
\REexample%
   {Filter $\mathfrak{F}$ is coarser than filter $\mathfrak{B}$.}%
   {Фильтр $\mathfrak{F}$ минорирует фильтр $\mathfrak{B}$.}%
   {}{}{}{}
\REexample%
   {Topology $\mathfrak{T}_1$ is coarser than topology $\mathfrak{T}_2$.}%
   {Топология $\mathfrak{T}_1$ минорирует топологию $\mathfrak{T}_2$.}%
   {}{}{}{}

\REitem%
   {Milky Way}%
   {Млечный Путь}%

\REitem%
   {polynomial}%
   {многочлен}%

\REitem%
   {range}%
   {множество значений}%

\REitem%
   {set of power of continuum}%
   {множество мощности континуум}%

\REitem%
   {factor}%
   {множитель}%

\REitem%
   {angular momentum}%
   {момент количества движения}%

\REitem%
   {monotone function; monotonic function}%
   {монотонная функция}%

\REitem%
   {identification morphism}%
   {морфизм отождествления}%

\REitem%
   {fibered map}%
   {морфизм расслоений}%

\REitem%
   {power of set}%
   {мощность множества}%

\REitem%
   {multiplicative group}%
   {мультипликативная группа}%

\REitem%
   {muon}%
   {мюон}%

\end{description}%

\section{Н}%
\begin{description}%

\REitem%
   {at first glance; at first sight}%
   {на первый взгляд}%

\REitem%
   {highest common factor of $p$ and $q$}%
   {наибольший общий делитель $p$ и $q$}%

\REitem%
   {covering space}%
   {накрытие}%
\REexample%
   {Consider the covering space $R\rightarrow S^1$ of the circle $S^1$ defined by $p(t)=(\sin t,\cos t)$ for any $t\in R$.}%
   {Рассмотрим накрытие $R\rightarrow S^1$ окружности $S^1$, определённое формулой $p(t)=(\sin t,\cos t)$ для любого $t\in R$.}%
   {}{}{}{}

\REitem%
   {positive integer}%
   {натуральное число}%

\REitem%
   {tail of vector}%
   {начало вектора}%

\REitem%
   {without loss of generality}%
   {не нарушая общности}%

\REitem%
   {without loss of generality}%
   {не уменьшая общности}%

\REitem%
   {non-Abelian group}%
   {неабелевая группа}%

\REitem%
   {nondegenerate form}%
   {невырожденная форма}%

\REitem%
   {anholonomity}%
   {неголономность}%

\REitem%
   {anholonomic coordinates}%
   {неголономные координаты}%

\REitem%
   {neutrino}%
   {нейтрино}%

\REitem%
   {neutron}%
   {нейтрон}%

\REitem%
   {neutron star}%
   {нейтронная звезда}%

\REitem%
   {necessary and sufficient}%
   {необходимо и достаточно}%
\REexample%
   {In order that any $x\in A$ is root of the system of linear equations\[a_i^jx^i=0\]necessary and sufficient $a^{\gij}_{\gii}=0$.}%
   {Для того, чтобы любое $x\in A$ было корнем системы линейных уравнений\[a_i^jx^i=0\]необходимо и достаточно, чтобы $a^{\gij}_{\gii}=0$.}%
   {}{}{}{}

\REitem%
   {inhomogeneous Lorentz group}%
   {неоднородная группа Лоренца}%

\REitem%
   {inhomogeneous}%
   {неоднородный}%

\REitem%
   {verify directly}%
   {непосредственная проверка доказывает}%
\REexample%
   {We verify directly that $A$  is linear map.}%
   {Непосредственная проверка показывает, что $A$ - линейный оператор.}%
   {}{}{}{}
\REexample%
   {We verify the statement of the theorem directly.}%
   {Мы можем доказать утверждение теоремы непосредственной проверкой.}%
   {}{}{}{}

\REitem%
   {continuous in neighborhood}%
   {непрерывен в окрестности}%

\REitem%
   {continuous in $x$}%
   {непрерывный по $x$}%

\REitem%
   {irreducible representation}%
   {неприводимое представление}%

\REitem%
   {inequation}%
   {неравенство}%

\REitem%
   {nontrivial}%
   {нетривиальный}%

\REitem%
   {lower index}%
   {нижний индекс}%

\REitem%
   {absolute value; norm}%
   {норма}%
\REexample%
   {Norm on $\Omega$\Hyph group $A$ is a map\[d\in A\rightarrow \|d\|\in R\]such that\begin{itemize}\item $\|a\|\ge 0$\item $\|a\|=0$ iff, $a=0$\item $\|a+b\|\le\|a\|+\|b\|$\item $\|-a\|=\|a\|$\end{itemize}}%
   {Норма на $\Omega$\Hyph группе $A$ - это отображение\[d\in A\rightarrow \|d\|\in R\]такое, что\begin{itemize}\item $\|a\|\ge 0$\item $\|a\|=0$ равносильно $a=0$\item $\|a+b\|\le\|a\|+\|b\|$\item $\|-a\|=\|a\|$\end{itemize}}%
   {}{}{}{}
\REexample%
   {Absolute value on skew field $D$ is a map\[d\in D\rightarrow |d|\in R\]which satisfies the following axioms\begin{itemize}\item $|a|\ge 0$\item $|a|=0$ if, and only if, $a=0$\item $|ab|=|a|\ |b|$\item $|a+b|\le|a|+|b|$\end{itemize}}%
   {Норма на теле $D$ - это отображение\[d\in D\rightarrow |d|\in R\]такое, что\begin{itemize}\item $|a|\ge 0$\item $|a|=0$ равносильно $a=0$\item $|ab|=|a|\ |b|$\item $|a+b|\le|a|+|b|$\end{itemize}}%
   {}{}{}{}

\REitem%
   {valued field}%
   {нормированное поле}%

\REitem%
   {normed space}%
   {нормированное пространство}%

\REitem%
   {valued division ring; valued skew field}%
   {нормированное тело}%

\end{description}%

\section{О}%
\begin{description}%

\REitem%
   {enveloping algebra}%
   {обвёртывающая алгебра}%

\REitem%
   {range}%
   {область значений}%

\REitem%
   {domain}%
   {область определения}%

\REitem%
   {image under map}%
   {образ при отображении}%
\REexample%
   {We define the image of the set $A$ under correspondence $\Phi$ according to law\[A\Phi=\{b:(a,b)\in\Phi,a\in A\}\]}%
   {Мы определим образ множества $A$ при соответствии $\Phi$ согласно равенству\[A\Phi=\{b:(a,b)\in\Phi,a\in A\}\]}%
   {}{}{}{}

\REitem%
   {generator}%
   {образующая}%

\REitem%
   {conversely}%
   {обратно}%

\REitem%
   {inverse transformation}%
   {обратное преобразование}%

\REitem%
   {reciprocal image}%
   {обратный образ}%

\REitem%
   {poolback bundle}%
   {обратный образ расслоения}%

\REitem%
   {Hadamard inverse}%
   {обращение Адамара}%

\REitem%
   {general relativity}%
   {общая теория относительности}%

\REitem%
   {anholonomity object}%
   {объект неголономности}%

\REitem%
   {envelope of a family of plane curves}%
   {огибающая семейства плоских кривых}%

\REitem%
   {homogeneous Lorentz group}%
   {однородная группа Лоренца}%

\REitem%
   {homogeneous}%
   {однородный}%

\REitem%
   {simply connected}%
   {односвязный}%

\REitem%
   {single transitive representation}%
   {однотранзитивное представление}%

\REitem%
   {monomial}%
   {одночлен}%

\REitem%
   {neighborhood}%
   {окрестность}%

\REitem%
   {closure operator}%
   {оператор замыкания}%

\REitem%
   {momentum operator}%
   {оператор импульса}%

\REitem%
   {creation operator}%
   {оператор рождения}%

\REitem%
   {annihilation operator}%
   {оператор уничтожения}%

\REitem%
   {define}%
   {определяет}%
\REexample%
   {This equation defines the inverse transformation.}%
   {Это уравнение определяет обратное преобразование.}%
   {}{}{}{}

\REitem%
   {mutually orthogonal}%
   {ортогональные друг другу}%

\REitem%
   {orthonormal basis}%
   {ортонормированный базис}%

\REitem%
   {remainder}%
   {остаток}%

\REitem%
   {acute angle}%
   {острый угол}%

\REitem%
   {relation}%
   {отношение}%

\REitem%
   {mapping}%
   {отображение}%

\REitem%
   {identification}%
   {отождествление}%

\REitem%
   {segment}%
   {отрезок}%

\REitem%
   {evidently}%
   {очевидно}%
\REexample%
   {Evidently $x=1$ is the root of the equation.}%
   {Очевидно, $x=1$ является корнем уравнения.}%
   {}{}{}{}

\REitem%
   {it is evident that}%
   {очевидно, что}%
\REexample%
   {From (2.2), it is evident that any solution of (2.7) satisfies (2.9).}%
   {На основании (2.2) очевидно, что любое решение уравнения (2.7) удовлетворяет (2.9).}%
   {}{}{}{}

\REitem%
   {evidence}%
   {очевидность}%

\end{description}%

\section{П}%
\begin{description}%

\REitem%
   {parallelepiped}%
   {параллелепипед}%

\REitem%
   {parallel transport}%
   {параллельный перенос}%

\REitem%
   {twin representation}%
   {парное представление}%

\REitem%
   {variable}%
   {переменная}%

\REitem%
   {Bott periodicity}%
   {периодичность Ботта}%

\REitem%
   {pericentre}%
   {перицентр}%

\REitem%
   {mutually perpendicular}%
   {перпендикулярные друг другу}%

\REitem%
   {at least}%
   {по крайней мере}%
\REexample%
   {At least in the neighborhood of the identity.}%
   {По крайней мере, в окрестности единичного элемента.}%
   {}{}{}{}

\REitem%
   {behavior}%
   {поведение}%

\REitem%
   {absorption of photon}%
   {поглощение фотона}%

\REitem%
   {similar triangles}%
   {подобные треугольники}%

\REitem%
   {in a similar way}%
   {подобным образом}%
\REexample%
   {In a similar way, we can introduce a coordinate reference frame.}%
   {Подобным образом мы можем определить координатную систему отсчёта.}%
   {}{}{}{}

\REitem%
   {evaluating by equating $x$ to the $a$}%
   {подстановка $a$ вместо $x$}%

\REitem%
   {subordinate}%
   {подчинённый}%
\REexample%
   {For each open cover $U_a$ of $X$ there is a partition of unity $\{\varphi_b\}$ subordinate to the cover.}%
   {Для любого открытого покрытия $U_a$ многообразия $X$ существует разложение единицы $\{\varphi_b\}$, подчинённое покрытию.}%
   {}{}{}{}

\REitem%
   {lift of vector field}%
   {подъём векторного поля}%

\REitem%
   {lift of morphism}%
   {подъём морфизма}%

\REitem%
   {lift of correspondence}%
   {подъём соответствия}%

\REitem%
   {componentwise}%
   {покомпонентно}%

\REitem%
   {cover}%
   {покрытие}%

\REitem%
   {complex field}%
   {поле комплексных чисел}%

\REitem%
   {rational field}%
   {поле рациональных чисел}%

\REitem%
   {polyadditive map}%
   {полиаддитивное отображение}%

\REitem%
   {polyvector}%
   {поливектор}%

\REitem%
   {polylinear form}%
   {полилинейная форма}%

\REitem%
   {polynomial}%
   {полином}%

\REitem%
   {complete system}%
   {полная система}%

\REitem%
   {complete lattice}%
   {полная структура}%

\REitem%
   {complete field}%
   {полное поле}%

\REitem%
   {complete space}%
   {полное пространство}%

\REitem%
   {complete division ring}%
   {полное тело}%

\REitem%
   {total differential}%
   {полный дифференциал}%

\REitem%
   {positive definite form}%
   {положительно определённая форма}%

\REitem%
   {polology}%
   {полология}%

\REitem%
   {semigroup}%
   {полугруппа}%

\REitem%
   {semiconductor}%
   {полупроводник}%

\REitem%
   {obtain by differentiating}%
   {получить дифференцированием}%

\REitem%
   {mutually orthogonal}%
   {попарно ортогональные}%

\REitem%
   {generated}%
   {порождённый}%
\REexample%
   {Algebra $A$ generated by the set $S$ is a $K$\Hyph algebra}%
   {Алгебра $A$, порождённая множеством $S$, является $K$\Hyph алгеброй.}%
   {}{}{}{}

\REitem%
   {Cauchy sequence}%
   {последовательность Коши}%

\REitem%
   {right side of equation}%
   {правая часть равенства}%

\REitem%
   {chain rule}%
   {правило дифференцирования сложной функции}%

\REitem%
   {Cramer's Rule}%
   {правило Крамера}%

\REitem%
   {L'H\^ospital's rule}%
   {правило Лопиталя}%

\REitem%
   {summation convention}%
   {правило суммирования}%

\REitem%
   {limit}%
   {предел}%

\REitem%
   {limit of sequence}%
   {предел последовательности}%

\REitem%
   {limit of correspondence with respect to the filter}%
   {предел соответствия по фильтру}%

\REitem%
   {limit point}%
   {предельная точка}%

\REitem%
   {limit set}%
   {предельное множество}%

\REitem%
   {passage to the limit}%
   {предельный переход}%

\REitem%
   {interaction picture}%
   {представление взаимодействия}%

\REitem%
   {represent}%
   {представлять}%

\REitem%
   {preordering}%
   {предупорядоченность}%

\REitem%
   {Lorentz transformation}%
   {преобразование Лоренца}%

\REitem%
   {precession}%
   {прецессия}%

\REitem%
   {approximation}%
   {приближение}%

\REitem%
   {reduction of similar terms}%
   {приведение подобных}%

\REitem%
   {reduced Cartesian product of bundles}%
   {приведенное декартово произведение расслоений}%

\REitem%
   {reduced fibered correspondence}%
   {приведенное расслоенное соответствие}%

\REitem%
   {monic polynomial}%
   {приведенный многочлен}%

\REitem%
   {tidal acceleration}%
   {приливное ускорение}%

\REitem%
   {the Mach principle}%
   {принцип Маха}%

\REitem%
   {adjoint group}%
   {присоединённая группа}%

\REitem%
   {adjoin}%
   {присоединить}%
\REexample%
   {To derive equations of motion of a charged particle we adjoin Lorentz equations to Maxwell equations.}%
   {Для того, чтобы вывести уравнения движения заряда, мы присоединим уравнения Лоренца к уравнения Максвелла.}%
   {}{}{}{}

\REitem%
   {causal vector field}%
   {причинное векторное поле}%

\REitem%
   {causal scalar field}%
   {причинное скалярное поле}%

\REitem%
   {causal relationship}%
   {причинно-следственная связь}%

\REitem%
   {proceeding in this way}%
   {продолжая таким образом}%

\REitem%
   {extension of correspondence}%
   {продолжение соответствия}%

\REitem%
   {proceeding in this way}%
   {продолжив этот процесс}%

\REitem%
   {projective plane}%
   {проективная плоскость}%

\REitem%
   {projection}%
   {проекция}%

\REitem%
   {derivative of second or greater order with respect}%
   {производная второго или более высокого порядка по }%

\REitem%
   {the G\^ateaux derivative}%
   {производная Гато}%

\REitem%
   {the Fr\'echet derivative}%
   {производная Фреше}%

\REitem%
   {poolback bundle}%
   {прообраз расслоения}%

\REitem%
   {simple ring}%
   {простое кольцо}%

\REitem%
   {prime ideal}%
   {простой идеал}%

\REitem%
   {simple root}%
   {простой корень}%

\REitem%
   {simple polyvector}%
   {простой поливектор}%

\REitem%
   {spacelike vector}%
   {пространственноподобный вектор}%

\REitem%
   {manifolds with affine connections}%
   {пространство аффинной связности}%

\REitem%
   {event space}%
   {пространство событий}%

\REitem%
   {opposite preordering}%
   {противоположная предупорядоченность}%

\REitem%
   {contradiction}%
   {противоречие}%
\REexample%
   {The contradiction completes the proof of the theorem.}%
   {Полученное противоречие доказывает теорему.}%
   {}{}{}{}

\REitem%
   {synchronization procedure}%
   {процедура синхронизации}%

\REitem%
   {Gram-Schmidt orthogonalization procedure}%
   {процесс ортогонализации Грама--Шмидта}%

\REitem%
   {right angle}%
   {прямой угол}%

\REitem%
   {pseudo\Hyph Euclidean space}%
   {псевдоевклидовое пространство}%

\REitem%
   {pulsar}%
   {пульсар}%

\REitem%
   {pfaffian derivative}%
   {пфаффова производная}%

\end{description}%

\section{Р}%
\begin{description}%

\REitem%
   {uniform continuity}%
   {равномерная непрерывность}%

\REitem%
   {uniformly continuous function}%
   {равномерно непрерывная функция}%

\REitem%
   {uniform space}%
   {равномерное пространство}%

\REitem%
   {radiation belt}%
   {радиационный пояс}%

\REitem%
   {partition of unity}%
   {разбиение единицы}%

\REitem%
   {partition of unity}%
   {разложение единицы}%

\REitem%
   {factorization}%
   {разложение на множители}%

\REitem%
   {decomposition of map}%
   {разложение отображения}%

\REitem%
   {factor}%
   {разложить на множители}%
\REexample%
   {To factor a polynomial means to find two or more po\-ly\-no\-mials whose product is the given polynomial.}%
   {Чтобы разложить многочлен на множители, необходимо найти два или более многочленов, произведение которых есть данный многочлен.}%
   {}{}{}{}

\REitem%
   {difference}%
   {разность}%

\REitem%
   {solve for the $c$}%
   {разрешить относительно $c$}%
\REexample%
   {Equation may be solved for the $c$.}%
   {Уравнение можно разрешить  относительно $c$.}%
   {}{}{}{}

\REitem%
   {propagation}%
   {распространение}%

\REitem%
   {scattering}%
   {рассеяние}%

\REitem%
   {scattering from crystal}%
   {рассеяние на кристалле}%

\REitem%
   {algebra bundle}%
   {расслоенная алгебра}%

\REitem%
   {group bundle}%
   {расслоенная группа}%

\REitem%
   {fibered product}%
   {расслоенное произведение}%

\REitem%
   {fibered correspondence}%
   {расслоенное соответствие}%

\REitem%
   {consider}%
   {рассматривать}%
\REexample%
   {Consider correspondence from set $A$ to set $B$.}%
   {Рассмотрим соответствие $\Phi$ из множества $A$ в множество $B$.}%
   {}{}{}{}

\REitem%
   {extension field}%
   {расширение поля}%

\REitem%
   {enhanced}%
   {расширенный}%

\REitem%
   {splittable algebra}%
   {расщепляемая алгебра}%

\REitem%
   {regression}%
   {регрессия}%

\REitem%
   {resonance}%
   {резонанс}%

\REitem%
   {retract}%
   {ретракт}%

\REitem%
   {retraction}%
   {ретракция}%

\REitem%
   {reflexive}%
   {рефлексивный}%

\REitem%
   {solution of differential equation}%
   {решение дифференциального уравнения}%

\REitem%
   {sieve of Eratosthenes}%
   {решето Эратосфена}%

\REitem%
   {series}%
   {ряд (бесконечная сумма)}%

\REitem%
   {Taylor series}%
   {Ряд Тейлора}%

\end{description}%

\section{С}%
\begin{description}%

\REitem%
   {up to notation}%
   {с точностью до обозначений}%
\REexample%
   {I can repeat, up to notation, proof of theorem 2.1.}%
   {Я могу с точностью до обозначений повторить доказательство теоремы 2.1.}%
   {}{}{}{}

\REitem%
   {finest topology}%
   {самая сильная топология}%

\REitem%
   {coarsest topology}%
   {самая слабая топология}%

\REitem%
   {hyperfine splitting}%
   {сверхтонкое расщепление}%

\REitem%
   {free representation}%
   {свободное представление}%

\REitem%
   {connected group}%
   {связная группа}%

\REitem%
   {seismology}%
   {сейсмология}%

\REitem%
   {secant}%
   {секанс}%

\REitem%
   {set of functions}%
   {семейство функций}%

\REitem%
   {section}%
   {сечение}%

\REitem%
   {force}%
   {сила}%

\REitem%
   {friction}%
   {сила трения}%

\REitem%
   {symmetric}%
   {симметричный}%

\REitem%
   {symmetry}%
   {симметрия}%

\REitem%
   {simplex}%
   {симплекс}%

\REitem%
   {sine}%
   {синус}%

\REitem%
   {synchrotron}%
   {синхротрон}%

\REitem%
   {closure system}%
   {система замыканий}%

\REitem%
   {reference frame}%
   {система отсчёта}%

\REitem%
   {system of total differential equations}%
   {система уравнений в полных дифференциалах}%

\REitem%
   {as small as we please}%
   {сколь угодно малый}%

\REitem%
   {fiber}%
   {слой}%

\REitem%
   {adjacent angle}%
   {смежный угол}%

\REitem%
   {mixed system}%
   {смешанная система}%

\REitem%
   {proper value}%
   {собственное значение}%

\REitem%
   {proper state}%
   {собственное состояние}%

\REitem%
   {according to theorem 2.1; By Theorem 2,1}%
   {согласно теореме 2.1}%
\REexample%
   {According to theorem 2.1, triangles $ABC$ and $DBC$ equal.}%
   {Согласно теореме 2.1 треугольники $ABC$ и $DBC$ равны.}%
   {}{}{}{}
\REexample%
   {By Theorem 2.1, $a=b$.}%
   {Согласно теореме 2.1, $a=b$.}%
   {}{}{}{}

\REitem%
   {congruence}%
   {согласованность}%

\REitem%
   {convention}%
   {соглашение}%
\REexample%
   {We use the convention that we present any set of vectors of the vector space as a row.}%
   {Мы пользуемся соглашением, что в заданном векторном пространстве мы представляем любое семейство векторов в виде строки.}%
   {}{}{}{}

\REitem%
   {factor}%
   {сомножитель}%

\REitem%
   {correspondence from $A$ to $B$}%
   {соответствие из $A$ в $B$}%

\REitem%
   {resistance}%
   {сопротивление}%

\REitem%
   {conjugate quaternion}%
   {сопряжённый кватернион}%

\REitem%
   {special relativity}%
   {специальная теория относительности}%

\REitem%
   {helical structure}%
   {спиральная структура}%

\REitem%
   {helicity}%
   {спиральность}%

\REitem%
   {comparable topology}%
   {сравнимые топологии}%

\REitem%
   {stationary state}%
   {стационарное состояние}%

\REitem%
   {degree of map}%
   {степень отображения}%

\REitem%
   {strictly monotone function; strictly monotonic function; strongly monotone function; strongly monotonic function}%
   {строго монотонная функция}%

\REitem%
   {lattice}%
   {структура (алгебраическая система)}%

\REitem%
   {structural constants; structure constants}%
   {структурные константы}%

\REitem%
   {restriction of correspondence $\Phi$ to set $C$}%
   {сужение соответствия $\Phi$ на множество $C$}%

\REitem%
   {relevance}%
   {существенность}%

\REitem%
   {essential parameters in a set of functions}%
   {существенные параметры семейства функций}%

\REitem%
   {there exist}%
   {существует}%
\REexample%
   {There exists a positive integer $N$ such that equations $F_1$, ..., $F_N$ are compatible.}%
   {Существует такое положительное целое число $N$, что уравнения $F_1$, ..., $F_N$ совместны.}%
   {}{}{}{}

\REitem%
   {sphere}%
   {сфера}%

\REitem%
   {spherical triangle}%
   {сферический треугольник}%

\REitem%
   {converge}%
   {сходиться}%
\REexample%
   {Filter $\mathfrak{F}$ converges to $x$.}%
   {Фильтр $\mathfrak{F}$ сходится к $x$.}%
   {}{}{}{}

\REitem%
   {countable set}%
   {счётное множество}%

\REitem%
   {counter}%
   {счётчик}%

\REitem%
   {surjection}%
   {сюръекция}%

\end{description}%

\section{Т}%
\begin{description}%

\REitem%
   {multiplication table}%
   {таблица умножения}%

\REitem%
   {tangent}%
   {тангенс}%

\REitem%
   {division ring; skew field}%
   {тело (кольцо с делением)}%

\REitem%
   {tensor of order $2$}%
   {тензор валентности $2$}%

\REitem%
   {field-strength tensor}%
   {тензор напряжённости поля}%

\REitem%
   {tensor product}%
   {тензорное произведение}%

\REitem%
   {Bell's theorem}%
   {теорема Белла}%

\REitem%
   {mean value theorem}%
   {теорема о конечных приращениях}%

\REitem%
   {Zermelo proposition}%
   {теорема Цермело}%

\REitem%
   {graph theory}%
   {теория графов}%

\REitem%
   {scattering theory}%
   {теория рассеяния}%

\REitem%
   {iff}%
   {тогда и только тогда, когда}%
\REexample%
   {$a=0$ iff $a_{\gi i}^{\gi j}=0$ for any ${\gi i}$, ${\gi j}$.}%
   {$a=0$ тогда и только тогда, когда $a_{\gi i}^{\gi j}=0$ для любых ${\gi i}$, ${\gi j}$.}%
   {}{}{}{}

\REitem%
   {identical particles}%
   {тождественные частицы}%

\REitem%
   {topological space}%
   {топологическое пространство}%

\REitem%
   {topology}%
   {топология}%

\REitem%
   {torus}%
   {тор}%

\REitem%
   {total space}%
   {тотальное пространство расслоения}%

\REitem%
   {point}%
   {точечный}%

\REitem%
   {point}%
   {точка}%

\REitem%
   {cluster point}%
   {точка прикосновения}%

\REitem%
   {exact sequence of modules}%
   {точная последовательность модулей}%

\REitem%
   {trajectory}%
   {траектория}%

\REitem%
   {transitive}%
   {транзитивный}%

\REitem%
   {friction}%
   {трение}%

\REitem%
   {triangle}%
   {треугольник}%

\REitem%
   {chart over $U$}%
   {тривиализация над $U$}%

\REitem%
   {trivial}%
   {тривиальный}%

\REitem%
   {trigonometrical}%
   {тригонометрический}%

\REitem%
   {trigonometry}%
   {тригонометрия}%

\REitem%
   {Tunguska Cosmic Body}%
   {тунгусский метеорит}%

\REitem%
   {obtuse angle}%
   {тупой угол}%

\REitem%
   {turbulence}%
   {турбулентность}%

\end{description}%

\section{У}%
\begin{description}%

\REitem%
   {angle}%
   {угол}%

\REitem%
   {angle of reflection}%
   {угол отражения}%

\REitem%
   {angle of incidence}%
   {угол падения}%

\REitem%
   {angle of refraction}%
   {угол преломления}%

\REitem%
   {resistivity}%
   {удельное сопротивление}%

\REitem%
   {knot}%
   {узел}%

\REitem%
   {multiplication}%
   {умножение}%

\REitem%
   {multiply by $2$}%
   {умножить на $2$}%

\REitem%
   {multiply by $b$}%
   {умножить на $b$}%

\REitem%
   {monic polynomial}%
   {унитарный многочлен}%

\REitem%
   {ordered set}%
   {упорядоченное множество}%

\REitem%
   {ordering}%
   {упорядоченность}%

\REitem%
   {partial differential equation}%
   {уравнение в частных производных}%

\REitem%
   {equation is satisfied identically}%
   {уравнение удовлетворяется тождественно}%

\REitem%
   {acceleration}%
   {ускорение}%

\REitem%
   {accelerator}%
   {ускоритель}%

\REitem%
   {conditions of integrability}%
   {условия интегрируемости}%

\end{description}%

\section{Ф}%
\begin{description}%

\REitem%
   {quotient set}%
   {фактор множество}%

\REitem%
   {difference module}%
   {фактор модуль}%

\REitem%
   {quotient bundle}%
   {фактор расслоение}%

\REitem%
   {factor group; quotient group}%
   {факторгруппа}%

\REitem%
   {quotient ring}%
   {факторкольцо}%

\REitem%
   {quotient topology}%
   {фактортопология}%

\REitem%
   {physicist}%
   {физик}%

\REitem%
   {physics}%
   {физика}%

\REitem%
   {physical}%
   {физический}%

\REitem%
   {filter}%
   {фильтр}%

\REitem%
   {Finsler metric; Finslerian metric}%
   {финслеровая метрика}%

\REitem%
   {photon}%
   {фотон}%

\REitem%
   {fundamental sequence}%
   {фундаментальная последовательность}%

\REitem%
   {functor}%
   {функтор}%

\REitem%
   {gluing functions}%
   {функции склеивания}%

\REitem%
   {functional}%
   {функционал}%

\REitem%
   {function $f$ of $x$}%
   {функция $f$ от $x$}%

\REitem%
   {transition function}%
   {функция перехода}%

\end{description}%

\section{Х}%
\begin{description}%

\REitem%
   {chaos}%
   {хаос}%

\end{description}%

\section{Ц}%
\begin{description}%

\REitem%
   {cycle}%
   {цикл}%

\REitem%
   {cyclic group}%
   {циклическая группа}%

\end{description}%

\section{Ч}%
\begin{description}%

\REitem%
   {partial ordering}%
   {частичная упорядоченность}%

\REitem%
   {frequency}%
   {частота}%

\REitem%
   {parity}%
   {чётность}%

\REitem%
   {numerator}%
   {числитель}%

\REitem%
   {real function; real valued function}%
   {числовая функция}%

\end{description}%

\section{Э}%
\begin{description}%

\REitem%
   {Euclidean metric}%
   {эвклидова метрика}%

\REitem%
   {Euclidean space}%
   {эвклидово пространство}%

\REitem%
   {equivalence relation}%
   {эквивалентность}%

\REitem%
   {extreme line}%
   {экстремальная кривая}%

\REitem%
   {extremal; extreme}%
   {экстремальный}%

\REitem%
   {eccentricity}%
   {эксцентриситет}%

\REitem%
   {elementary particle}%
   {элементарная частица}%

\REitem%
   {endomorphism}%
   {эндоморфизм}%

\REitem%
   {energy}%
   {энергия}%

\REitem%
   {entropy}%
   {энтропия}%

\REitem%
   {Erlanger Program}%
   {Эрлангенская программа}%

\REitem%
   {hermitian form}%
   {эрмитова форма}%

\REitem%
   {Doppler shift}%
   {эффект Доплера}%

\end{description}%

\section{Я}%
\begin{description}%

\REitem%
   {phenomenon}%
   {явление}%

\REitem%
   {phenomena}%
   {явления}%

\REitem%
   {nucleus}%
   {ядро (атома)}%

\REitem%
   {kernel}%
   {ядро (отображения)}%

\REitem%
   {Jacobian}%
   {якобиан}%

\end{description}%

\def\indexname{Name index}%
\begin{theindex}%
\indexspace%
\EName
 {Niels Henrik}%
 {Abel}%
 {Нильс Хенрик}%
 {Абель}%
\EName
 {Yakir}%
 {Aharonov}%
 {Якир}%
 {Ааронов}%
\EName
 {James Waddell}%
 {Alexander}%
 {Джеймс Уэдделл}%
 {Александер}%
\EName
 {Victor Amazaspovich}%
 {Ambartsumian}%
 {Виктор Амазаспович}%
 {Амбарцумян}%
\EName
 {Andr\'e-Marie }%
 {Amp\`ere}%
 {Андре-Мари}%
 {Ампер}%
\EName
 {Carl David}%
 {Anderson}%
 {Карл Дейвид}%
 {Андерсон}%
\EName
 {}%
 {Archimedes}%
 {}%
 {Архимед}%
\EName
 {Emil}%
 {Artin}%
 {Эмиль}%
 {Артин}%
\EName
 {Cesare}%
 {Arzel\`a}%
 {Чезаре}%
 {Арцела}%
\EName
 {Abhay}%
 {Ashtekar}%
 {Абэй}%
 {Аштекар}%
\EName
 {Alain}%
 {Aspect}%
 {Ален}%
 {Аспект}%
\EName
 {Michael}%
 {Atiyah}%
 {Майкл}%
 {Атья}%
\EName
 {Pierre Victor}%
 {Auger}%
 {Пьер Виктор}%
 {Оже}%
\EName
 {Amedeo}%
 {Avogadro}%
 {Амедео}%
 {Авогадро}%
\indexspace%
\EName
 {\'Etienne}%
 {B\'ezout}%
 {Этьенн}%
 {Безу}%
\EName
 {Walter}%
 {Baade}%
 {Вальтер}%
 {Бааде}%
\EName
 {John C.}%
 {Baez}%
 {Джон С.}%
 {Баез}%
\EName
 {Stefan}%
 {Banach}%
 {Стефан}%
 {Банах}%
\EName
 {John}%
 {Bardeen}%
 {Джон}%
 {Бардин}%
\EName
 {Asim Orhan}%
 {Barut}%
 {Асим Орхан}%
 {Барут}%
\EName
 {Tatyana}%
 {Baturina}%
 {Татьяна}%
 {Батурина}%
\EName
 {Eric Temple}%
 {Bell}%
 {Эрик Темпл}%
 {Белл}%
\EName
 {John Stewart}%
 {Bell}%
 {Джон Стюарт}%
 {Белл}%
\EName
 {Felix Alexandrovich}%
 {Berezin}%
 {Феликс Александрович}%
 {Березин}%
\EName
 {Peter Gabriel}%
 {Bergmann}%
 {П. Г.}%
 {Бергман}%
\EName
 {Jacob}%
 {Bernoulli}%
 {Яков}%
 {Бернулли}%
\EName
 {George}%
 {Birkhoff}%
 {Джордж}%
 {Биркхоф}%
\EName
 {Marietta}%
 {Blau}%
 {Мариетта}%
 {Блау}%
\EName
 {Nikolai}%
 {Bogoliubov}%
 {Николай}%
 {Боголюбов}%
\EName
 {David Joseph}%
 {Bohm}%
 {Дэвид Джозеф}%
 {Бом}%
\EName
 {Niels}%
 {Bohr}%
 {Нильс}%
 {Бор}%
\EName
 {Martin}%
 {Bojowald}%
 {Мартин}%
 {Боджовалд}%
\EName
 {Ludwig}%
 {Boltzmann}%
 {Людвиг}%
 {Больцман}%
\EName
 {J\'anos}%
 {Bolyai}%
 {Янош}%
 {Больяй}%
\EName
 {Bernhard}%
 {Bolzano}%
 {Бернард}%
 {Больцано}%
\EName
 {George}%
 {Boole}%
 {Джордж}%
 {Буль}%
\EName
 {F\'elix}%
 {Borel}%
 {Феликс}%
 {Борель}%
\EName
 {Max}%
 {Born}%
 {Макс}%
 {Борн}%
\EName
 {Satyendra Nath}%
 {Bose}%
 {Шатьендранат}%
 {Бозе}%
\EName
 {Raoul}%
 {Bott}%
 {Рауль}%
 {Ботт}%
\EName
 {Nicolas}%
 {Bourbaki}%
 {Никола}%
 {Бурбаки}%
\EName
 {Ferdinand}%
 {Brickwedde}%
 {Фердинанд}%
 {Брикведде}%
\indexspace%
\EName
 {Eugenio}%
 {Calabi}%
 {}%
 {Калаби}%
\EName
 {Georg}%
 {Cantor}%
 {Георг}%
 {Кантор}%
\EName
 {Constantin}%
 {Caratheodory}%
 {Константин}%
 {Каратеодори}%
\EName
 {Sadi}%
 {Carnot}%
 {Сади}%
 {Карно}%
\EName
 {Richard Christopher}%
 {Carrington}%
 {Ричард Кристофер}%
 {Кэррингтон}%
\EName
 {Lewis}%
 {Carroll}%
 {Льюис}%
 {Кэррол}%
\EName
 {Elie Joseph}%
 {Cartan}%
 {Эли Жозеф}%
 {Картан}%
\EName
 {Henri Paul}%
 {Cartan}%
 {Анри Поль}%
 {Картaн}%
\EName
 {Hendrik}%
 {Casimir}%
 {Хендрик}%
 {Казимир}%
\EName
 {Guido}%
 {Castelnuovo}%
 {Гвидо}%
 {Кастельнуово}%
\EName
 {Augustin Louis}%
 {Cauchy}%
 {Августин}%
 {Коши}%
\EName
 {James}%
 {Chadwick}%
 {Джеймс}%
 {Чедвик}%
\EName
 {Owen}%
 {Chamberlain}%
 {Оуэн}%
 {Чемберлен}%
\EName
 {Nikolai Grigorievich}%
 {Chebotaryov}%
 {Николай Григорьевич}%
 {Чеботарёв}%
\EName
 {Pafnuty Lvovich}%
 {Chebyshev}%
 {Пафнутий Львович}%
 {Чебышев}%
\EName
 {Shiing-Shen}%
 {Chern}%
 {Шиинг-Шен}%
 {Черн}%
\EName
 {N. A.}%
 {Chernikov}%
 {Н. А.}%
 {Черников}%
\EName
 {Claude}%
 {Chevalley}%
 {Клод}%
 {Шевалле}%
\EName
 {Geoffrey Foucar}%
 {Chew}%
 {Джеффри Фаукар}%
 {Чу}%
\EName
 {Elwin Bruno}%
 {Christoffel}%
 {Элвин Бруно}%
 {Кристоффель}%
\EName
 {Elwin Bruno}%
 {Chritoffel}%
 {Эльвин Бруно}%
 {Кристоффель}%
\EName
 {William}%
 {Clifford}%
 {Уильям}%
 {Клиффорд}%
\EName
 {Paul Moritz}%
 {Cohn}%
 {Пол Мориц}%
 {Кон}%
\EName
 {Alain}%
 {Connes}%
 {Ален}%
 {Конн}%
\EName
 {John Horton}%
 {Conway}%
 {Джон Хортон}%
 {Конвей}%
\EName
 {Leon Niels}%
 {Cooper}%
 {Леон Нил}%
 {Купер}%
\EName
 {Gustave\Hyph Gaspard}%
 {Coriolis}%
 {Гюстав Гаспар}%
 {Кориолис}%
\EName
 {Charles-Augustin de}%
 {Coulomb}%
 {Шарль Огюстен де}%
 {Кулон}%
\EName
 {Gabriel}%
 {Cramer}%
 {Габриель}%
 {Крамер}%
\EName
 {Luigi}%
 {Cremona}%
 {Луиджи}%
 {Кремона}%
\indexspace%
\EName
 {Jean Le Rond}%
 {d'Alembert}%
 {Жан Лерон}%
 {Д’Аламбер}%
\EName
 {Leonardo}%
 {da Vinci}%
 {Леонардо}%
 {да Винчи}%
\EName
 {Louis}%
 {de Broglie}%
 {Луи}%
 {де Бройль}%
\EName
 {Jean Fr\'ed\'eric Auguste}%
 {Delsarte}%
 {Жан Фредерик Огюст}%
 {Дельсарт}%
\EName
 {Ren\'e}%
 {Descartes}%
 {Рене}%
 {Декарт}%
\EName
 {James}%
 {Dewar}%
 {Джеймс}%
 {Дьюар}%
\EName
 {Bryce S.}%
 {DeWitt}%
 {Брайс С.}%
 {Девитт}%
\EName
 {Jean Alexandre Eug\`ene }%
 {Dieudonn\'e}%
 {Жан Александр Эжен}%
 {Дьедонне}%
\EName
 {Paul}%
 {Dirac}%
 {Пауль}%
 {Дирак}%
\EName
 {Johann Peter Gustav Lejeune}%
 {Dirichlet}%
 {Иоганн Петер Густав Лежён}%
 {Дирихле}%
\EName
 {Christian}%
 {Doppler}%
 {Христиан}%
 {Допплер}%
\EName
 {Paul}%
 {Dubreil}%
 {Поль}%
 {Дюбрейль}%
\EName
 {Frank Watson}%
 {Dyson}%
 {Франк Уотсон}%
 {Дайсон}%
\EName
 {Freeman John}%
 {Dyson}%
 {Фримен Джон}%
 {Дайсон}%
\indexspace%
\EName
 {Arthur Stanley}%
 {Eddington}%
 {Артур Стэнли}%
 {Эддингтон}%
\EName
 {Dmitri Fyodorovich}%
 {Egorov}%
 {Дмитрий Фёдорович}%
 {Егоров}%
\EName
 {Werner}%
 {Ehrenberg}%
 {Вернер}%
 {Эренберг}%
\EName
 {Paul}%
 {Ehrenfest}%
 {Пауль}%
 {Эренфест}%
\EName
 {Samuel}%
 {Eilenberg}%
 {Самуил}%
 {Эйленберг}%
\EName
 {Albert}%
 {Einstein}%
 {Альберт}%
 {Эйнштейн}%
\EName
 {Luther Pfahler}%
 {Eisenhart}%
 {Л. П.}%
 {Эйзенхарт}%
\EName
 {Lev}%
 {Elsgolts}%
 {Л. Э.}%
 {Эльсгольц}%
\EName
 {Federigo}%
 {Enriques}%
 {Федериго}%
 {Энрикес}%
\EName
 {Paul Sophus}%
 {Epstein}%
 {Пауль Софус}%
 {Эпштейн}%
\EName
 {}%
 {Euclid}%
 {}%
 {Эвклид}%
\EName
 {Leonhard}%
 {Euler}%
 {Леонард}%
 {Эйлер}%
\indexspace%
\EName
 {Carl}%
 {Faith}%
 {Карл}%
 {Фейс}%
\EName
 {Michael}%
 {Faraday}%
 {Майкл}%
 {Фарадей}%
\EName
 {Pierre}%
 {Fermat}%
 {Пьер}%
 {Ферма}%
\EName
 {Enrico}%
 {Fermi}%
 {Энрико}%
 {Ферми}%
\EName
 {Richard}%
 {Feynman}%
 {Ричард}%
 {Фейнман}%
\EName
 {Grigorii Mikhailovich}%
 {Fikhtengolts}%
 {Григорий Михайлович}%
 {Фихтенгольц}%
\EName
 {Serge\u{i} Pavlovich}%
 {Finikov}%
 {Сергей Павлович}%
 {Фиников}%
\EName
 {Paul}%
 {Finsler}%
 {Пауль}%
 {Финслер}%
\EName
 {Armand\Hyph Hippolyte\Hyph Louis}%
 {Fizeau}%
 {Арман Ипполит Луи}%
 {Физо}%
\EName
 {John Ambrose}%
 {Fleming}%
 {Джон Амброз}%
 {Флеминг}%
\EName
 {Serge\u{i} Vasil'evich}%
 {Fomin}%
 {Сергей Васильевич}%
 {Фомин}%
\EName
 {James}%
 {Forbes}%
 {Джеймс}%
 {Форбс}%
\EName
 {Joseph}%
 {Fourier}%
 {Джозеф}%
 {Фурье}%
\EName
 {Maurice Ren\'e}%
 {Fr\'echet}%
 {Морис Рене}%
 {Фреше}%
\EName
 {George}%
 {Francis}%
 {Джордж}%
 {Франсис}%
\EName
 {Benjamin}%
 {Franklin}%
 {Франклин}%
 {Бенджамин}%
\EName
 {Erik Ivar}%
 {Fredholm}%
 {Эрик Ивар}%
 {Фредгольм}%
\EName
 {Alexander}%
 {Friedmann}%
 {Александр}%
 {Фридман}%
\EName
 {Ferdinand Georg}%
 {Frobenius}%
 {Фердинанд Георг}%
 {Фробениус}%
\EName
 {Guido}%
 {Fubini}%
 {Гвидо}%
 {Фубини}%
\indexspace%
\EName
 {Kurt Friedrich}%
 {G\"odel}%
 {Курт Фридрих}%
 {Гёдель}%
\EName
 {Ren\'e Eug\`ene}%
 {G\^ateaux}%
 {Рене Эжен}%
 {Гато}%
\EName
 {Galileo}%
 {Galilei}%
 {Галилео}%
 {Галилей}%
\EName
 {}%
 {Galois}%
 {}%
 {Галуа}%
\EName
 {Carl Friedrich}%
 {Gauss}%
 {Карл Фридрих}%
 {Гаусс}%
\EName
 {Hans}%
 {Geiger}%
 {Ханс}%
 {Гейгер}%
\EName
 {Israel Moiseevich}%
 {Gelfand}%
 {Израиль Моисеевич}%
 {Гельфанд}%
\EName
 {Murray}%
 {Gell-Mann}%
 {Мюррей}%
 {Гелл-Манн}%
\EName
 {Josiah Willard}%
 {Gibbs}%
 {Джозайя Уиллард}%
 {Гиббс}%
\EName
 {Donald Arthur}%
 {Glaser}%
 {Дональд Артур}%
 {Глазер}%
\EName
 {Marvin Leonard}%
 {Goldberger}%
 {Мэрвин Леонард}%
 {Гольдбергер}%
\EName
 {Maurice}%
 {Goldhaber}%
 {Морис}%
 {Гольдхабер}%
\EName
 {Walter}%
 {Gordon}%
 {Вальтер}%
 {Гордон}%
\EName
 {Samuel Abraham}%
 {Goudsmit}%
 {Сэмюэл Абрахам}%
 {Гаудсмит}%
\EName
 {Hermann G\"unther}%
 {Grassmann }%
 {Герман Гюнтер}%
 {Грассман}%
\EName
 {Dmitry Aleksandrovich}%
 {Grave}%
 {Дмитрий Александрович}%
 {Граве}%
\EName
 {Marcel}%
 {Grossman}%
 {Марсель}%
 {Гроссман}%
\EName
 {Alexander}%
 {Grothendieck}%
 {Александр}%
 {Гротендик}%
\indexspace%
\EName
 {Jacques}%
 {Hadamard}%
 {Жак}%
 {Адамар}%
\EName
 {Georg Karl Wilhelm}%
 {Hamel}%
 {Георг Карл Вильгельм}%
 {Гамель}%
\EName
 {William Rowan}%
 {Hamilton}%
 {Вильям Роуэн}%
 {Гамильтон}%
\EName
 {Serge}%
 {Haroche}%
 {Серж}%
 {Арош}%
\EName
 {Felix}%
 {Hausdorff}%
 {Феликс}%
 {Хаусдорф}%
\EName
 {S. W.}%
 {Hawking}%
 {С. В.}%
 {Хокинг}%
\EName
 {Friedrich W.}%
 {Hehl}%
 {Фридрих}%
 {Хель}%
\EName
 {Werner}%
 {Heisenberg}%
 {Вернер}%
 {Гейзенберг}%
\EName
 {Kurt}%
 {Hensel}%
 {Курт}%
 {Хензель}%
\EName
 {William}%
 {Herschel}%
 {Вильям}%
 {Гершель}%
\EName
 {Heinrich Rudolf}%
 {Hertz}%
 {Генрих Рудольф}%
 {Герц}%
\EName
 {Ejnar}%
 {Hertzsprung}%
 {Эйнар}%
 {Герцшпрунг}%
\EName
 {Victor Franz}%
 {Hess}%
 {Виктор Франц}%
 {Гесс}%
\EName
 {Antony}%
 {Hewish}%
 {Энтони}%
 {Хьюиш}%
\EName
 {David}%
 {Hilbert}%
 {Давид}%
 {Гильберт}%
\EName
 {Friedrich Ernst Peter}%
 {Hirzebruch}%
 {Фридрих Эрнст Петер}%
 {Хирцебрух}%
\EName
 {Heinz}%
 {Hopf}%
 {Хайнц}%
 {Хопф}%
\EName
 {Fred}%
 {Hoyle}%
 {Фред}%
 {Хойл}%
\EName
 {Edwin}%
 {Hubble}%
 {Эдвин}%
 {Хаббл}%
\EName
 {Adolf}%
 {Hurwitz}%
 {Адольф}%
 {Гурвиц}%
\indexspace%
\EName
 {Leopold}%
 {Infeld}%
 {Леопольд}%
 {Инфельд}%
\EName
 {Dmitri}%
 {Iwanenko}%
 {Дмитрий Дмитриевич}%
 {Иваненко}%
\indexspace%
\EName
 {Carl}%
 {Jacobi}%
 {Карл}%
 {Якоби}%
\EName
 {Nathan}%
 {Jacobson}%
 {Натан}%
 {Джекобсон}%
\EName
 {Pascual}%
 {Jordan}%
 {Паскуаль}%
 {Иордан}%
\indexspace%
\EName
 {Heike}%
 {Kamerlingh Onnes}%
 {Хейке}%
 {Камерлинг-Оннес}%
\EName
 {Pyotr}%
 {Kapitsa}%
 {Пётр}%
 {Капица}%
\EName
 {Max}%
 {Karoubi}%
 {Макс}%
 {Каруби}%
\EName
 {Johann}%
 {Kepler}%
 {Иоганн}%
 {Кеплер}%
\EName
 {Roy Patrick}%
 {Kerr}%
 {Рой Патрик}%
 {Керр}%
\EName
 {Wilhelm}%
 {Killing}%
 {Вильгельм}%
 {Киллинг}%
\EName
 {Gustav Robert}%
 {Kirchhoff}%
 {Густав Роберт}%
 {Кирхгоф}%
\EName
 {Felix}%
 {Klein}%
 {Феликс}%
 {Клейн}%
\EName
 {Oskar}%
 {Klein}%
 {Оскар}%
 {Клейн}%
\EName
 {Shoshichi}%
 {Kobayashi}%
 {С.}%
 {Кобаяси}%
\EName
 {Andre\u{i} Nikolaevich}%
 {Kolmogorov}%
 {Андрей Николаевич}%
 {Колмогоров}%
\EName
 {Sofia Vasilyevna}%
 {Kovalevskaya}%
 {Софья Васильевна}%
 {Ковалевская}%
\EName
 {Robert}%
 {Kraichnan}%
 {Роберт}%
 {Крайчнан}%
\EName
 {Leopold}%
 {Kronecker}%
 {Леопольд}%
 {Кронекер}%
\EName
 {Martin}%
 {Kruskal}%
 {Мартин}%
 {Крускал}%
\EName
 {Kazimierz}%
 {Kuratowski}%
 {Казимир}%
 {Куратовский}%
\EName
 {Martin Wilhelm}%
 {Kutta}%
 {Мартин Вильгельм}%
 {Кутта}%
\indexspace%
\EName
 {Guillaume Fran\c cois de}%
 {L'H\^opital}%
 {Гиймон Франсуа}%
 {де Лопиталь}%
\EName
 {Joseph Louis}%
 {Lagrange}%
 {Жозеф Луи}%
 {Лагранж}%
\EName
 {Tsit Yuen}%
 {Lam}%
 {Цит Юань}%
 {Лам}%
\EName
 {Lev}%
 {Landau}%
 {Лев}%
 {Ландау}%
\EName
 {Serge}%
 {Lang}%
 {Серж}%
 {Ленг}%
\EName
 {Pierre-Simon}%
 {Laplace}%
 {Пьер-Симон}%
 {Лаплас}%
\EName
 {Henrietta Swan}%
 {Leavitt}%
 {Генриетта Суон}%
 {Ливитт}%
\EName
 {Henri}%
 {Lebesgue}%
 {Анри}%
 {Лебег}%
\EName
 {Joel}%
 {Lebowitz}%
 {Джоэль}%
 {Лейбовиц}%
\EName
 {Leon Max}%
 {Lederman}%
 {Леон Макс}%
 {Ледерман}%
\EName
 {Adrien-Marie}%
 {Legendre}%
 {Адриен Мари}%
 {Лежандр}%
\EName
 {Gottfried Wilhelm}%
 {Leibniz}%
 {Готфрид Вильгельм}%
 {Лейбниц}%
\EName
 {Joseph}%
 {Lense}%
 {Джозеф}%
 {Лензе}%
\EName
 {Beppo}%
 {Levi}%
 {Беппо}%
 {Леви}%
\EName
 {Andre}%
 {Lichnerowicz}%
 {Андрэ}%
 {Лихнерович}%
\EName
 {Andrei}%
 {Linde}%
 {Андрей}%
 {Линде}%
\EName
 {Joseph}%
 {Liouville}%
 {Жозеф}%
 {Лиувилль}%
\EName
 {Nikolai Ivanovich}%
 {Lobachevsky}%
 {Николай Иванович}%
 {Лобачевский}%
\EName
 {Hendrik Antoon}%
 {Lorentz}%
 {Хендрик Антон}%
 {Лоренц}%
\EName
 {Edward Norton}%
 {Lorenz}%
 {Эдвард Нортон}%
 {Лоренц}%
\indexspace%
\EName
 {Ettore}%
 {Majorana}%
 {Этторе}%
 {Майорана}%
\EName
 {James Clerk}%
 {Maxwell}%
 {Джеймс Клерк}%
 {Максвелл}%
\EName
 {Ulf}%
 {Meissner}%
 {Ульф}%
 {Мейснер}%
\EName
 {Lise}%
 {Meitner}%
 {Лиза}%
 {Мейтнер}%
\EName
 {Eugene}%
 {Merzbacher}%
 {Юджин}%
 {Мерцбахер}%
\EName
 {Robert}%
 {Millikan}%
 {Роберт}%
 {Миликен}%
\EName
 {John Willard}%
 {Milnor}%
 {Джон Уиллард}%
 {Милнор}%
\EName
 {Hermann}%
 {Minkowski}%
 {Герман}%
 {Минковский}%
\EName
 {Charles W.}%
 {Misner}%
 {Чарлз}%
 {Мизнер}%
\EName
 {Eliakim Hastings}%
 {Moore}%
 {Элиаким Гастингс}%
 {Мур}%
\EName
 {Holger}%
 {Mueller}%
 {Холгер}%
 {Мюллер}%
\EName
 {David Bryant}%
 {Mumford}%
 {Дэвид Брайант}%
 {Мамфорд}%
\indexspace%
\EName
 {}%
 {Newton}%
 {}%
 {Ньютон}%
\EName
 {John von}%
 {Neyman}%
 {Джон фон}%
 {Нейман}%
\EName
 {Emmy}%
 {Noether}%
 {Эмми}%
 {Нётер}%
\EName
 {Katsumi}%
 {Nomizu}%
 {К.}%
 {Номидзу}%
\indexspace%
\EName
 {Yuri}%
 {Obukhov}%
 {Юрий}%
 {Обухов}%
\EName
 {William of}%
 {Ockham}%
 {Уильям}%
 {Оккам}%
\EName
 {Peter}%
 {Olver}%
 {Питер}%
 {Олвер}%
\EName
 {Robert}%
 {Oppenheimer}%
 {Роберт}%
 {Опенгеймер}%
\indexspace%
\EName
 {Moritz}%
 {Pasch}%
 {Мориц}%
 {Паш}%
\EName
 {Blaise}%
 {Paskal}%
 {Блез}%
 {Паскаль}%
\EName
 {Wolfgang}%
 {Pauli}%
 {Вольфганг}%
 {Паули}%
\EName
 {Giuseppe}%
 {Peano}%
 {Джузеппе}%
 {Пеано}%
\EName
 {James}%
 {Peebles}%
 {Джеймс}%
 {Пиблс}%
\EName
 {Rudolf Ernst}%
 {Peierls}%
 {Рудольф Эрнст}%
 {Пайерлс}%
\EName
 {Christian Heinrich Friedrich}%
 {Peters}%
 {Кристиан Генрих Фридрих}%
 {Петерс}%
\EName
 {Julius}%
 {Pl\"ucker}%
 {Юлий}%
 {Плюкер}%
\EName
 {Max}%
 {Planck}%
 {Макс}%
 {Планк}%
\EName
 {Boris}%
 {Podolsky}%
 {Борис Яковлевич}%
 {Подольский}%
\EName
 {Jules Henri}%
 {Poincar\'e}%
 {Жюль Анри}%
 {Пуанкаре}%
\EName
 {Sim\'eon-Denis}%
 {Poisson}%
 {Симеон-Дени}%
 {Пуассон}%
\EName
 {George}%
 {Polya}%
 {Джордж}%
 {Пойа}%
\EName
 {Isaac}%
 {Pomeranchuk}%
 {Исаак Яковлевич}%
 {Померанчук}%
\EName
 {Jean-Fran\c{c}ois}%
 {Pommaret}%
 {Жан-Франсуа}%
 {Поммаре}%
\EName
 {Lev Semenovich}%
 {Pontryagin}%
 {Лев Семёнович}%
 {Понтрягин}%
\EName
 {Cecil Frank}%
 {Powell}%
 {Сесил Фрэнк}%
 {Пауэлл}%
\indexspace%
\EName
 {Ryszard}%
 {R\d aczka}%
 {Р}%
 {Рончка}%
\EName
 {Isidor Isaac}%
 {Rabi}%
 {Исидор Айзек}%
 {Раби}%
\EName
 {Martin}%
 {Rees}%
 {Мартин}%
 {Рис}%
\EName
 {Kurt Werner Friedrich}%
 {Reidemeister}%
 {Курт Вернер Фридрих}%
 {Рейдемейстер}%
\EName
 {William}%
 {Reynolds}%
 {Вильям}%
 {Рейнольдс}%
\EName
 {Gregorio}%
 {Ricci}%
 {Грегорио}%
 {Риччи}%
\EName
 {Bernhard}%
 {Riemann}%
 {Бернхард}%
 {Риман}%
\EName
 {George Neil}%
 {Robertson}%
 {Джордж Нил}%
 {Робертсон}%
\EName
 {Howard Percy}%
 {Robertson}%
 {Говард Перси}%
 {Робертсон}%
\EName
 {Ole Christensen}%
 {Roemer}%
 {Оле Кристенсен}%
 {Рёмер}%
\EName
 {Michel}%
 {Rolle}%
 {Мишель}%
 {Ролль}%
\EName
 {Nathan}%
 {Rosen}%
 {Натан}%
 {Розен}%
\EName
 {Marshall Nicholas}%
 {Rosenbluth}%
 {Маршалл Николас}%
 {Розенблют}%
\EName
 {Carlo}%
 {Rovelli}%
 {Карло}%
 {Ровелли}%
\EName
 {Warren de la}%
 {Rue}%
 {Уоррен де ла}%
 {Рю}%
\EName
 {Hanno}%
 {Rund}%
 {Ханно}%
 {Рунд}%
\EName
 {Karl David Tolme}%
 {Runge}%
 {Карл Давид Тольме}%
 {Рунге}%
\EName
 {Bertrand}%
 {Russell}%
 {Бертран}%
 {Рассел}%
\EName
 {Henry Norris}%
 {Russell}%
 {Генри Норрис}%
 {Расселл}%
\EName
 {Ernest}%
 {Rutherford}%
 {Эрнест}%
 {Резерфорд}%
\indexspace%
\EName
 {Edwin Ernest}%
 {Salpeter}%
 {Эдвин Эрнест}%
 {Сальпетер}%
\EName
 {}%
 {Schmidt}%
 {}%
 {Шмидт}%
\EName
 {Erwin}%
 {Schr$\ddot{o}$dinger}%
 {Эрвин}%
 {Шредингер}%
\EName
 {John}%
 {Schrieffer}%
 {Джон}%
 {Шриффер}%
\EName
 {Jacob}%
 {Schwartz}%
 {Дж.}%
 {Шварц}%
\EName
 {Laurent Moise}%
 {Schwartz}%
 {Лоран Моиз}%
 {Шварц}%
\EName
 {Karl}%
 {Schwarzschild}%
 {Карл}%
 {Шварцшильд}%
\EName
 {Julian}%
 {Schwinger}%
 {Джулиус}%
 {Швингер}%
\EName
 {Corrado}%
 {Segre}%
 {Коррадо}%
 {Сегре}%
\EName
 {Jean-Pierre}%
 {Serre}%
 {Жан Пьер}%
 {Серр}%
\EName
 {Waclaw}%
 {Sierpinski}%
 {Вацлав}%
 {Серпинский}%
\EName
 {James Harris}%
 {Simons}%
 {Джеймс Харрис}%
 {Симонс}%
\EName
 {Dmitriy}%
 {Skobeltsin}%
 {Дмитрий}%
 {Скобельцин}%
\EName
 {Tony}%
 {Skyrme}%
 {Тони}%
 {Скайрм}%
\EName
 {Lee}%
 {Smolin}%
 {Ли}%
 {Смолин}%
\EName
 {Hartland Sweet}%
 {Snyder}%
 {Хартланд Свит}%
 {Снайдер}%
\EName
 {Arnold}%
 {Sommerfeld}%
 {Арнольд}%
 {Зоммерфельд}%
\EName
 {Norman}%
 {Steenrod}%
 {Норман}%
 {Стинрод}%
\EName
 {Paul}%
 {Steinhardt}%
 {Поль}%
 {Штейнхардт}%
\EName
 {Shlomo}%
 {Sternberg}%
 {Соломон}%
 {Стернберг}%
\EName
 {Eduard}%
 {Stiefel}%
 {Эдуард}%
 {Штифель}%
\EName
 {Thomas Joannes}%
 {Stieltjes}%
 {Томас Иоаннес}%
 {Стилтьес}%
\EName
 {James}%
 {Stirling}%
 {Джеймс}%
 {Стирлинг}%
\EName
 {George Gabriel}%
 {Stokes}%
 {Джордж Габриель}%
 {Стокс}%
\EName
 {Eduard}%
 {Study}%
 {Эдуард}%
 {Стади}%
\EName
 {Jacques Charles Fracois}%
 {Sturm}%
 {Жак Шарль Франсуа}%
 {Штурм}%
\EName
 {Anthony}%
 {Sudbery}%
 {Энтони}%
 {Садбери}%
\EName
 {Leonard}%
 {Susskind}%
 {Леонард}%
 {Сасскинд}%
\EName
 {Leo}%
 {Szilard}%
 {Лео}%
 {Сцилард}%
\indexspace%
\EName
 {Alfred}%
 {Tarski}%
 {Альфред}%
 {Тарский}%
\EName
 {Hans}%
 {Thirring}%
 {Ганс}%
 {Тирринг}%
\EName
 {Walter }%
 {Thirring}%
 {Вальтер}%
 {Тирринг}%
\EName
 {Ren\'e Fr\'ed\'eric}%
 {Thom}%
 {Рене Фредерик}%
 {Том}%
\EName
 {George Paget}%
 {Thomson}%
 {Джордж Паджет}%
 {Томсон}%
\EName
 {Kip S.}%
 {Thorne}%
 {Кип}%
 {Торн}%
\EName
 {Richard}%
 {Tolman}%
 {Ричард}%
 {Толмен}%
\EName
 {Evangelista}%
 {Torricelli}%
 {Евангелиста}%
 {Торричелли}%
\indexspace%
\EName
 {George Eugene}%
 {Uhlenbeck}%
 {Джордж Юджин}%
 {Уленбек}%
\EName
 {Harold Clayton}%
 {Urey}%
 {Гарольд Клейтон}%
 {Ури}%
\indexspace%
\EName
 {James}%
 {Van Allen}%
 {Джеймс}%
 {Ван Аллен}%
\EName
 {John}%
 {Venn}%
 {Джон}%
 {Венн }%
\EName
 {Fran\c cois}%
 {Vi\'ete}%
 {Франсуа}%
 {Виет}%
\EName
 {Leopold}%
 {Vietoris}%
 {Леопольд}%
 {Фиторис}%
\EName
 {Alex}%
 {Vilenkin}%
 {Алекс}%
 {Виленкин}%
\EName
 {Valerii}%
 {Vinokur}%
 {Валерий}%
 {Винокур}%
\EName
 {Max}%
 {von Laue}%
 {Макс}%
 {фон Лауэ}%
\indexspace%
\EName
 {Wilhelm Eduard}%
 {Weber}%
 {Вильгельм Эдуард}%
 {Вебер}%
\EName
 {Karl Theodor Wilhelm}%
 {Weierstrass}%
 {Карл Теодор Вильгельм}%
 {Вейерштрасс}%
\EName
 {Andr\'e}%
 {Weil}%
 {Андре}%
 {Вейль}%
\EName
 {Steven}%
 {Weinberg}%
 {Стивен}%
 {Вайнберг}%
\EName
 {Hermann}%
 {Weyl}%
 {Герман}%
 {Вейль}%
\EName
 {John Archibald}%
 {Wheeler}%
 {Джон Арчибалд}%
 {Уилер}%
\EName
 {Alfred North}%
 {Whitehead}%
 {Альфред Норт}%
 {Уайтхед}%
\EName
 {Hassler}%
 {Whitney}%
 {Хасслер}%
 {Уитни}%
\EName
 {Eugene}%
 {Wigner}%
 {Юджин}%
 {Вигнер}%
\EName
 {Frank}%
 {Wilczek}%
 {Франк}%
 {Вилчек}%
\EName
 {Robert}%
 {Wilson}%
 {Роберт}%
 {Вильсон}%
\indexspace%
\EName
 {Chen Ning}%
 {Yang}%
 {Чжэньнин}%
 {Янг}%
\EName
 {Shing-Tung}%
 {Yau}%
 {Шин-Тан}%
 {Яу}%
\EName
 {Hideki}%
 {Yukawa}%
 {Хидеки}%
 {Юкава}%
\indexspace%
\EName
 {Oscar}%
 {Zariski}%
 {Оскар}%
 {Зарисский}%
\EName
 {Yakov Borisovich}%
 {Zel’dovich}%
 {Яков Борисович}%
 {Зельдович}%
\EName
 {Ernst Friedrich Ferdinand}%
 {Zermelo}%
 {Эрнст Фридрих Фердинанд}%
 {Цермело}%
\end{theindex}%

\def\indexname{Именной указатель}%
\begin{theindex}%
\indexspace%
\RName
 {Yakir}%
 {Aharonov}%
 {Якир}%
 {Ааронов}%
\RName
 {Niels Henrik}%
 {Abel}%
 {Нильс Хенрик}%
 {Абель}%
\RName
 {Amedeo}%
 {Avogadro}%
 {Амедео}%
 {Авогадро}%
\RName
 {Jacques}%
 {Hadamard}%
 {Жак}%
 {Адамар}%
\RName
 {James Waddell}%
 {Alexander}%
 {Джеймс Уэдделл}%
 {Александер}%
\RName
 {Victor Amazaspovich}%
 {Ambartsumian}%
 {Виктор Амазаспович}%
 {Амбарцумян}%
\RName
 {Andr\'e-Marie }%
 {Amp\`ere}%
 {Андре-Мари}%
 {Ампер}%
\RName
 {Carl David}%
 {Anderson}%
 {Карл Дейвид}%
 {Андерсон}%
\RName
 {Serge}%
 {Haroche}%
 {Серж}%
 {Арош}%
\RName
 {Emil}%
 {Artin}%
 {Эмиль}%
 {Артин}%
\RName
 {}%
 {Archimedes}%
 {}%
 {Архимед}%
\RName
 {Cesare}%
 {Arzel\`a}%
 {Чезаре}%
 {Арцела}%
\RName
 {Alain}%
 {Aspect}%
 {Ален}%
 {Аспект}%
\RName
 {Michael}%
 {Atiyah}%
 {Майкл}%
 {Атья}%
\RName
 {Abhay}%
 {Ashtekar}%
 {Абэй}%
 {Аштекар}%
\indexspace%
\RName
 {Walter}%
 {Baade}%
 {Вальтер}%
 {Бааде}%
\RName
 {John C.}%
 {Baez}%
 {Джон С.}%
 {Баез}%
\RName
 {Stefan}%
 {Banach}%
 {Стефан}%
 {Банах}%
\RName
 {John}%
 {Bardeen}%
 {Джон}%
 {Бардин}%
\RName
 {Asim Orhan}%
 {Barut}%
 {Асим Орхан}%
 {Барут}%
\RName
 {Tatyana}%
 {Baturina}%
 {Татьяна}%
 {Батурина}%
\RName
 {\'Etienne}%
 {B\'ezout}%
 {Этьенн}%
 {Безу}%
\RName
 {John Stewart}%
 {Bell}%
 {Джон Стюарт}%
 {Белл}%
\RName
 {Eric Temple}%
 {Bell}%
 {Эрик Темпл}%
 {Белл}%
\RName
 {Benjamin}%
 {Franklin}%
 {Франклин}%
 {Бенджамин}%
\RName
 {Peter Gabriel}%
 {Bergmann}%
 {П. Г.}%
 {Бергман}%
\RName
 {Felix Alexandrovich}%
 {Berezin}%
 {Феликс Александрович}%
 {Березин}%
\RName
 {Jacob}%
 {Bernoulli}%
 {Яков}%
 {Бернулли}%
\RName
 {George}%
 {Birkhoff}%
 {Джордж}%
 {Биркхоф}%
\RName
 {Marietta}%
 {Blau}%
 {Мариетта}%
 {Блау}%
\RName
 {Nikolai}%
 {Bogoliubov}%
 {Николай}%
 {Боголюбов}%
\RName
 {Martin}%
 {Bojowald}%
 {Мартин}%
 {Боджовалд}%
\RName
 {Satyendra Nath}%
 {Bose}%
 {Шатьендранат}%
 {Бозе}%
\RName
 {Bernhard}%
 {Bolzano}%
 {Бернард}%
 {Больцано}%
\RName
 {Ludwig}%
 {Boltzmann}%
 {Людвиг}%
 {Больцман}%
\RName
 {J\'anos}%
 {Bolyai}%
 {Янош}%
 {Больяй}%
\RName
 {David Joseph}%
 {Bohm}%
 {Дэвид Джозеф}%
 {Бом}%
\RName
 {Niels}%
 {Bohr}%
 {Нильс}%
 {Бор}%
\RName
 {F\'elix}%
 {Borel}%
 {Феликс}%
 {Борель}%
\RName
 {Max}%
 {Born}%
 {Макс}%
 {Борн}%
\RName
 {Raoul}%
 {Bott}%
 {Рауль}%
 {Ботт}%
\RName
 {Ferdinand}%
 {Brickwedde}%
 {Фердинанд}%
 {Брикведде}%
\RName
 {George}%
 {Boole}%
 {Джордж}%
 {Буль}%
\RName
 {Nicolas}%
 {Bourbaki}%
 {Никола}%
 {Бурбаки}%
\indexspace%
\RName
 {Steven}%
 {Weinberg}%
 {Стивен}%
 {Вайнберг}%
\RName
 {James}%
 {Van Allen}%
 {Джеймс}%
 {Ван Аллен}%
\RName
 {Wilhelm Eduard}%
 {Weber}%
 {Вильгельм Эдуард}%
 {Вебер}%
\RName
 {Karl Theodor Wilhelm}%
 {Weierstrass}%
 {Карл Теодор Вильгельм}%
 {Вейерштрасс}%
\RName
 {Andr\'e}%
 {Weil}%
 {Андре}%
 {Вейль}%
\RName
 {Hermann}%
 {Weyl}%
 {Герман}%
 {Вейль}%
\RName
 {John}%
 {Venn}%
 {Джон}%
 {Венн }%
\RName
 {Eugene}%
 {Wigner}%
 {Юджин}%
 {Вигнер}%
\RName
 {Fran\c cois}%
 {Vi\'ete}%
 {Франсуа}%
 {Виет}%
\RName
 {Alex}%
 {Vilenkin}%
 {Алекс}%
 {Виленкин}%
\RName
 {Frank}%
 {Wilczek}%
 {Франк}%
 {Вилчек}%
\RName
 {Robert}%
 {Wilson}%
 {Роберт}%
 {Вильсон}%
\RName
 {Valerii}%
 {Vinokur}%
 {Валерий}%
 {Винокур}%
\indexspace%
\RName
 {Galileo}%
 {Galilei}%
 {Галилео}%
 {Галилей}%
\RName
 {}%
 {Galois}%
 {}%
 {Галуа}%
\RName
 {Georg Karl Wilhelm}%
 {Hamel}%
 {Георг Карл Вильгельм}%
 {Гамель}%
\RName
 {William Rowan}%
 {Hamilton}%
 {Вильям Роуэн}%
 {Гамильтон}%
\RName
 {Ren\'e Eug\`ene}%
 {G\^ateaux}%
 {Рене Эжен}%
 {Гато}%
\RName
 {Samuel Abraham}%
 {Goudsmit}%
 {Сэмюэл Абрахам}%
 {Гаудсмит}%
\RName
 {Carl Friedrich}%
 {Gauss}%
 {Карл Фридрих}%
 {Гаусс}%
\RName
 {Kurt Friedrich}%
 {G\"odel}%
 {Курт Фридрих}%
 {Гёдель}%
\RName
 {Hans}%
 {Geiger}%
 {Ханс}%
 {Гейгер}%
\RName
 {Werner}%
 {Heisenberg}%
 {Вернер}%
 {Гейзенберг}%
\RName
 {Murray}%
 {Gell-Mann}%
 {Мюррей}%
 {Гелл-Манн}%
\RName
 {Israel Moiseevich}%
 {Gelfand}%
 {Израиль Моисеевич}%
 {Гельфанд}%
\RName
 {Heinrich Rudolf}%
 {Hertz}%
 {Генрих Рудольф}%
 {Герц}%
\RName
 {Ejnar}%
 {Hertzsprung}%
 {Эйнар}%
 {Герцшпрунг}%
\RName
 {William}%
 {Herschel}%
 {Вильям}%
 {Гершель}%
\RName
 {Victor Franz}%
 {Hess}%
 {Виктор Франц}%
 {Гесс}%
\RName
 {Josiah Willard}%
 {Gibbs}%
 {Джозайя Уиллард}%
 {Гиббс}%
\RName
 {David}%
 {Hilbert}%
 {Давид}%
 {Гильберт}%
\RName
 {Donald Arthur}%
 {Glaser}%
 {Дональд Артур}%
 {Глазер}%
\RName
 {Marvin Leonard}%
 {Goldberger}%
 {Мэрвин Леонард}%
 {Гольдбергер}%
\RName
 {Maurice}%
 {Goldhaber}%
 {Морис}%
 {Гольдхабер}%
\RName
 {Walter}%
 {Gordon}%
 {Вальтер}%
 {Гордон}%
\RName
 {Dmitry Aleksandrovich}%
 {Grave}%
 {Дмитрий Александрович}%
 {Граве}%
\RName
 {Hermann G\"unther}%
 {Grassmann }%
 {Герман Гюнтер}%
 {Грассман}%
\RName
 {Marcel}%
 {Grossman}%
 {Марсель}%
 {Гроссман}%
\RName
 {Alexander}%
 {Grothendieck}%
 {Александр}%
 {Гротендик}%
\RName
 {Adolf}%
 {Hurwitz}%
 {Адольф}%
 {Гурвиц}%
\indexspace%
\RName
 {Jean Le Rond}%
 {d'Alembert}%
 {Жан Лерон}%
 {Д’Аламбер}%
\RName
 {Leonardo}%
 {da Vinci}%
 {Леонардо}%
 {да Винчи}%
\RName
 {Frank Watson}%
 {Dyson}%
 {Франк Уотсон}%
 {Дайсон}%
\RName
 {Freeman John}%
 {Dyson}%
 {Фримен Джон}%
 {Дайсон}%
\RName
 {Louis}%
 {de Broglie}%
 {Луи}%
 {де Бройль}%
\RName
 {Bryce S.}%
 {DeWitt}%
 {Брайс С.}%
 {Девитт}%
\RName
 {Ren\'e}%
 {Descartes}%
 {Рене}%
 {Декарт}%
\RName
 {Jean Fr\'ed\'eric Auguste}%
 {Delsarte}%
 {Жан Фредерик Огюст}%
 {Дельсарт}%
\RName
 {Nathan}%
 {Jacobson}%
 {Натан}%
 {Джекобсон}%
\RName
 {Paul}%
 {Dirac}%
 {Пауль}%
 {Дирак}%
\RName
 {Johann Peter Gustav Lejeune}%
 {Dirichlet}%
 {Иоганн Петер Густав Лежён}%
 {Дирихле}%
\RName
 {Christian}%
 {Doppler}%
 {Христиан}%
 {Допплер}%
\RName
 {Jean Alexandre Eug\`ene }%
 {Dieudonn\'e}%
 {Жан Александр Эжен}%
 {Дьедонне}%
\RName
 {James}%
 {Dewar}%
 {Джеймс}%
 {Дьюар}%
\RName
 {Paul}%
 {Dubreil}%
 {Поль}%
 {Дюбрейль}%
\indexspace%
\RName
 {Dmitri Fyodorovich}%
 {Egorov}%
 {Дмитрий Фёдорович}%
 {Егоров}%
\indexspace%
\RName
 {Oscar}%
 {Zariski}%
 {Оскар}%
 {Зарисский}%
\RName
 {Yakov Borisovich}%
 {Zel’dovich}%
 {Яков Борисович}%
 {Зельдович}%
\RName
 {Arnold}%
 {Sommerfeld}%
 {Арнольд}%
 {Зоммерфельд}%
\indexspace%
\RName
 {Dmitri}%
 {Iwanenko}%
 {Дмитрий Дмитриевич}%
 {Иваненко}%
\RName
 {Leopold}%
 {Infeld}%
 {Леопольд}%
 {Инфельд}%
\RName
 {Pascual}%
 {Jordan}%
 {Паскуаль}%
 {Иордан}%
\indexspace%
\RName
 {Hendrik}%
 {Casimir}%
 {Хендрик}%
 {Казимир}%
\RName
 {Eugenio}%
 {Calabi}%
 {}%
 {Калаби}%
\RName
 {Heike}%
 {Kamerlingh Onnes}%
 {Хейке}%
 {Камерлинг-Оннес}%
\RName
 {Georg}%
 {Cantor}%
 {Георг}%
 {Кантор}%
\RName
 {Pyotr}%
 {Kapitsa}%
 {Пётр}%
 {Капица}%
\RName
 {Constantin}%
 {Caratheodory}%
 {Константин}%
 {Каратеодори}%
\RName
 {Sadi}%
 {Carnot}%
 {Сади}%
 {Карно}%
\RName
 {Henri Paul}%
 {Cartan}%
 {Анри Поль}%
 {Картaн}%
\RName
 {Elie Joseph}%
 {Cartan}%
 {Эли Жозеф}%
 {Картан}%
\RName
 {Max}%
 {Karoubi}%
 {Макс}%
 {Каруби}%
\RName
 {Guido}%
 {Castelnuovo}%
 {Гвидо}%
 {Кастельнуово}%
\RName
 {Johann}%
 {Kepler}%
 {Иоганн}%
 {Кеплер}%
\RName
 {Roy Patrick}%
 {Kerr}%
 {Рой Патрик}%
 {Керр}%
\RName
 {Wilhelm}%
 {Killing}%
 {Вильгельм}%
 {Киллинг}%
\RName
 {Gustav Robert}%
 {Kirchhoff}%
 {Густав Роберт}%
 {Кирхгоф}%
\RName
 {Oskar}%
 {Klein}%
 {Оскар}%
 {Клейн}%
\RName
 {Felix}%
 {Klein}%
 {Феликс}%
 {Клейн}%
\RName
 {William}%
 {Clifford}%
 {Уильям}%
 {Клиффорд}%
\RName
 {Shoshichi}%
 {Kobayashi}%
 {С.}%
 {Кобаяси}%
\RName
 {Sofia Vasilyevna}%
 {Kovalevskaya}%
 {Софья Васильевна}%
 {Ковалевская}%
\RName
 {Andre\u{i} Nikolaevich}%
 {Kolmogorov}%
 {Андрей Николаевич}%
 {Колмогоров}%
\RName
 {Paul Moritz}%
 {Cohn}%
 {Пол Мориц}%
 {Кон}%
\RName
 {John Horton}%
 {Conway}%
 {Джон Хортон}%
 {Конвей}%
\RName
 {Alain}%
 {Connes}%
 {Ален}%
 {Конн}%
\RName
 {Gustave\Hyph Gaspard}%
 {Coriolis}%
 {Гюстав Гаспар}%
 {Кориолис}%
\RName
 {Augustin Louis}%
 {Cauchy}%
 {Августин}%
 {Коши}%
\RName
 {Robert}%
 {Kraichnan}%
 {Роберт}%
 {Крайчнан}%
\RName
 {Gabriel}%
 {Cramer}%
 {Габриель}%
 {Крамер}%
\RName
 {Luigi}%
 {Cremona}%
 {Луиджи}%
 {Кремона}%
\RName
 {Elwin Bruno}%
 {Christoffel}%
 {Элвин Бруно}%
 {Кристоффель}%
\RName
 {Elwin Bruno}%
 {Chritoffel}%
 {Эльвин Бруно}%
 {Кристоффель}%
\RName
 {Leopold}%
 {Kronecker}%
 {Леопольд}%
 {Кронекер}%
\RName
 {Martin}%
 {Kruskal}%
 {Мартин}%
 {Крускал}%
\RName
 {Charles-Augustin de}%
 {Coulomb}%
 {Шарль Огюстен де}%
 {Кулон}%
\RName
 {Leon Niels}%
 {Cooper}%
 {Леон Нил}%
 {Купер}%
\RName
 {Kazimierz}%
 {Kuratowski}%
 {Казимир}%
 {Куратовский}%
\RName
 {Martin Wilhelm}%
 {Kutta}%
 {Мартин Вильгельм}%
 {Кутта}%
\RName
 {Richard Christopher}%
 {Carrington}%
 {Ричард Кристофер}%
 {Кэррингтон}%
\RName
 {Lewis}%
 {Carroll}%
 {Льюис}%
 {Кэррол}%
\indexspace%
\RName
 {Guillaume Fran\c cois de}%
 {L'H\^opital}%
 {Гиймон Франсуа}%
 {де Лопиталь}%
\RName
 {Joseph Louis}%
 {Lagrange}%
 {Жозеф Луи}%
 {Лагранж}%
\RName
 {Tsit Yuen}%
 {Lam}%
 {Цит Юань}%
 {Лам}%
\RName
 {Lev}%
 {Landau}%
 {Лев}%
 {Ландау}%
\RName
 {Pierre-Simon}%
 {Laplace}%
 {Пьер-Симон}%
 {Лаплас}%
\RName
 {Henri}%
 {Lebesgue}%
 {Анри}%
 {Лебег}%
\RName
 {Beppo}%
 {Levi}%
 {Беппо}%
 {Леви}%
\RName
 {Leon Max}%
 {Lederman}%
 {Леон Макс}%
 {Ледерман}%
\RName
 {Adrien-Marie}%
 {Legendre}%
 {Адриен Мари}%
 {Лежандр}%
\RName
 {Gottfried Wilhelm}%
 {Leibniz}%
 {Готфрид Вильгельм}%
 {Лейбниц}%
\RName
 {Joel}%
 {Lebowitz}%
 {Джоэль}%
 {Лейбовиц}%
\RName
 {Serge}%
 {Lang}%
 {Серж}%
 {Ленг}%
\RName
 {Joseph}%
 {Lense}%
 {Джозеф}%
 {Лензе}%
\RName
 {Henrietta Swan}%
 {Leavitt}%
 {Генриетта Суон}%
 {Ливитт}%
\RName
 {Andrei}%
 {Linde}%
 {Андрей}%
 {Линде}%
\RName
 {Joseph}%
 {Liouville}%
 {Жозеф}%
 {Лиувилль}%
\RName
 {Andre}%
 {Lichnerowicz}%
 {Андрэ}%
 {Лихнерович}%
\RName
 {Nikolai Ivanovich}%
 {Lobachevsky}%
 {Николай Иванович}%
 {Лобачевский}%
\RName
 {Hendrik Antoon}%
 {Lorentz}%
 {Хендрик Антон}%
 {Лоренц}%
\RName
 {Edward Norton}%
 {Lorenz}%
 {Эдвард Нортон}%
 {Лоренц}%
\indexspace%
\RName
 {Ettore}%
 {Majorana}%
 {Этторе}%
 {Майорана}%
\RName
 {James Clerk}%
 {Maxwell}%
 {Джеймс Клерк}%
 {Максвелл}%
\RName
 {David Bryant}%
 {Mumford}%
 {Дэвид Брайант}%
 {Мамфорд}%
\RName
 {Ulf}%
 {Meissner}%
 {Ульф}%
 {Мейснер}%
\RName
 {Lise}%
 {Meitner}%
 {Лиза}%
 {Мейтнер}%
\RName
 {Eugene}%
 {Merzbacher}%
 {Юджин}%
 {Мерцбахер}%
\RName
 {Charles W.}%
 {Misner}%
 {Чарлз}%
 {Мизнер}%
\RName
 {Robert}%
 {Millikan}%
 {Роберт}%
 {Миликен}%
\RName
 {John Willard}%
 {Milnor}%
 {Джон Уиллард}%
 {Милнор}%
\RName
 {Hermann}%
 {Minkowski}%
 {Герман}%
 {Минковский}%
\RName
 {Eliakim Hastings}%
 {Moore}%
 {Элиаким Гастингс}%
 {Мур}%
\RName
 {Holger}%
 {Mueller}%
 {Холгер}%
 {Мюллер}%
\indexspace%
\RName
 {John von}%
 {Neyman}%
 {Джон фон}%
 {Нейман}%
\RName
 {Emmy}%
 {Noether}%
 {Эмми}%
 {Нётер}%
\RName
 {Katsumi}%
 {Nomizu}%
 {К.}%
 {Номидзу}%
\RName
 {}%
 {Newton}%
 {}%
 {Ньютон}%
\indexspace%
\RName
 {Yuri}%
 {Obukhov}%
 {Юрий}%
 {Обухов}%
\RName
 {Pierre Victor}%
 {Auger}%
 {Пьер Виктор}%
 {Оже}%
\RName
 {William of}%
 {Ockham}%
 {Уильям}%
 {Оккам}%
\RName
 {Peter}%
 {Olver}%
 {Питер}%
 {Олвер}%
\RName
 {Robert}%
 {Oppenheimer}%
 {Роберт}%
 {Опенгеймер}%
\indexspace%
\RName
 {Rudolf Ernst}%
 {Peierls}%
 {Рудольф Эрнст}%
 {Пайерлс}%
\RName
 {Blaise}%
 {Paskal}%
 {Блез}%
 {Паскаль}%
\RName
 {Wolfgang}%
 {Pauli}%
 {Вольфганг}%
 {Паули}%
\RName
 {Cecil Frank}%
 {Powell}%
 {Сесил Фрэнк}%
 {Пауэлл}%
\RName
 {Moritz}%
 {Pasch}%
 {Мориц}%
 {Паш}%
\RName
 {Giuseppe}%
 {Peano}%
 {Джузеппе}%
 {Пеано}%
\RName
 {Christian Heinrich Friedrich}%
 {Peters}%
 {Кристиан Генрих Фридрих}%
 {Петерс}%
\RName
 {James}%
 {Peebles}%
 {Джеймс}%
 {Пиблс}%
\RName
 {Max}%
 {Planck}%
 {Макс}%
 {Планк}%
\RName
 {Julius}%
 {Pl\"ucker}%
 {Юлий}%
 {Плюкер}%
\RName
 {Boris}%
 {Podolsky}%
 {Борис Яковлевич}%
 {Подольский}%
\RName
 {George}%
 {Polya}%
 {Джордж}%
 {Пойа}%
\RName
 {Isaac}%
 {Pomeranchuk}%
 {Исаак Яковлевич}%
 {Померанчук}%
\RName
 {Jean-Fran\c{c}ois}%
 {Pommaret}%
 {Жан-Франсуа}%
 {Поммаре}%
\RName
 {Lev Semenovich}%
 {Pontryagin}%
 {Лев Семёнович}%
 {Понтрягин}%
\RName
 {Jules Henri}%
 {Poincar\'e}%
 {Жюль Анри}%
 {Пуанкаре}%
\RName
 {Sim\'eon-Denis}%
 {Poisson}%
 {Симеон-Дени}%
 {Пуассон}%
\indexspace%
\RName
 {Isidor Isaac}%
 {Rabi}%
 {Исидор Айзек}%
 {Раби}%
\RName
 {Bertrand}%
 {Russell}%
 {Бертран}%
 {Рассел}%
\RName
 {Henry Norris}%
 {Russell}%
 {Генри Норрис}%
 {Расселл}%
\RName
 {Ernest}%
 {Rutherford}%
 {Эрнест}%
 {Резерфорд}%
\RName
 {Kurt Werner Friedrich}%
 {Reidemeister}%
 {Курт Вернер Фридрих}%
 {Рейдемейстер}%
\RName
 {William}%
 {Reynolds}%
 {Вильям}%
 {Рейнольдс}%
\RName
 {Ole Christensen}%
 {Roemer}%
 {Оле Кристенсен}%
 {Рёмер}%
\RName
 {Bernhard}%
 {Riemann}%
 {Бернхард}%
 {Риман}%
\RName
 {Martin}%
 {Rees}%
 {Мартин}%
 {Рис}%
\RName
 {Gregorio}%
 {Ricci}%
 {Грегорио}%
 {Риччи}%
\RName
 {Howard Percy}%
 {Robertson}%
 {Говард Перси}%
 {Робертсон}%
\RName
 {George Neil}%
 {Robertson}%
 {Джордж Нил}%
 {Робертсон}%
\RName
 {Carlo}%
 {Rovelli}%
 {Карло}%
 {Ровелли}%
\RName
 {Nathan}%
 {Rosen}%
 {Натан}%
 {Розен}%
\RName
 {Marshall Nicholas}%
 {Rosenbluth}%
 {Маршалл Николас}%
 {Розенблют}%
\RName
 {Michel}%
 {Rolle}%
 {Мишель}%
 {Ролль}%
\RName
 {Ryszard}%
 {R\d aczka}%
 {Р}%
 {Рончка}%
\RName
 {Karl David Tolme}%
 {Runge}%
 {Карл Давид Тольме}%
 {Рунге}%
\RName
 {Hanno}%
 {Rund}%
 {Ханно}%
 {Рунд}%
\RName
 {Warren de la}%
 {Rue}%
 {Уоррен де ла}%
 {Рю}%
\indexspace%
\RName
 {Anthony}%
 {Sudbery}%
 {Энтони}%
 {Садбери}%
\RName
 {Edwin Ernest}%
 {Salpeter}%
 {Эдвин Эрнест}%
 {Сальпетер}%
\RName
 {Leonard}%
 {Susskind}%
 {Леонард}%
 {Сасскинд}%
\RName
 {Corrado}%
 {Segre}%
 {Коррадо}%
 {Сегре}%
\RName
 {Waclaw}%
 {Sierpinski}%
 {Вацлав}%
 {Серпинский}%
\RName
 {Jean-Pierre}%
 {Serre}%
 {Жан Пьер}%
 {Серр}%
\RName
 {James Harris}%
 {Simons}%
 {Джеймс Харрис}%
 {Симонс}%
\RName
 {Tony}%
 {Skyrme}%
 {Тони}%
 {Скайрм}%
\RName
 {Dmitriy}%
 {Skobeltsin}%
 {Дмитрий}%
 {Скобельцин}%
\RName
 {Lee}%
 {Smolin}%
 {Ли}%
 {Смолин}%
\RName
 {Hartland Sweet}%
 {Snyder}%
 {Хартланд Свит}%
 {Снайдер}%
\RName
 {Eduard}%
 {Study}%
 {Эдуард}%
 {Стади}%
\RName
 {Shlomo}%
 {Sternberg}%
 {Соломон}%
 {Стернберг}%
\RName
 {Thomas Joannes}%
 {Stieltjes}%
 {Томас Иоаннес}%
 {Стилтьес}%
\RName
 {Norman}%
 {Steenrod}%
 {Норман}%
 {Стинрод}%
\RName
 {James}%
 {Stirling}%
 {Джеймс}%
 {Стирлинг}%
\RName
 {George Gabriel}%
 {Stokes}%
 {Джордж Габриель}%
 {Стокс}%
\RName
 {Leo}%
 {Szilard}%
 {Лео}%
 {Сцилард}%
\indexspace%
\RName
 {Alfred}%
 {Tarski}%
 {Альфред}%
 {Тарский}%
\RName
 {Walter }%
 {Thirring}%
 {Вальтер}%
 {Тирринг}%
\RName
 {Hans}%
 {Thirring}%
 {Ганс}%
 {Тирринг}%
\RName
 {Richard}%
 {Tolman}%
 {Ричард}%
 {Толмен}%
\RName
 {Ren\'e Fr\'ed\'eric}%
 {Thom}%
 {Рене Фредерик}%
 {Том}%
\RName
 {George Paget}%
 {Thomson}%
 {Джордж Паджет}%
 {Томсон}%
\RName
 {Kip S.}%
 {Thorne}%
 {Кип}%
 {Торн}%
\RName
 {Evangelista}%
 {Torricelli}%
 {Евангелиста}%
 {Торричелли}%
\indexspace%
\RName
 {Alfred North}%
 {Whitehead}%
 {Альфред Норт}%
 {Уайтхед}%
\RName
 {John Archibald}%
 {Wheeler}%
 {Джон Арчибалд}%
 {Уилер}%
\RName
 {Hassler}%
 {Whitney}%
 {Хасслер}%
 {Уитни}%
\RName
 {George Eugene}%
 {Uhlenbeck}%
 {Джордж Юджин}%
 {Уленбек}%
\RName
 {Harold Clayton}%
 {Urey}%
 {Гарольд Клейтон}%
 {Ури}%
\indexspace%
\RName
 {Michael}%
 {Faraday}%
 {Майкл}%
 {Фарадей}%
\RName
 {Richard}%
 {Feynman}%
 {Ричард}%
 {Фейнман}%
\RName
 {Carl}%
 {Faith}%
 {Карл}%
 {Фейс}%
\RName
 {Pierre}%
 {Fermat}%
 {Пьер}%
 {Ферма}%
\RName
 {Enrico}%
 {Fermi}%
 {Энрико}%
 {Ферми}%
\RName
 {Armand\Hyph Hippolyte\Hyph Louis}%
 {Fizeau}%
 {Арман Ипполит Луи}%
 {Физо}%
\RName
 {Serge\u{i} Pavlovich}%
 {Finikov}%
 {Сергей Павлович}%
 {Фиников}%
\RName
 {Paul}%
 {Finsler}%
 {Пауль}%
 {Финслер}%
\RName
 {Leopold}%
 {Vietoris}%
 {Леопольд}%
 {Фиторис}%
\RName
 {Grigorii Mikhailovich}%
 {Fikhtengolts}%
 {Григорий Михайлович}%
 {Фихтенгольц}%
\RName
 {John Ambrose}%
 {Fleming}%
 {Джон Амброз}%
 {Флеминг}%
\RName
 {Serge\u{i} Vasil'evich}%
 {Fomin}%
 {Сергей Васильевич}%
 {Фомин}%
\RName
 {Max}%
 {von Laue}%
 {Макс}%
 {фон Лауэ}%
\RName
 {James}%
 {Forbes}%
 {Джеймс}%
 {Форбс}%
\RName
 {George}%
 {Francis}%
 {Джордж}%
 {Франсис}%
\RName
 {Erik Ivar}%
 {Fredholm}%
 {Эрик Ивар}%
 {Фредгольм}%
\RName
 {Maurice Ren\'e}%
 {Fr\'echet}%
 {Морис Рене}%
 {Фреше}%
\RName
 {Alexander}%
 {Friedmann}%
 {Александр}%
 {Фридман}%
\RName
 {Ferdinand Georg}%
 {Frobenius}%
 {Фердинанд Георг}%
 {Фробениус}%
\RName
 {Guido}%
 {Fubini}%
 {Гвидо}%
 {Фубини}%
\RName
 {Joseph}%
 {Fourier}%
 {Джозеф}%
 {Фурье}%
\indexspace%
\RName
 {Edwin}%
 {Hubble}%
 {Эдвин}%
 {Хаббл}%
\RName
 {Felix}%
 {Hausdorff}%
 {Феликс}%
 {Хаусдорф}%
\RName
 {Friedrich W.}%
 {Hehl}%
 {Фридрих}%
 {Хель}%
\RName
 {Kurt}%
 {Hensel}%
 {Курт}%
 {Хензель}%
\RName
 {Friedrich Ernst Peter}%
 {Hirzebruch}%
 {Фридрих Эрнст Петер}%
 {Хирцебрух}%
\RName
 {Fred}%
 {Hoyle}%
 {Фред}%
 {Хойл}%
\RName
 {S. W.}%
 {Hawking}%
 {С. В.}%
 {Хокинг}%
\RName
 {Heinz}%
 {Hopf}%
 {Хайнц}%
 {Хопф}%
\RName
 {Antony}%
 {Hewish}%
 {Энтони}%
 {Хьюиш}%
\indexspace%
\RName
 {Ernst Friedrich Ferdinand}%
 {Zermelo}%
 {Эрнст Фридрих Фердинанд}%
 {Цермело}%
\indexspace%
\RName
 {Nikolai Grigorievich}%
 {Chebotaryov}%
 {Николай Григорьевич}%
 {Чеботарёв}%
\RName
 {Pafnuty Lvovich}%
 {Chebyshev}%
 {Пафнутий Львович}%
 {Чебышев}%
\RName
 {James}%
 {Chadwick}%
 {Джеймс}%
 {Чедвик}%
\RName
 {Owen}%
 {Chamberlain}%
 {Оуэн}%
 {Чемберлен}%
\RName
 {Shiing-Shen}%
 {Chern}%
 {Шиинг-Шен}%
 {Черн}%
\RName
 {N. A.}%
 {Chernikov}%
 {Н. А.}%
 {Черников}%
\RName
 {Geoffrey Foucar}%
 {Chew}%
 {Джеффри Фаукар}%
 {Чу}%
\indexspace%
\RName
 {Jacob}%
 {Schwartz}%
 {Дж.}%
 {Шварц}%
\RName
 {Laurent Moise}%
 {Schwartz}%
 {Лоран Моиз}%
 {Шварц}%
\RName
 {Karl}%
 {Schwarzschild}%
 {Карл}%
 {Шварцшильд}%
\RName
 {Julian}%
 {Schwinger}%
 {Джулиус}%
 {Швингер}%
\RName
 {Claude}%
 {Chevalley}%
 {Клод}%
 {Шевалле}%
\RName
 {}%
 {Schmidt}%
 {}%
 {Шмидт}%
\RName
 {Erwin}%
 {Schr$\ddot{o}$dinger}%
 {Эрвин}%
 {Шредингер}%
\RName
 {John}%
 {Schrieffer}%
 {Джон}%
 {Шриффер}%
\RName
 {Paul}%
 {Steinhardt}%
 {Поль}%
 {Штейнхардт}%
\RName
 {Eduard}%
 {Stiefel}%
 {Эдуард}%
 {Штифель}%
\RName
 {Jacques Charles Fracois}%
 {Sturm}%
 {Жак Шарль Франсуа}%
 {Штурм}%
\indexspace%
\RName
 {}%
 {Euclid}%
 {}%
 {Эвклид}%
\RName
 {Arthur Stanley}%
 {Eddington}%
 {Артур Стэнли}%
 {Эддингтон}%
\RName
 {Luther Pfahler}%
 {Eisenhart}%
 {Л. П.}%
 {Эйзенхарт}%
\RName
 {Samuel}%
 {Eilenberg}%
 {Самуил}%
 {Эйленберг}%
\RName
 {Leonhard}%
 {Euler}%
 {Леонард}%
 {Эйлер}%
\RName
 {Albert}%
 {Einstein}%
 {Альберт}%
 {Эйнштейн}%
\RName
 {Lev}%
 {Elsgolts}%
 {Л. Э.}%
 {Эльсгольц}%
\RName
 {Federigo}%
 {Enriques}%
 {Федериго}%
 {Энрикес}%
\RName
 {Paul Sophus}%
 {Epstein}%
 {Пауль Софус}%
 {Эпштейн}%
\RName
 {Werner}%
 {Ehrenberg}%
 {Вернер}%
 {Эренберг}%
\RName
 {Paul}%
 {Ehrenfest}%
 {Пауль}%
 {Эренфест}%
\indexspace%
\RName
 {Hideki}%
 {Yukawa}%
 {Хидеки}%
 {Юкава}%
\indexspace%
\RName
 {Carl}%
 {Jacobi}%
 {Карл}%
 {Якоби}%
\RName
 {Chen Ning}%
 {Yang}%
 {Чжэньнин}%
 {Янг}%
\RName
 {Shing-Tung}%
 {Yau}%
 {Шин-Тан}%
 {Яу}%

\end{theindex}%

\def\texTitle{ English Russian Scientific Dictionary }
\def\texAuthors{ Aleks Kleyn }
\def\texComments{ 67 pages. }
\def\texAbstract
{
English Russian and Russian English dictionaries presented
in this paper are dedicated to help translate a text
from one language to another. I also included the bilingual
name index into this book.
}



\end{document}